 	\def\version{\today}
\newif\ifpdf
\ifx\pdfoutput\undefined
\pdffalse % we are not running PDFLaTeX
\else
\pdfoutput=1 % we are running PDFLaTeX
\pdftrue \fi

\newif\iffinal
\finalfalse	% Not final version
\finaltrue % Final version

\documentclass[reqno,twoside,11pt]{amsart}
\iffinal\else\usepackage[notref,notcite]{showkeys}\fi
\iffinal\else\IfFileExists{pdfsync.sty}{\usepackage{pdfsync}}{}\fi
\usepackage{cite}
\usepackage{amsmath}
\usepackage{amsfonts}
\usepackage{amssymb}
\usepackage{epsfig}
\usepackage{verbatim}

\usepackage{mathrsfs}
\usepackage{mathpazo}
\usepackage{mathabx}

\usepackage{amsthm,mathtools}
\usepackage{xcolor}
\usepackage{esint}

\DeclareFontFamily{OT1}{eusb}{} \DeclareFontShape{OT1}{eusb}{m}{n} {<5> <6> <7> <8> <9> <10> <11> <12> <14.4> eusb10}{}
\DeclareMathAlphabet{\eusb}{OT1}{eusb}{m}{n}

\DeclareFontFamily{OT1}{eusm}{} \DeclareFontShape{OT1}{eusm}{m}{n} {<5> <6> <7> <8> <9> <10> <11> <12> <14.4> eusm10}{}
\DeclareMathAlphabet{\eusm}{OT1}{eusm}{m}{n}

\DeclareFontFamily{OT1}{eufm}{} \DeclareFontShape{OT1}{eufm}{m}{n} {<5> <6> <7> <8> <9> <10> <11> <12> <14.4> eufm10}{}
\DeclareMathAlphabet{\mathfrak}{OT1}{eufm}{m}{n}

\DeclareFontFamily{OT1}{fraktura}{}
\DeclareFontShape{OT1}{fraktura}{m}{n} {<5> <6> <7> <8> <9> <10> <11> <12> <13> <14.4> [1.1] eufm10}{}
\DeclareMathAlphabet{\fraktura}{OT1}{fraktura}{m}{n}

\DeclareFontFamily{OT1}{cmfi}{} \DeclareFontShape{OT1}{cmfi}{m}{n} {<5> <6> <7> <8> <9> <10> <11> <12> <13> <14.4> [0.9] cmfi10}{}
\DeclareMathAlphabet{\cmfi}{OT1}{cmfi}{b}{n}

\DeclareFontFamily{OT1}{cmss}{} \DeclareFontShape{OT1}{cmss}{m}{n} {<5> <6> <7> <8> <9> <10> <11> <12> <13> <14.4> cmss10}{}
\DeclareMathAlphabet{\cmss}{OT1}{cmss}{m}{n}
\newtheoremstyle{thm}{1.5ex}{1.5ex}{\itshape\rmfamily}{} {\bfseries\rmfamily}{}{2ex}{}

\newtheoremstyle{def}{1.5ex}{1.5ex}{\rmfamily\sl}{} {\bfseries\rmfamily}{}{2ex}{}

\newtheoremstyle{rem}{1.3ex}{1.3ex}{\rmfamily}{} {\itshape}
{} {1.5ex}{}

\newenvironment{proofsect}[1] {\vskip0.1cm\noindent{\rmfamily\itshape#1.}}{\qed\vspace{0.15cm}}%{\newline\vspace{0.15cm}}

\theoremstyle{thm}
\newtheorem{theorem}{Theorem}[section]
\newtheorem{lemma}[theorem]{Lemma}
\newtheorem{proposition}[theorem]{Proposition}

\newtheorem*{Main Theorem}{Main Theorem.}
\newtheorem{corollary}[theorem]{Corollary}

\theoremstyle{def}
\newtheorem{definition}[theorem]{Definition}

\theoremstyle{rem}

\numberwithin{equation}{section}

\renewcommand{\section}{\secdef\sct\sect}
\newcommand{\sct}[2][default]{\refstepcounter{section}
\addcontentsline{toc}{section}
{{\tocsection {}{\thesection}{\!\!\!\!#1\dotfill}}{}}
\vspace{0.7cm}
\centerline{ %\large
\scshape\arabic{section}.\ #1} \nopagebreak \vspace{0.2cm}}
\newcommand{\sect}[1]{
\vspace{0.4cm} \centerline{\large\scshape\rmfamily #1}
\vspace{0.2cm}}

\renewcommand{\subsection}{\secdef\subsct\sbsect}
\newcommand{\subsct}[2][default]{\refstepcounter{subsection}
\addcontentsline{toc}{subsection}
{{\tocsection{\!\!}{\hspace{1.2em}\thesubsection}{\!\!\!\!#1\dotfill}}{}}
\vspace{0.45\baselineskip} {\flushleft\bf
\arabic{section}.\arabic{subsection}~\bf #1.~}
\\*[3mm]\noindent
\nopagebreak}
\newcommand{\sbsect}[1]{\vspace{0.1cm}\noindent
\textbf{#1.~}\vspace{0.1cm}}

\renewcommand{\subsubsection}{%
\secdef \subsubsect\sbsbsect}
\newcommand{\subsubsect}[2][default]{%
\refstepcounter{subsubsection} 
\addcontentsline{toc}{subsubsection}{{\tocsection{\!\!}
{\hspace{3.05em}\thesubsubsection}{\!\!\!\!#1\dotfill}}{}}
\nopagebreak
\vspace{0.15\baselineskip} \nopagebreak {\flushleft\rmfamily
\itshape\arabic{section}.\arabic{subsection}.\arabic{subsubsection}
\ \rmfamily #1\/.}\ }
\newcommand{\sbsbsect}[1]{\vspace{0.1cm}\noindent
\rmfamily \itshape
\arabic{section}.\arabic{subsection}.\arabic{subsubsection} \
\sffamily #1\/.\ }

\iffinal

\else

\fi

\renewcommand{\caption}[1]{%
\vglue0.5cm
\refstepcounter{figure}
\begin{minipage}{0.9\textwidth}\small {\sc Figure~\thefigure. }#1\end{minipage}}

\newcommand{\supp}{\operatorname{supp}}

\newcommand{\textd}{\text{\rm d}\mkern0.5mu}

\newcommand{\texte}{\text{\rm  e}\mkern0.7mu}

\newcommand{\1}{{1\mkern-4.5mu\textrm{l}}}
\renewcommand{\1}{\text{\sf 1}}

\renewcommand{\AA}{\mathcal A}
\newcommand{\BB}{\mathcal B}

\newcommand{\DD}{\mathcal D}
\newcommand{\EE}{\mathcal E}

\newcommand{\GG}{\mathcal G}
\newcommand{\HH}{\mathcal H}
\newcommand{\II}{\mathcal I}
\newcommand{\JJ}{\mathcal J}

\newcommand{\LL}{\mathcal L}
\newcommand{\MM}{\mathcal M}
\newcommand{\NN}{\mathcal N}

\newcommand{\CalS}{\mathcal S}

\newcommand{\ZZ}{\mathcal Z}

\newcommand{\E}{\mathbb E}

\newcommand{\BbbL}{\mathbb L}

\newcommand{\N}{\mathbb N}

\newcommand{\BbbP}{\mathbb P}

\newcommand{\R}{\mathbb R}

\newcommand{\T}{\mathbb T}

\newcommand{\V}{\mathbb V}

\newcommand{\Z}{\mathbb Z}

\newcommand{\scrF}{\mathscr{F}}

\newcommand{\twoeqref}[2]{(\ref{#1}--\ref{#2})}

\def\myffrac#1#2 in #3{\raise 2.6pt\hbox{$#3 #1$}\mkern-1.5mu\raise 0.8pt\hbox{$#3/$}\mkern-1.1mu\lower 1.5pt\hbox{$#3 #2$}}
\newcommand{\ffrac}[2]{\mathchoice%
	{\myffrac{#1}{#2} in \scriptstyle}
	{\myffrac{#1}{#2} in \scriptstyle}
	{\myffrac{#1}{#2} in \scriptscriptstyle}
	{\myffrac{#1}{#2} in \scriptscriptstyle}
}

\newcommand{\cc}{\text{\rm c}}

\newcommand{\wt}{\widetilde}
\newcommand{\wh}{\widehat}
\newcommand{\Cloc}{\fraktura C_{\text{\rm loc}}}
\newcommand{\frakp}{\mathfrak p}
\newcommand{\frakm}{\mathfrak m}

\newcommand{\laweq}{\,\overset{\text{\rm law}}=\,}
\newcommand{\lawarrow}{{\overset{\text{\rm law}}\longrightarrow}}
\newcommand{\slb}{\sqrt{\log b}}
\newcommand{\PPP}{\text{\rm PPP}}
\newcommand{\ZZZ}{\eusb Z}

\newcommand{\fraka}{\fraktura a}

\newcommand{\zero}{{\underline{0}}}

\newcommand\independent{\protect\mathpalette{\protect\independenT}{\perp}}
\def\independenT#1#2{\mathrel{\rlap{$#1#2$}\mkern3mu{#1#2}}}

\begin{document}

\vglue-2mm

\title[Local-time extremal process\hfill]{% \hfill\version\hfill]{
Extremal process of the local time of simple\\random walk on a regular tree}
%\author[\hfill\version\hfill Y.~Abe, M.~Biskup]
\author[\hfill Y.~Abe, M.~Biskup]
{Yoshihiro Abe$^1$ \,{\tiny and}\, Marek Biskup$^2$}
\thanks{\hglue-4.5mm\fontsize{9.6}{9.6}\selectfont\copyright\,\textrm{2025}\ \ \textrm{Y.~Abe, M.~Biskup. Reproduction for non-commercial use is permitted.\\Typeset on \version}}
\maketitle

\vglue-5mm
\centerline{\it $^1$Mathematical Institute, Tohoku University, Sendai, Japan}
\centerline{\it $^2$Department of Mathematics, UCLA, Los  Angeles, California, USA}
%\smallskip

%\smallskip
%\centerline{\version}

\begin{abstract}
We study a continuous-time simple random walk on a regular rooted tree of depth~$n$ in two  settings: either the walk is started from a leaf vertex and run until the tree root is first hit or it is started from the root and run until it has spent a prescribed amount of time  there. In both cases we show that the extremal process associated with centered square-root local time on the leaves tends, as $n\to\infty$, to a decorated Poisson point process with a random intensity measure. 
While the intensity measure is specific to the local-time problem at hand, the decorations are exactly those for the tree-indexed Markov chain (a.k.a.\ Branching Random Walk or Gaussian Free Field) with normal step distribution. The proof demonstrates the latter by way of a Lindeberg-type swap of the decorations of the two processes which itself relies on a well-known isomorphism theorem.
\end{abstract}

%

%\tableofcontents\vglue-1cm

\section{Introduction and results}
\label{sec1}\noindent
%\vglue-2mm
%\subsection{Background}
Consider a regular rooted tree~$\T_n$ of depth~$n\ge1$ and forward degree $b\ge2$. Pick any of its~$b^n$ leaf vertices and
use it to start a continuous-time Markov chain with state space~$\T_n$ and unit jump rate across each edge of~$\T_n$. Write $\ell_t(x)$ for the time spent by the chain at~$x$ up to time~$t$ and~$\tau_\varrho$ for the first time the chain hits the root~$\varrho$ of~$\T_n$. We are interested in the extremal properties of~$\ell_{\tau_\varrho}$ on the set~$\BbbL_n$ of the leaves~of~$\T_n$.

The limit distribution of the maximal time spent at any leaf-vertex has been identified in a recent study by the second author and O.~Louidor~\cite{BL4} drawing on earlier work of the first author~\cite{A18}. One way to state the conclusion is
\begin{equation}
\label{E:1.1}
\frac1n\biggl(\max_{x\in\BbbL_n}\ell_{\tau_\varrho}(x) -\bigl(n^2\log b -2n\log n\bigr)\biggr) 
\,\,\underset{n\to\infty}\lawarrow\,\,\log\ZZZ+G,
\end{equation}
where~$\ZZZ$ is an a.s.-positive random variable whose law can be characterized and~$G$ is a normalized Gumbel random variable independent of~$\ZZZ$. Perhaps more familiar, albeit equivalent, way to put this is by saying that, for all~$u\in\R$ and $x_n\in\BbbL_n$,
\begin{equation}
\label{E:1.2}
P^{x_n}\biggl(\max_{x\in\BbbL_n}\sqrt{\ell_{\tau_\varrho}(x)}\le \slb\, n-\frac1{\slb}\log n +u\biggr)
\,\,\underset{n\to\infty}\longrightarrow\,\,\E\bigl(\texte^{-\ZZZ \,\texte^{-2u\sqrt{\log b}}}\bigr),
\end{equation}
where~$P^{x_n}$ is the law of the chain started at~$x_n$ and the expectation is with respect to the law of~$\ZZZ$. The structure of the limit law places this problem in the universality class of log-correlated models; see Section~\ref{sec-2} for more discussion and references.

With the maximal time spent by the chain identified, a natural list of follow-up questions arises; for instance: How is the maximizer distributed relative to the starting point of the chain? What is the law of the second, third, etc maximum? Is there spatial clustering? What is the dependency structure of individual clusters? As is standard in extreme-order statistics, such questions are conveniently encoded by an associated empirical extremal process. The present paper aims to extract a weak limit of this process as the depth of the tree tends to infinity. We will treat two settings: the chain started from a leaf and the chain started from the root.

\subsection{Random walk started from a leaf}
We start with some notation. Observe that each leaf-vertex $x\in\BbbL_n$ can be identified with a sequence $(x_1,\dots,x_n)\in\{0,\dots,b-1\}^n$ of ``instructions'' indicating the ``turns'' in the unique path in~$\T_n$ from the root to~$x$. Relying on this representation, we define an injection~$\theta_n\colon\BbbL_n\to[0,1]$ by
\begin{equation}
%\label{}
\theta_n(x):=\sum_{i=1}^n b^{-i}x_i\quad\text{ when }\quad x=(x_1,\dots,x_n).
\end{equation}
We will write~$\zero$ for the vertex represented by the sequence~$(0,\dots,0)$. 

We are interested in simple random walk on~$\T_n$ that, technically, is a continuous-time Markov chain $\{X_t\colon t\ge0\}$ on~$\T_n$ with infinitesimal generator acting on functions $f\colon\T_n\to\R$ as $\LL f(x):=\sum_{(x,y)\in E(\T_n)}[f(y)-f(x)]$, where~$E(\T_n)$ is the set of undirected edges of~$\T_n$. Let $P^x$ denote for the law of the walk started from~$x\in\T_n$. For~$x\in\T_n$ and~$t\ge0$, let~$\ell_t(x):=\int_0^t 1_{\{X_s=x\}}\textd s$ be the total time spent by~$X$ at~$x$ by time~$t$. Write 
\begin{equation}
\label{E:1.4}
m_n:=\slb\, n-\frac1{\slb}\log n
\end{equation}
for the centering sequence from \eqref{E:1.2} and denote~$a^+:=\max\{a,0\}$. We then have:

\begin{theorem}
\label{thm-1}
There exists a random Borel measure~$\ZZ$ on~$[0,1]$ and a (deterministic) law~$\DD$ on infinite, locally finite point processes on~$(-\infty,0]$ such that, for~$\ell_{\tau_\varrho}$  sampled under~$P^\zero$,
\begin{equation}
\label{E:1.5}
\sum_{x\in\BbbL_n}\delta_{\theta_n(x)}\otimes\delta_{\sqrt{\ell_{\tau_\varrho}(x)}-m_n}\,\,\,\underset{n\to\infty}\lawarrow\,\,\,\sum_{i\ge1}\sum_{j\ge1}\delta_{x_i}\otimes\delta_{h_i+d_j^{(i)}},
\end{equation}
where~$\{(x_i,h_i)\}_{i\ge1}$ enumerates points in a sample from the Poisson point process
\begin{equation}
\label{E:1.6}
\PPP\bigl(\,\ZZ(\textd x)\otimes\texte^{-2h\slb}\textd h\bigr)
\end{equation}
and~$\{d_j^{(i)}\colon j\ge1\}_{i\ge1}$ are i.i.d.\ samples from~$\DD$ independent of ~$\{(x_i,h_i)\}_{i\ge1}$ and~$\ZZ$. 
Moreover, a.e.\ sample of $\ZZ$ is such that $\ZZ([0,1]) <\infty$, $\ZZ([0,\epsilon))>0$ for each~$\epsilon>0$ yet $\ZZ(\{0\})=0$ and there exists a constant $\wt C_\star\in(0,\infty)$ such that, for~$\ell_{\tau_\varrho}$ sampled from~$P^\zero$,
\begin{equation}
\label{E:1.7}
\wt C_\star\, b^{-2n} \sum_{x \in \BbbL_n} \Bigl(n\slb - 
\sqrt{\ell_{\tau_\varrho}(x)}\,\Bigr)^+ \ell_{\tau_\varrho}(x)^{1/4}\,
	\texte^{2\slb\,\sqrt{\ell_{\tau_\varrho}(x)}}\delta_{\theta_n(x)}\,\,\,\underset{n\to\infty}\lawarrow\,\,\,\ZZ.
\end{equation}
A.e.~sample of $\{d_j\colon j\ge1\}$ from~$\DD$ has a point at the origin.
\end{theorem}

The convergence in law of the point measures in \eqref{E:1.5} is relative to the vague topology on Radon measures on~$[0,1]\times\R$. The sampling in \eqref{E:1.6} is done conditional on~$\ZZ$ (i.e., $\ZZ$ is sampled first and the Poisson points second). An equivalent way to state \eqref{E:1.5} is by saying that, for all continuous~$f\colon[0,1]\times\R\to[0,\infty)$ with compact support, 
\begin{equation}
\label{E:1.8}
\begin{aligned}
E^\zero\biggl(\exp&\Bigl\{-\sum_{x\in\BbbL_n}f\bigl(\theta_n(x),\sqrt{\ell_{\tau_\varrho}(x)}-m_n\bigr)\Bigr\}\biggr)
\\
&\underset{n\to\infty}\longrightarrow\,\,\E\biggl(\exp\Bigl\{-\int\ZZ(\textd x)\otimes\texte^{-2h\slb}\textd h\otimes\DD(\textd\xi)\bigl(1-\texte^{-\langle\xi, f(x,h+\cdot)\rangle}\bigr)\Bigr\}\biggr),
\end{aligned}
\end{equation}
where  $E^\zero$ denotes the expectation with respect to~$P^\zero$ and~$\E$ is the expectation with respect to the law of~$\ZZ$. The expression~$\langle\xi, f(x,h+\cdot)\rangle$ abbreviates the integral of $s\mapsto f(x,h+s)$ with respect to the point measure~$\xi$.

The limit law \eqref{E:1.5} has the structure of a randomly-shifted, decorated Poisson point process. 
Indeed, writing the intensity in \eqref{E:1.6} as
\begin{equation}
%\label{}
\wh\ZZ(\textd x)\otimes\Bigl(\ZZ\bigl([0,1]\bigr)\texte^{-2h\slb}\textd h\Bigr),
\end{equation}
 where~$\wh\ZZ$ is~$\ZZ$ normalized by its total mass~$\ZZ([0,1])$, we can realize the objects on the right of \eqref{E:1.5} as follows: First draw $\ZZ$ and then draw i.i.d.\ samples~$\{x_i\}_{i\ge1}$ from~$\wh\ZZ$. Then, given an independent sample~$\{h_i'\}_{i\ge1}$ from the Poisson point process of Gumbel intensity $\texte^{-2u\slb}\textd u$, let~$\{h_i\}_{i\ge1}$ be defined by
\begin{equation}
%\label{}
h_i:=h_i'+(2\slb)^{-1}\log\ZZ([0,1]),\quad i\ge1.
\end{equation}
Finally, ``attach'' to each~$h_i$ an independent sample (a ``decoration'') from~$\DD$ drawn independently of the points $\{(x_i,h_i)\}_{i\ge1}$. 
 
The stated properties of the law~$\DD$ ensure that the ``cluster'' of points ``attached'' to~$h_i$ has its maximal point at~$h_i$. It follows that, for the process on the right of \eqref{E:1.5} to not charge~$[0,1]\times(u,\infty)$, we need to have~$h_i\le u$ for all~$i\ge1$. Taking~$f$ along a sequence of approximations of $1_{[0,1]\times(u,\infty)}$ in \eqref{E:1.8} then forces~$\ZZZ$ from \eqref{E:1.2} to obey
\begin{equation}
\label{E:1.9}
\ZZZ\laweq \frac1{2\slb}\,\ZZ\bigl([0,1]\bigr).
\end{equation}
The equality in law is confirmed independently from the fact that the total mass of the measure on the left of \eqref{E:1.7} is known to converge weakly to~$\ZZZ$; see~\cite[Theorem~1.5]{BL4}. (The change in normalization due to the prefactor in \eqref{E:1.9} accounts for the difference between~$\wt C_\star$ above and~$C_\star$ in~\cite[Theorem~1.5]{BL4}.)

\subsection{Random walk started from the root}
\label{sec-1.2}\noindent
The setting of the random walk started from a leaf and killed upon first visit to the root has been introduced in order to mimic the exit problem from a lattice domain; see the discussion in Section~\ref{sec-2}. For the random walk on~$\T_n$, another possible setting of interest is that of the walk started at the root~$\varrho$. This becomes particularly neat if we parametrize the process by the time spent at~$\varrho$. To this end we set, for each~$t\ge0$, 
\begin{equation}
\label{E:1.10}
\wt\tau_\varrho(t):=\inf\bigl\{s\ge0\colon\ell_s(\varrho)\ge t\bigr\}
\end{equation}
and abbreviate
\begin{equation}
\label{E:1.11}
L_t(x):=\ell_{\wt\tau_\varrho(t)}(x),\quad x\in\T_n.
\end{equation}
A key technical advantage of this representation is that $\{L_t(x)\colon x\in\T_n\}$ is a time-homo\-gene\-ous tree-indexed Markov chain; see Lemma~\ref{lemma-M}. Note that~$L_t(\varrho)=t$ deterministically while $E^\varrho(L_t(x))=t$ for all~$x\in\T_n$. 

The present setting does not necessitate that the walk ever visits the leaves by the time it has accumulated time~$t$ at the root and we in fact have~$L_t(x)=0$ for all~$x\in\BbbL_n$ with uniformly positive probability. For this reason we introduce
\begin{equation}
\label{E:1.14b}
\tau_{\BbbL_n}:=\inf\Bigl\{t\ge0\colon\max_{x\in\BbbL_n}L_t(x)>0\Bigr\}
\end{equation}
which a.s.\ coincides with the first time~$X$ visits~$\BbbL_n$.
Recall the notation~$m_n$ for the sequence from \eqref{E:1.4}. We then claim: 

\begin{theorem}
\label{thm-2}
For all $t > 0$, there exists an a.s.-finite random Borel measure~$Z_t$ on~$[0,1]$ such that the following holds:
\begin{equation}
\label{E:1.13}
\lim_{n\to\infty}P^\varrho(\tau_{\BbbL_n}<t) = \BbbP\bigl(Z_t([0,1])>0\bigr)\in(0,1)
\end{equation}
and, for~$L_t$ sampled under the conditional law $P^\varrho(\cdot\,|\tau_{\BbbL_n}<t)$, 
\begin{equation}
\label{E:1.14}
\sum_{x\in\BbbL_n}\delta_{\theta_n(x)}\otimes\delta_{\sqrt{L_t(x)}-m_n}\,\,\,\underset{n\to\infty}\lawarrow\,\,\,\sum_{i\ge1}\sum_{j\ge1}\delta_{x_i}\otimes\delta_{h_i+d_j^{(i)}}.
\end{equation}
Here~$\{(x_i,h_i)\}_{i\ge1}$ enumerates points in a sample of the Poisson point process
\begin{equation}
\label{E:1.15}
\PPP\bigl(\,Z_t(\textd x)\otimes\texte^{-2h\slb}\textd h\bigr)
\end{equation}
with~$Z_t$ sampled conditional on~$\{Z_t([0,1])>0\}$
and~$\{d_j^{(i)}\colon j\ge1\}_{i\ge1}$ denoting i.i.d.\ samples from~$\DD$ in Theorem~\ref{thm-1}, drawn independently of ~$\{(x_i,h_i)\}_{i\ge1}$ and~$Z_t$. 
Moreover, with $\wt C_\star \in (0,\infty)$ as in Theorem~\ref{thm-1} and~$L_t$ sampled under~$P^\varrho$,
\begin{equation}
\label{E:1.18}
\wt C_\star\,b^{-2n} \sum_{x \in \BbbL_n} \Bigl(n\slb - 
\sqrt{L_t(x)}\,\Bigr)^+ L_t(x)^{1/4}\,
	\texte^{2\slb\,\sqrt{L_t(x)}}\delta_{\theta_n(x)}\,\,\underset{n\to\infty}\lawarrow\,\,Z_t.
\end{equation}
A.e.\ sample of~$Z_t$ is diffuse (i.e., does not charge singletons). 
\end{theorem}

The form and structure of the limit law \eqref{E:1.14} is quite similar to that in \eqref{E:1.5} and so the discussion after Theorem~\ref{thm-1} applies here as well. In particular, the maximum of~$L_t$ obeys the analogue of \eqref{E:1.1},  where the random variable~$\ZZZ$ is replaced by $(2\slb)^{-1}Z_t([0,1])$. As shown in~\cite{BL4}, we have $P(Z_t([0,1])>0)\to 0$ as~$t\downarrow0$. We will set~$Z_0(\cdot):=0$. 

\subsection{Connecting the intensity measures}
Theorem~\ref{thm-1} will be derived from Theorem~\ref{thm-2} and the proof expresses the law of~$\ZZ$ via those of the random measures $\{Z_t\colon t\ge0\}$. To state this precisely, note that the weak convergence \eqref{E:1.18} along with the continuity of the left-hand side in~$t$ implies that, for any continuous~$f\colon[0,1]\to[0,\infty)$,
\begin{equation}
%\label{}
\phi_f(t):=\E\bigl(\texte^{-\langle Z_t,f\rangle}\bigr)
\end{equation}
defines a Borel measurable function~$t\mapsto\phi_f(t)$. This is enough to couple any random variable~$T\ge0$ with random Borel measure~$Z_T$ so that
\begin{equation}
\label{E:1.20z}
\E\bigl(\texte^{-\lambda T-\langle Z_T,f\rangle}\bigr) = \E \bigl(\texte^{-\lambda T}\phi_f(T)\bigr)
\end{equation}
holds for all continuous~$f\colon[0,1]\to[0,\infty)$ and all~$\lambda\ge0$.
Under the joint law~$(T,Z_T)$, the conditional law of~$Z_T$ given~$T$ is that of~$Z_t$ with~$t:=T$. (If we had control of the law of the process~$t\mapsto Z_t$, we would simply sample~$T$ independently and plug it for~$t$.)

A similar coupling of random Borel measures~$\{Z_{T_k}\}_{k\ge0}$ exists for any sequence of non-negative random ``times'' $\{T_k\}_{k\ge0}$, independent or not, so that
\begin{equation}
\label{E:1.21a}
\E\biggl(\,\prod_{k\ge0}\texte^{-\lambda_k T_k-\langle Z_{T_k},f_k\rangle}\biggr)=\E\biggl(\,\prod_{k\ge1}\bigl[\texte^{-\lambda_k T_k}\phi_{f_k}(T_k)\bigr]\biggr)
\end{equation}
for all continuous $f_k\colon[0,1]\to[0,\infty)$ and all~$\lambda_k\ge0$. Here, conditionally on~$\{T_k\}_{k\ge0}$, the random variables $\{Z_{T_k}\}_{k\ge0}$ are independent with $Z_{T_k}$ having the law defined in \eqref{E:1.20z}.
In the proof of Theorem~\ref{thm-1} we characterize the law of~$\ZZ$ as follows:

\begin{corollary}
\label{cor-1.3}
Let~$\{T_k\}_{k\ge0}$ have the law of $\{\frac12|B_k|^2\}_{k\ge0}$ for~$B$ denoting a standard two-dimen\-sional Brownian motion. Let~$\{Z_{T_k}\}_{k\ge0}$  be a sequence of random Borel measures on~$[0,1]$ coupled with~$\{T_k\}_{k\ge0}$ so that \eqref{E:1.21a} holds as stated. Then
\begin{equation}
\label{E:1.23a}
\ZZ(\textd x)\laweq \sum_{k\ge1}1_{ [b^{-k-1},b^{-k})}(x)\, b^{-2k} Z_{T_k}\bigl(b^k\textd x\bigr). 
\end{equation}
A.e.\ sample of~$\ZZ$ is diffuse.
\end{corollary}

Note that the measure in \eqref{E:1.23a} is concentrated on~$[0,b^{-1})$. In light of \eqref{E:1.7}, this is no surprise as, in order for the walk started at~$\zero$ to reach the parts of $\T_n$ where~$\theta_n(x)\ge b^{-1}$, it must hit the root of~$\T_n$ first. The representation also yields a limit law for the location of the maximizer of~$\ell_{\tau_\varrho}$:

\begin{corollary}
\label{cor-1.4}
Given~$n\ge1$ and a sample~$\ell_{\tau_\varrho}$ on~$\T_n$ under~$P^\zero$, let~$Y_n$ denote the a.s.-unique maximizer of~$x\mapsto\ell_{\tau_\varrho}(x)$ on~$\BbbL_n$. Then~$\theta_n(Y_n)$ converges in law to a random variable~$U$ on~$[0,1]$ described as follows: Given a sample~$\{Z_{T_k}\}_{k\ge1}$ of the measures in Corollary~\ref{cor-1.3}, let~$K$ be the smallest~$k\ge1$ maximizing
\begin{equation}
%\label{}
k\mapsto \log Z_{T_k}\bigl([b^{-1},1]\bigr)-2k\log b+G_k
\end{equation}
where~$\log0:=-\infty$ and $\{G_k\}_{k\ge1}$ are i.i.d.\ standard Gumbel 
independent of $\{Z_{T_k}\colon k\ge1\}$.
Then $K<\infty$ a.s.\ and $U$ is uniform on~$[b^{-k-1},b^{-k})$ conditional on~$\{K=k\}$.
\end{corollary}

%\begin{remark}
We expect the conclusion of Theorem~\ref{thm-2} to hold even for~$t$ that increases with~$n$ (albeit slower than~$n^2$) provided that the centering is done by the sequence $\sqrt t+a_n(t)$ with
\begin{equation}
\label{E:1.24a}
a_n(t):=n\slb -\frac34\frac1{\slb}\log n -\frac1{4\slb}\log\Bigl(\frac{n+\sqrt t}{\sqrt t}\Bigr)
\end{equation}
instead of~$m_n$ in \eqref{E:1.4}.
 This would require performing the relevant estimates uniformly in~$t\in(0,t_n]$, for any sequence~$t_n=o(n^2)$, similarly as was done for the maximal local time in~\cite[Theorem~3.1]{BL4}. As explained in~\cite[Remark~3.8]{BL4}, 
for~$t$ growing as, or faster than, order~$n^2$, additional corrections arise; see~\cite[Theorem 1.1]{A18}. (Roughly speaking, the ``constant'' $\wt C_\star$ becomes dependent on the asymptotic value of~$t/n^2$.) 
 
We also expect that, under~$P^\varrho$,
\begin{equation}
\label{E:1.25}
t^{-1/4}\texte^{-2\slb\sqrt t}Z_t\,\,\,\,\underset{t\to\infty}\lawarrow\,\,\,W,
\end{equation}
where $W$ is a measure governing the extrema of a tree-indexed Branching Random Walk; see \eqref{E:3.9c} and \eqref{E:1-step} below.
In \cite[Theorem~1.3]{BL4} this was checked for the convergence of the total mass of these measures. The prefactors on the left of \eqref{E:1.25} arise from the differences between the centering sequences~$\sqrt t+a_n(t)$ and~$m_n$. 
%\end{remark}

\subsection{The cluster law}
\label{sec-1.3}\noindent
In the above limit results, we have not yet identified the cluster law~$\DD$ beyond its basic properties. Perhaps the most important part of our conclusions is that we can characterize~$\DD$ quite explicitly. We recall a definition:

\begin{definition}
Given a locally-finite random point measure $\eta$ on~$\R$, a Branching Random Walk with step distribution~$\eta$ started at $x\in\R$ is a sequence
$\{\xi_k \}_{k \ge 0}$
of random locally-finite point measures on~$\R$ such that $\xi_0=\delta_x$ a.s.\ and, denoting
$\scrF_k:=\sigma(\xi_j\colon j\le k)$,
\begin{equation}
\label{E:3.9c}
E\bigl(\texte^{-\langle\xi_{k+1},f\rangle}\,|\,\scrF_k\bigr) =\texte^{-\langle\xi_k,\tilde f\rangle}\quad\text{where}\quad\tilde f(h):=
- \log E \bigl(\texte^{-\langle \eta, f(h+\cdot) \rangle} \bigr)
\end{equation}
holds for all $k\ge0$ and all continuous $f\colon\R\to[0,\infty)$  with compact support.
\end{definition}

The formula \eqref{E:3.9c} expresses that, to get~$\xi_{k+1}$ from~$\xi_k$, we just replace each point~$x$ of the process~$\xi_k$ by an independent copy of~$\eta$ shifted by~$x$. The choice of the step distribution~$\eta$ relevant for the local-time problem at hand is
\begin{equation}
\label{E:1-step}
\eta:=\sum_{i=1}^b\delta_{Y_i}\quad\text{\rm for}\quad Y_1,\dots,Y_b\text{\rm\ \ i.i.d.\ \ }\NN(0,1/2),
\end{equation}
where~$\NN(\mu,\sigma^2)$ denotes normal distribution with mean~$\mu$ and variance~$\sigma^2$.
Since the family-tree of the associated branching process is the regular rooted tree of forward degree~$b$, this Branching Random Walk can alternatively be viewed as a tree-indexed Markov chain with step distribution~$\NN(0,1/2)$.

Under suitable moment conditions on~$\eta$, which comfortably include \eqref{E:1-step}, significant efforts culminating in the work of A\"id\'ekon~\cite{Aidekon} and Madaule~\cite{Madaule2} showed that there exist constants~$c_1,c_2\in\R$ and $\alpha>0$ such that, for $\wt m_n:=c_1 n- c_2\log n$ and any continuous $f\colon\R\to[0,\infty)$ with compact support,
\begin{equation}
\label{E:6.1}
E\bigl(\texte^{-\langle\xi_n,f(\cdot-\wt m_n)\rangle}\bigr)
\,\,\underset{n\to\infty}\longrightarrow\,\,
E\biggl(\exp\Bigl\{-\eusb W\int \texte^{-\alpha h}\textd h\otimes\DD'(\textd\chi)\bigl[1-\texte^{-\int f(h+\cdot))\textd\chi}\bigr]\Bigr\}\biggr).
\end{equation}
Here~$\eusb W$ is an a.s.\ finite and positive random variable and~$\DD'$ is a deterministic law on locally finite point processes such that a.e.~sample~$\chi$ from~$\DD'$ obeys $\supp(\chi)\subseteq(-\infty,0]$ and~$\sup\supp(\chi)=0$. 

The interpretation of \eqref{E:6.1} is that, as~$n\to\infty$, the point process~$\xi_n$ shifted by~$\wt m_n$ tends in law (relative to the vague topology) to a point process of the form
\begin{equation}
%\label{}
\sum_{i\ge1}\sum_{j\ge1}\delta_{\alpha^{-1}\log \eusb W+h_i+d_j^{(i)}}
\end{equation}
where~$\{h_i\}_{i\ge1}$ are the points of~$\PPP(\texte^{-\alpha h}\textd h)$ while $\{d_j^{(i)}\colon j\ge1\}_{i\ge1}$ are the points in i.i.d.\ samples from~$\DD'$, with the random objects $\eusb W$, $\{h_i\}_{i\ge1}$ and $\{d_j^{(i)}\colon j\ge1\}_{i\ge1}$ sampled independently. Both~$\alpha$ and~$\DD'$ depend on the law of~$\eta$.

Our proof of Theorem~\ref{thm-2} is readily modified to give an independent proof of \eqref{E:6.1} for our Branching Random Walk and thus show: 

\begin{corollary}
\label{cor-1.5}
The cluster law~$\DD$ from Theorems~\ref{thm-1} and~\ref{thm-2} is the cluster law~$\DD'$ from \eqref{E:6.1} for the Branching Random Walk with step distribution \eqref{E:1-step}.
\end{corollary}

In conclusion, we have found that, while the intensity measure governing the spatial position of the extremal values of the local time is specific to the local time process at hand, the local structure of the configuration near the extremal points is universal. This can be attributed to the fact that, for each~$n\ge1$ fixed,
\begin{equation}
\label{E:1.30}
\sum_{x\in\BbbL_n}\delta_{\sqrt{L_t(x)}-\sqrt t}\,\,\,\,\underset{t\to\infty}\lawarrow\,\,\,\xi_n,
\end{equation}
where~$\xi_n$ is the state at time~$n$ of the Branching Random Walk with step distribution \eqref{E:1-step}. The statement \eqref{E:1.30} is deduced from a Multivariate CLT; the proof of above results relies on a stronger version encoded in the form of the Second Ray-Knight Theorem of Eisenbaum, Kaspi, Marcus, Rosen and Shi~\cite{EKMRS}; see~Lemma~\ref{lemma-D}.

%\newpage
\section{Connections, discussion and outline}
\label{sec-2}\noindent
We proceed by discussing a broader context of extremal behavior of logarithmically correlated processes. We also outline the main steps of our proof and mention alternative approaches to local convergence we considered. 

\subsection{Extremal properties of random walks}
The study of extremal properties of random walks dates back to an influential study from 1960 by Erd\H os and Taylor~\cite{ET60}. There, among other things, the authors addressed the time~$T_n^\star$ that a $d$-dimensional simple symmetric random walk of~$n$ steps spends at its most visited site, dubbed a \emph{frequent point}. The answer turns out to be most interesting in dimension $d=2$ where a resolution of just the leading order term 
\begin{equation}
%\label{}
\frac{T_n^\star}{(\log n)^2}\,\underset{n\to\infty}{\overset{\text{\rm a.s.}}\longrightarrow}\,\, \frac1\pi 
\end{equation}
 was given full 40 years later by Dembo, Peres, Rosen and Zeitouni~\cite{DPRZ01}. What makes this case hard is the (asymptotic) scale invariance of the random walk which causes~$T_n^\star$ to collect non-trivial contributions from excursions on all spatial scales.

Further progress occurred on this and related questions over the last decade. For instance, Abe~\cite{A15} proved a result analogous to~\cite{DPRZ01} for the walk on an $N\times N$ torus in~$\Z^2$ run for times of order~$N^2(\log N)^2$. Jego~\cite{J18} generalized the conclusions of \cite{DPRZ01} to a large class of random walks. Aiming for an asymptotic limit law, Rosen~\cite{Rosen-new} established tightness around an explicit deterministic centering sequence of the kind \eqref{E:1.4} for a two-dimensional Brownian motion stopped upon its first exit from a bounded domain. Jego~\cite{J20} in turn constructed a candidate for the measure that should govern the law of the maximum via a formula of the kind~\eqref{E:1.2}. Still, proving such a formula rigorously seems elusive at present.

Motivated by this, the present authors and O.~Louidor considered the problem of simple random walk on a $b$-ary tree. This setting bears a number of close connections to the random walk on~$\Z^2$. Indeed, for~$b=4$ the set of leaves~$\BbbL_n$ can be identified with points in the square
\begin{equation}
%\label{}
\{0,\dots,2^n-1\}\times \{0,\dots,2^n-1\}
\end{equation}
 in~$\Z^2$; the simple random walk on this square is then imitated by the Markov chain on~$\BbbL_n$ obtained by recording the successive visits to the leaves of simple random walk on~$\T_n$. (A more quantitative connection is seen in, e.g., the behavior of the Green function.) Killing the walk on~$\T_n$ on its first visit to~$\varrho$ then corresponds to killing the lattice random walk upon its first exit from the square.
 
In~\cite{BL4}, O.~Louidor and the second author proved a limit result for the maximal local time in a (variable speed) continuous time setting and the walk started from both the leaves and the root. The paper \cite{BL4} draws on earlier work of the first author~\cite{A18} where the whole process of extreme local maxima was controlled for the walk started from and parametrized by time spent at the root, albeit only for times~$t$ that grow at least as constant times $n\log n$. The present paper extends this to all times and adds control of the clusters ``hanging off'' the local maxima. 

Taken from a larger perspective, our results add another instance to the universality class associated with extremal behavior of logarithmically correlated processes. Other contexts in which similar conclusions have been proved include:
\settowidth{\leftmargini}{(11)}
\begin{itemize}
\item
Branching Brownian motion (Arguin, Bovier and Kistler~\cite{ABK1,ABK2,ABK3}, 
A\"id\'ekon, Berestycki, Brunet and Shi~\cite{ABBS}),
\item
critical Branching Random Walks (A\"id\'ekon~\cite{Aidekon}, Madaule~\cite{Madaule2}),
\item
Gaussian Free Field in finite subsets of~$\Z^2$ (Bramson, Ding and Zeitouni~\cite{BDingZ}, Biskup and Louidor~\cite{BL1,BL2,BL3}),
\item more general logarithmically correlated Gaussian processes (Madaule~\cite{Madaule1}, Ding, Roy and Zeitouni~\cite{DRZ}, Schweiger and Zeitouni~\cite{SZ}) including the four-dimensional membrane model (Schweiger~\cite{Schweiger}),
\item
characteristic polynomial of a random matrix ensemble (Paquette and Zeitouni~\cite{PZ}),
\item
a class of non-Gaussian fields on a torus (Bauerschmidt-Hofstetter~\cite{BHo}, Hofstetter~\cite{H}, Barashkov, Gunaratnam and Hofstetter~\cite{BGH}),
\item
subcritical hierarchical DG-model (Biskup and Huang~\cite{BH}).
\end{itemize}
Tightness of the maximum has recently been shown for uniformly-convex Ginzburg-Landau 
models (Wu and Zeitouni~\cite{WZ}, Schweiger, Wu and Zeitouni~\cite{SWZ}). 

A universal feature of these problems is that a suitably centered maximum tends, as the system size increases to infinity, to a randomly shifted Gumbel law while the extremal process tends to a decorated Cox process. The law of the random shift and the decorations are typically (but, as our results show, not always) specific to the problem at hand as they depend on the nature of global and local correlations.

\subsection{Proof baseline}
\label{sec-2.2}\noindent
We will now move to a discussion of our proofs. As noted in the introduction, the main result to be proved is Theorem~\ref{thm-2} from which the rest of the conclusions follow by relatively soft means. The setting of Theorem~\ref{thm-2} is more tractable thanks to:

\begin{lemma}[Markov property]
\label{lemma-M}
Suppose that~$x_1,\dots,x_m\in\T_n$ are vertices such that the subtrees $\T^{(1)},\dots,\T^{(m)}$ of~$\T_n$ rooted at these vertices are vertex-disjoint. Write~$n_i$ for the depth of~$\T^{(i)}$ and denote
\begin{equation}
%\label{}
\V(x_1,\dots,x_m):=\{x_1,\dots,x_m\}\cup\Bigl(\T_n\smallsetminus\bigcup_{i=1}^m\T^{(i)}\Bigr).
\end{equation}
Then, for   all   $t>0$, conditional on $\{L_t(x)\colon x\in\V(x_1,\dots,x_m)\}$, the families
\begin{equation}
%\label{}
\{L_t(x)\colon x\in\T^{(i)}\}_{i=1}^m
\end{equation}
are independent with the~$i$-th family distributed as $\{L_u(x)\colon x\in\T_{n_i}\}$ for~$u:=L_t(x_i)$.
\end{lemma}

\begin{proofsect}{Proof}
This is a direct consequence of  the geometry of the tree and the reliance on exponential clocks to run the underlying random walk. See, e.g., \cite[Lemma~2.5]{BL4}.
\end{proofsect}

A short way to state the above is that $\{L_t(x)\colon x\in\T_n\}$ is a tree-indexed Markov chain. The step distribution is quite explicit; indeed, 
given the value~$L_t(x)$, 
the values of~$L_t$ at the ``children''~$z_1,\dots,z_b$ of~$x$ are independent with each $L_t(z_i)$ having the law of the sum of a Poisson$(L_t(x))$-number of i.i.d.\ Exponentials with parameter~$1$.

Another convenient fact (proved in, e.g., Belius, Rosen and Zeitouni~\cite[Lemma 3.1]{BRZ}) is that the law of $\sqrt{L_t}$ on a path $x_0=\varrho,\dots,x_n=x$ from the root to ~$x\in\BbbL_n$ can be encoded via the zero-dimensional Bessel process. This permits squeezing $\sqrt{L_t}$ on the path by a barrier event (see~\cite[Lemma~3.2]{A18} and \cite[Lemma~2.9]{BL4}) which then shows that, for $\sqrt{L_t(x)}=m_n+O(1)$, the value $L_t(x_k)$ for~$k$ and~$n-k$ large is unusually large and the value~$L_t(x)$ is unusually large given the value $L_t(x_k)$. 

The above shows that, if we are after the absolute maximum of~$L_t$ on the leaves, all we need is the leading-order asymptotic of the conditional law of $\max_{x\in\BbbL_{n-k}(x_k)}L_t(x)$, where
$\BbbL_{n-k}(x_k)$ is the set of leaves of the subtree rooted at $x_k$,
given that~$L_t(x_k)$ is large. In~\cite{BL4} this was supplied by:

\begin{proposition}
\label{prop-3}
There exists $c_\star\in(0,\infty)$ such that the quantity $\epsilon_{n,t,u}$ defined for integer $n\ge1$ and real $t>0$ and~$u>0$ by
\begin{equation}
\label{E:3.34u}
P^\varrho \Bigl(\,\max_{x\in\BbbL_n} \sqrt{L_t(x)} - \sqrt{t} - a_n(t) > u \Bigr)
= c_\star \bigl(1+\epsilon_{n,t,u}\bigr) u\texte^{-2u\slb} 
\end{equation}
obeys
\begin{equation}
\label{E:3.29ui}
\lim_{m\to\infty}\sup_{t,u\ge m} \limsup_{n\to\infty}\,\bigl|\,\epsilon_{n,t,u}\bigr|=0.
\end{equation}
\end{proposition}

\begin{proofsect}{Proof}
This is a restatement of~\cite[Proposition~3.5]{BL4}.
\end{proofsect}

Here the centering must be done by $\sqrt{t}+a_n(t)$ as we need to use the result for~$t$ replaced by $L_t(x_k)$ in the regime when that is large. The difference in the centering matters; indeed, the factor
\begin{equation}
%\label{}
L_t(x)^{1/4}\,
	\texte^{2\slb\,\sqrt{L_t(x)}}
\end{equation}
in \eqref{E:1.18} arises directly from $\sqrt{t}+a_{n-k}(t)-m_n$ for $t:=L_t(x_k)$.

The asymptotic \eqref{E:3.34u} is actually sufficient to extract a weak limit of the process of  \emph{extreme local maxima}, which are the leaf-vertices with $\sqrt{L_t}$-values near~$m_n$ that dominate a large neighborhood (relative to the intrinsic tree metric) thereof. (This has already been done in~\cite{A18} for the problem at hand in the regime when $t\ge c n\log n$. The control of extreme local maxima is an essential step in the analysis of the extremal process of the two-dimensional GFF; see~\cite{BL1}.) Each extreme local maximum comes with a cluster of nearby values, to be called  \emph{decorations}, which are strongly correlated but, as it turns out, with a limit law that is asymptotically independent of the local maximum. The control of the decorations constitutes the bulk of the proof. 

\subsection{Strategies for limit of the decorations}
We have actually developed and written up most of the details of three possible approaches to the limit of the decorations. Using the notation~$x_0=\varrho,x_1,\dots,x_n=x$ for a path from the root~$\varrho$ to~$x\in\BbbL_n$, our first argument is modeled on the proofs for Branching Brownian motion whose strategy is to show that, conditional on~$L_t(x_k)$ large, the whole process on the leaves of the subtree rooted at~$x_k$ converges as~$n\to\infty$ followed by~$k\to\infty$. This does work but, since the calculations have to be done for the local time process, the arguments quickly become very technical. In this approach, we do not get information about the law of the decorations.

Another approach we developed is based on the observation that any subsequential weak limit~$\zeta$ of the whole extremal process is invariant under post-composition with one step of the Branching Random Walk with step distribution
\begin{equation}
\label{E:2-step}
\tilde\eta:=\sum_{i=1}^b\delta_{Y_i}\quad\text{\rm for}\quad Y_1,\dots,Y_b\text{\rm\ \ i.i.d.\ \ }\NN(-\slb,1/2).
\end{equation}
This means that
\begin{equation}
\label{E:2.9}
E\bigl(\texte^{-\langle\zeta,f\rangle}\bigr) =E\bigl(\texte^{-\langle\zeta,\tilde f\rangle}\bigr)\quad\text{for}\quad\tilde f(h):=-\log E\bigl(\texte^{-\langle\tilde\eta,f(h+\cdot)\rangle}\bigr)
\end{equation}
holds for each non-negative continuous~$f$ with compact support.

In their work on pre-composition invariant processes, Maillard and Mallein~\cite{MM} conjectured that, for~$\tilde\eta$ satisfying certain moment assumptions, any point measure~$\zeta$  satisfying \eqref{E:2.9} has the law of a randomly-shifted decorated Poisson point process; see \cite[Conjecture~1.3]{MM}. This has so far been verified only for~$\tilde\eta$ corresponding to Branching Brownian motion by Kabluchko~\cite{Kabluchko} (supercritical drift) and by Chen, Garban and Shekhar~\cite{CGS} (critical drift). We did manage to do the same for our Branching Random Walk although the estimates become quite cumbersome here as well.

We will therefore follow yet another approach that capitalizes on the fact that for~$t$ very large, $\sqrt{L_t}$ is well approximated by the Gaussian process $\{h_x\colon x\in\T_n\}$ such that
\begin{equation}
\label{E:2.10a}
h_\varrho=0 \,\,\,\text{ and }\,\,\,\bigl\{h(x)-h(\frakm(x))\colon x\in\T_n\smallsetminus\{\varrho\}\bigr\} \text{ are i.i.d. }\NN(0,1/2),
\end{equation}
where~$\frakm(x)$ denotes the ``mother'' vertex of~$x$; i.e., the nearest vertex on the path from~$x$ to the root.
We will call this process the \emph{Gaussian Free Field} (GFF) on~$\T_n$, although it can also be viewed as a Branching Random Walk on its genealogical tree with~$\eta$ given by \eqref{E:1-step} and thus also as a tree-indexed Markov chain with step distribution~$\NN(0,1/2)$.

The said approximation can be formulated as a limit theorem that was noted already in~\eqref{E:1.30}. Explicitly, since~$L_t$ is the sum of a Poisson$(t)$-number of i.i.d.\ fields with mean one at each vertex and (as one has to check) covariance being~$4t$-times
that of~$h$, the multivariate CLT gives
\begin{equation}
%\label{} 
\bigl\{\sqrt{L_t(x)}-\sqrt t\colon x\in\T_n\bigr\}\,\,\,\underset{t\to\infty}\lawarrow\,\,\, \bigl\{h(x)\colon x\in\T_n\bigr\}.
\end{equation}
As also mentioned earlier, our proof will actually make use of a stronger connection whose advantage is that it works even for finite~$t$:

\begin{lemma}[Isomorphism theorem]
\label{lemma-D}
For each $n\ge0$ and each~$t>0$ there exists a coupling of~$L_t$ and two copies~$h$ and~$\tilde h$ of~GFF on~$\T_n$ such that
\begin{equation}
%\label{}
h\independent L_t 
\end{equation}
and
\begin{equation}
\label{E:3.20}
L_t(x)+h(x)^2 = \bigl(\tilde h(x)+\sqrt t\bigr )^2
\end{equation}
holds pointwise a.s.\ for all~$x\in\T_n$. (Note: $\tilde h$ and~$L_t$ are not independent.)
\end{lemma}

\begin{proofsect}{Proof}
The Second Ray-Knight Theorem of Eisenbaum, Kaspi, Marcus, Rosen and Shi~\cite{EKMRS} gives \eqref{E:3.20} as equality in distribution. (The reader alerted by absence of factors of~$1/2$ in \eqref{E:3.20}  beware that our GFF is scaled by $1/\sqrt2$ compared to the one usually used.) Zhai~\cite{Zhai} extended this to a pointwise equality by drawing the signs of $\tilde h(x)+\sqrt t$ from the corresponding conditional law.
\end{proofsect}

Our argument then goes as follows: We first show that the extremal process associated with~$\sqrt{L_t}$ is not significantly affected by swapping the increment of~$\sqrt{L_t}$ in the last~$k$ generations of the tree for an increment of the GFF. This relies on the above coupling and the observation that, conditional on~$L_t(z)$ with~$z:=\frakm^k(x)$,
\begin{equation}
\label{E:2.14}
\begin{aligned}
\sqrt{L_t(x)} &= \sqrt{L_t(x)+[h(x)-h(z)]^2} + O\biggl(\frac{[h(x)-h(z)]^2}{\sqrt{L_t(z)}}\biggr) 
\\
&= \tilde h(x)-\tilde h(z)+\sqrt{L_t(z)}+ O\biggl(\frac{[h(x)-h(z)]^2}{\sqrt{L_t(z)}}\biggr),
\end{aligned}
\end{equation}
see Proposition~\ref{prop-3.1}. Here the error term is negligible because, roughly speaking, $\sqrt{L_t(z)}$ scales with~$n$ but $h(x)-h(z)$ does not. Since~$\tilde h$ is known to have a positive ``gap'' between the first and the second largest value, we also get that the location of the extreme local maxima of~$\sqrt{L_t}$ and~$\tilde h$ agree with high probability.

The argument we just gave shows that the local extremal structure of~$\sqrt{L_t}$ is asymptotically that of~$\tilde h$. As it turns out, the law of the GFF decorations near an extremely large local maximum converges; see Proposition~\ref{prop-3.3} for a precise statement. This effectively reduces the problem to the process of extreme local maxima which, as explained earlier, is handled using Proposition~\ref{prop-3}.

We note that the connection between the local time and the Gaussian Free Field stated in Lemma~\ref{lemma-D} is a powerful tool that drove a number of earlier sharp conclusions about the extremal properties of the local time of random walks. This includes studies of the limit law of the cover time (Ding~\cite{Ding}, Cortines, Louidor and Saglietti~\cite{CLS}, Louidor and Saglietti~\cite{LS}) as well as the intermediate level sets, a.k.a.\ thick points, of random walks in planar domains (Abe and Biskup~\cite{AB}, Abe, Biskup and Lee~\cite{ABL}; see~\cite{B-Budapest} for a review).

\subsection{Outline}
The remainder of this paper is organized as follows. In Section~\ref{sec-3} we prove Theorem~\ref{thm-2} in a structured form  that records the extremal processes using the position of a nearby ``local maximum'' and the ``shape'' of the configuration relative to it; see Theorem~\ref{thm-3.1}. The proof relies on two technical propositions whose details are worked out in Section~\ref{sec-4} and Section~\ref{sec-6}, respectively. The remaining results are proved in Section~\ref{sec-5}. 

%\newpage
\section{The walk started at the root}
\label{sec-3}\noindent
The proofs of our main results commence with Theorem~\ref{thm-2} that we prove here assuming two technical propositions. These propositions encapsulate technical arguments whose immediate inclusion would detract from the main line of proof. 

\subsection{The structured extremal process}
Our proof of Theorem~\ref{thm-2} adopts the strategy developed for the two-dimensional Discrete Gaussian Free field by O.~Louidor and the second author~\cite{BL1,BL2,BL3}. In particular, we will describe the extremal process in a more structured way by keeping track of local maxima of the field along with the whole configuration nearby. 

\begin{figure}[t]
\refstepcounter{obrazek}
\label{fig1}
\vglue2mm
\centerline{\includegraphics[width=4.3in]{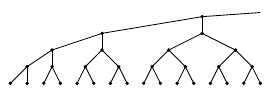}}
\vglue2mm
\begin{quote}
\fontsize{9}{5}\selectfont
{\sc Fig~\theobrazek.}\ The canopy-tree representing the portion of~$\T_n$ near the special vertex~$\zero$. The upper path is the path from the root to~$\zero$. The leaves are then identified with~$\N$.
\normalsize
\end{quote}
\end{figure}

We start with some definitions and notation.  Note that, given any $n\ge1$ and any leaf vertex $x=(x_1,\dots,x_n)\in\BbbL_n$, the map~$\pi_{n,x}$ assigning
\begin{equation}
\label{E:3.1}
\pi_{n,x}(y):=\sum_{k=1}^{n}(y_k-x_k \text{ mod }b)b^{n-k}
\end{equation}
to each $y=(y_1,\dots,y_n)\in\BbbL_n$ embeds~$\BbbL_n$ naturally into~$\N$.  (This map can be thought of as induced from the natural embedding of~$\T_n$ to a canopy tree taking~$x$ to the ``origin'' thereof; see Fig.~\ref{fig1}.) Given $h\colon\BbbL_n\to\R$ and~$x\in\mathbb{L}_n$, 
we write $h(x\cdot)\colon\N\to\R$ for the map
\begin{equation}
h(xj):=\begin{cases}
h\circ\pi_{n,x}^{-1}(j),\qquad&\text{if }j\in \pi_{n,x}(\BbbL_n),
\\
0,\qquad&\text{if }j\in\N\smallsetminus \pi_{n,x}(\BbbL_n),
\end{cases}
\end{equation}
that enumerates the values of~$h$ on~$\T_n$ increasingly in the graph-theoretical distance from~$x$ while breaking ties using the mod-$b$ rule in \eqref{E:3.1}.

Let~$B_k(x)$ denote the set of vertices in~$\BbbL_n$ that are at most~$2k$-steps away from~$x$ in the graph-theoretical distance on~$\T_n$. For $n\ge k\ge1$ and a function~$\varphi\colon\BbbL_n\to\R$, let
\begin{equation}
\label{E:3.2}
\MM_{n,k}(\varphi,x):=\Bigl\{\max_{y\in B_k(x)}\varphi(y) = \varphi(x)\Bigr\}
\end{equation}
be the event that~$x$ is a ``$k$-local maximum'' of~$\varphi$.

If~$x\in\BbbL_n$ is a $k$-local maximum of $\sqrt{L_t(\cdot)}$, then a natural way to encode the ``shape'' of the configuration near~$x$ is through the function $\sqrt{L_t(x)}-\sqrt{L_t(x\cdot)}$. Including also the position of and the value of~$\sqrt{L_t}$ at the local maximum, this leads us to a description based on the random Borel measure on~$[0,1]\times\R\times\R^\N$ defined by
\begin{equation}
%\label{}
\zeta_{n,k}^{(t)}:=\sum_{x\in\BbbL_n}1_{\MM_{n,k}(\sqrt{L_t},\,x)}\,\delta_{\theta_n(x)}\otimes\delta_{\sqrt{L_t(x)}-m_n}\otimes\delta_{\sqrt{L_t(x)}-\sqrt{L_t(x\cdot)}},
\end{equation}
where $m_n$ is as in \eqref{E:1.4}.
We will call $\zeta_{n,k}^{(t)}$ the \textit{structured extremal process}.

%%%% COUNTER FOR FIGURES
\newcounter{obrazek}

The structured extremal process may behave somewhat strange in the branches of~$\T_n$ where the random walk does not hit~$\BbbL_n$ by time~$t$. However, as $n\to\infty$, the shift by~$m_n$ effectively moves such points to~$-\infty$ in the second coordinate thus making them undetectable by integrals against compactly supported functions. Similarly, our convention about the values of the third component that do not come from~$\T_n$ forces us to probe $\zeta_{n,k}^{(t)}$ only by functions that depend on a finite number of coordinates of the third component. 
We are thus let to consider the space~$\Cloc^+$ of functions $f\colon[0,1]\times\R\times\R^\N\to[0,\infty)$ that depend only on coordinates in~$[0,1]\times\R\times\R^{\{0,\dots,j\}}$ for some~$j\ge0$ and, when restricted to these coordinates, are continuous with compact support. 

\smallskip
Our main result about the structured extremal process is~then:

\begin{theorem}
\label{thm-3.1}
For each~$t>0$, the weak limit \eqref{E:1.18} exists and defines an a.s.-finite random Borel measure~$Z_t$ on~$[0,1]$. Moreover, there exists a Borel probability measure~$\nu$ on~$\R^\N$ such that for all~$t>0$, all~$f\in\Cloc^+$ and all sequences $\{k_n\}_{n\in\N}$ with $k_n\to\infty$ and $n-k_n\to\infty$, 
\begin{equation}
\label{E:3.4}
\begin{aligned}
E^\varrho\bigl(&\texte^{-\langle\zeta_{n,k_n}^{(t)},f\rangle}\bigr)
\\
&\,\,\,\underset{n\to\infty}\longrightarrow\,\, \E\biggl(\exp\Bigl\{-\int Z_t(\textd x)\otimes\texte^{-2\slb\, h}\textd h\otimes\nu(\textd\phi)(1-\texte^{-f(x,h,\phi)})\Bigr\}\biggr),
\end{aligned}
\end{equation}
where the expectation on the right is with respect to the law of~$Z_t$. The measure $Z_t$ obeys \eqref{E:1.13},
$\nu$ is concentrated on~$\{\phi_0=0\}\cap
\bigcap_{x\in\N\smallsetminus \{0\}}\{\phi_x>0\}$ and the set $\{x\in\N\colon \phi_x\le a\}$ is finite~$\nu$-a.s.\ for each~$a>0$. A.e.\ sample of~$Z_t$ is diffuse (i.e., does not charge singletons). 
\end{theorem}

The statement \eqref{E:3.4} implies that $\zeta_{n,k_n}^{(t)}$ tends in law to a Poisson point process on $[0,1]\times\R\times\R^\N$ with the intensity measure
\begin{equation}
%\label{}
Z_t(\textd x)\otimes\texte^{-2\slb \,h}\textd h\otimes\nu(\textd\phi).
\end{equation}
This follows from the fact that \eqref{E:3.4} extends to all continuous compactly supported functions~$f$ on $[0,1]\times\R\times\R^\N$ endowed with product topology. 

\subsection{Proof of Theorem~\ref{thm-2}}
Before the delve into the proof of Theorem~\ref{thm-3.1}, let us see how it implies our main result from Section~\ref{sec-1.2}. We will need two lemmas.

\begin{lemma}
\label{lemma-3.2a}
Recall the definition of $a_n(t)$ from \eqref{E:1.24a}. There exist~$c,\tilde c >0$ such that
\begin{equation}
\label{E:3.6x}
P^\varrho\Bigl(\,\max_{x\in\BbbL_n}\sqrt{L_t(x)}\ge\sqrt t+a_n(t)+u\Bigr)\le c(1+u)\texte^{-2\slb\, u-\tilde c u^2/n}
\end{equation}
holds for all $n\ge1$,  $t>0$ and all~$u>0$. Moreover, there exist~$c_1,c_2>0$ such that \begin{equation}
\label{E:3.8w}
P^\varrho\biggl(\,\Bigl|\,\max_{x\in\BbbL_n}\sqrt{L_t(x)}-\sqrt t-a_n(t\vee1)\Bigr|>u\,\bigg|\,\max_{x\in\BbbL_n}L_t(x)>0\biggr)\le c_1\texte^{-c_2 u}
\end{equation}
holds for all~$n\ge1$, all $t>0$ and all $u\in(0,n]$.  Here $s\vee t:=\max\{s,t\}$.
\end{lemma}

\begin{proofsect}{Proof}
The inequality \eqref{E:3.6x} is a restatement of~\cite[Proposition~3.1]{A18}. The inequality~\eqref{E:3.8w} appeared in \cite[Theorem 2.1]{BL4}.
\end{proofsect}

\begin{lemma}
\label{lemma-3.3a}
For $n\ge1$,~$t>0$ and~$\lambda>0$, denote
\begin{equation}
\label{E:3.12}
\Gamma(\lambda):=\bigl\{x\in\BbbL_n\colon \sqrt{L_t(x)}\ge m_n-\lambda\bigr\}.
\end{equation}
Then for all~$\lambda>0$,
\begin{equation}
\label{E:3.8a}
\lim_{k\to\infty}\,\limsup_{n\to\infty}\,P^\varrho\Bigl(\exists x,y\in\Gamma(\lambda)\colon y\in B_{n-k}(x)\smallsetminus B_k(x)\Bigr) = 0.
\end{equation}
\end{lemma}

\begin{proofsect}{Proof}
This is~\cite[Proposition~4.1]{A18} with~$m_n$ instead of $a_n(t)+\sqrt t$, which in light of monotonicity in~$\lambda$ is immaterial as $|a_n(t)-m_n|$ is bounded by a~$t$-dependent constant uniformly in~$n\ge1$. (The statement in \cite[Proposition~4.1]{A18} actually bounds the probability with $\lambda:=c\log k$ by $c' k^{-1/8}$ for~$n$ large.)
\end{proofsect}

The above lemmas show that, with~$t>0$ fixed, the centered maximum of~$\sqrt{L_t}$ is tight and the level sets close to~$m_n$ are clustered; meaning that any two points of the level set are either within $O(1)$ or $n-O(1)$ in graph-theoretical distance on~$\T_n$ from each other. Elementary geometric consideration then yield: 

\begin{corollary}
\label{cor-3.4}
For any~$\lambda>0$,
\begin{equation}
%\label{}
\lim_{\ell\to\infty}\limsup_{n\to\infty}P^\varrho\Bigl(\bigl|\Gamma(\lambda)\bigr|>\ell\Bigr)=0.
\end{equation}
\end{corollary}

\begin{proofsect}{Proof}
The set~$\BbbL_n$ can be covered by~$b^k$ sets of the form $B_{n-k}(x)$ each of which contains at most~$b^k$ points from~$\Gamma(\lambda)$, unless the event in \eqref{E:3.8a} occurs. The claim thus follows from Lemma~\ref{lemma-3.3a}.
\end{proofsect}

We are now ready for:

\begin{proofsect}{Proof of Theorem~\ref{thm-2} from Theorem~\ref{thm-3.1}}
It suffices to show that, for each $t>0$ and each $f\colon[0,1]\times\R\to[0,\infty)$  continuous with compact support,
\begin{equation}
\label{E:1.8a}
\begin{aligned}
E^\varrho\biggl(\exp&\Bigl\{-\sum_{x\in\BbbL_n}f\bigl(\theta_n(x),\sqrt{L_t(x)}-m_n\bigr)\Bigr\}\biggr)
\\
&\underset{n\to\infty}\longrightarrow\,\,\E\biggl(\exp\Bigl\{-\int Z_t(\textd x)\otimes\texte^{-2h\slb}\textd h\otimes\DD(\textd\chi)\bigl(1-\texte^{-\int f(x,h+\cdot)\textd\chi}\bigr)\Bigr\}\biggr)
\end{aligned}
\end{equation}
holds with $Z_t$ as in \eqref{E:3.4} and~$\DD$ being the law of
\begin{equation}
\label{E:3.13ui}
\chi^\phi:=\sum_{x\in\N}\delta_{-\phi_x}
\end{equation}
 for~$\phi$ sampled from~$\nu$.
Note that the properties of~$\nu$ stated in Theorem~\ref{thm-3.1} ensure that $\chi^\phi$ is $\nu$-a.s.\ a Radon measure on~$\R$ with~$\supp(\chi^\phi)\subseteq(-\infty,0]$ and~$\chi^\phi(\{0\})=1$.

Let~$g\colon[0,\infty)\to[0,1]$ be continuous, non-increasing with~$\supp(g)\subseteq[0,2]$ and~$g=1$ on~$[0,1]$. For any~$r\ge1$, the function $\tilde f_r\colon[0,1]\times\R\times\R^\N\to[0,\infty)$ defined by
\begin{equation}
%\label{}
\tilde f_r(x,h,\phi):=\sum_{y=0}^{b^{r}-1} f(x,h-\phi_y)g\bigl(|h|/r\bigr)
\end{equation}
belongs to~$\Cloc^+$ and so \eqref{E:3.4} applies. Aiming to take~$r\to\infty$ limit of the resulting expression, note that~$\tilde f_r(x,h,\phi)$ increases to $\int f(x,h+\cdot)\textd\chi^\phi$ as $r\to\infty$. The Monotone Convergence Theorem yields
\begin{equation}
%\label{}
\begin{aligned}
\int Z_t(\textd x)&\otimes\texte^{-2\slb\, h}\textd h\otimes\nu(\textd\phi)(1-\texte^{-\tilde f_r(x,h,\phi)})
\\
&\underset{r\to\infty}\longrightarrow\,\,
\int Z_t(\textd x)\otimes\texte^{-2\slb\, h}\textd h\otimes\DD(\textd\chi)\bigl(1-\texte^{-\int f(x,h+\cdot)\textd\chi}\bigr)
\end{aligned}
\end{equation}
and the Bounded Convergence Theorem then shows that, as~$r\to\infty$, the right-hand side of \eqref{E:3.4} for~$f$ replaced by~$\tilde f_r$ tends to that of \eqref{E:1.8a}. Hence, as soon as we prove that
\begin{equation}
\label{E:3.10}
\lim_{r\to\infty}\limsup_{n\to\infty}\,P^\varrho\Biggl(\biggl|\langle\zeta_{n,k_n}^{(t)},\tilde f_r\rangle-\sum_{x\in\BbbL_n}f\bigl(\theta_n(x),\sqrt{L_t(x)}-m_n\bigr)\biggr|>0\Biggr)=0
\end{equation}
also the left-hand side of \eqref{E:3.4} asymptotically approaches that of \eqref{E:1.8a} and so \eqref{E:1.8a} is inferred from~\eqref{E:3.4}. 

For \eqref{E:3.10}, let  $\lambda>0$ be such that $\supp(f)\subseteq[0,1]\times[-\lambda,\lambda]$. Then $r>\lambda$ and $\max_{y\in\BbbL_n}\sqrt{L_t(y)}-m_n\le r$  imply $-r\le-\lambda\le\sqrt{L_t(x)}-m_n\le r$ at every $k_n$-local maximum possibly contributing to $\langle\zeta_{n,k_n}^{(t)},\tilde f_r\rangle$. The truncation by~$g$ then becomes immaterial and we get 
\begin{equation}
\label{E:3.11}
\langle\zeta_{n,k_n}^{(t)},\tilde f_r\rangle=\sum_{x\in\BbbL_n}1_{\MM_{n,k_n}(\sqrt{L_t},\,x)}\sum_{y\in B_r(x)}f\bigl(\theta_n(y),{\textstyle\sqrt{L_t(y)}}-m_n\bigr).
\end{equation}
Assuming $m_n>\lambda$, which guarantees that only terms with positive values of~$L_t$ can be local maxima contributing to \eqref{E:3.11}, and $2r<k_n$ along with (the full-measure event) that no two positive values of~$L_t$ coincide, each~$y\in\BbbL_n$ appears at most once in \eqref{E:3.11}. Moreover, for any~$y\in\BbbL_n$ with $f(\theta_n(y),\sqrt{L_t(y)}-m_n)\ne0$ that does not appear there exists~$x\in\BbbL_n$ with $L_t(x)\ge L_t(y)$ and $y\in B_{k_n}(x)\smallsetminus B_r(x)$. (Indeed, otherwise~$y$ would be a $k_n$-local maximum and would appear on the right.) 

These observations show that, if $m_n>\lambda$ and $k_n/2>r>\lambda$, then the probability in~\eqref{E:3.10} is bounded by
\begin{equation}
\label{E:3.15u}
P^\varrho\Bigl(\,\max_{x\in\BbbL_n}\sqrt{L_t(x)}-m_n>r\Bigr)+P^\varrho\Bigl(\exists x,y\in\Gamma(\lambda)\colon y\in B_{k_n}(x)\smallsetminus B_r(x)\Bigr).
\end{equation}
Our assumption that~$k_n\to\infty$ with $n-k_n \to \infty$
ensures that the above inequalities between~$m_n$, $k_n$, $r$ and~$\lambda$ are satisfied for~$n$ large and, by Lemmas~\ref{lemma-3.2a}--~\ref{lemma-3.3a}, that \eqref{E:3.15u} vanishes in the limit $n\to\infty$ followed by~$r\to\infty$. Hence we get \eqref{E:3.10} as desired. 
\end{proofsect}

\subsection{Two technical propositions}
We will now move to the proof of Theorem~\ref{thm-3.1}. As in~\cite{BL1,BL2,BL3}, different arguments are needed to describe the law of the position/value of $\sqrt{L_t}$ at its local maxima and the law of the ``shape'' of the configuration nearby. The former will conveniently be reduced to the arguments from the proof of convergence of the law of the maximum in~\cite{BL4} so most of our technical work will go towards controlling the ``shape.'' 

As discussed in Section~\ref{sec-2}, we will follow a Lindeberg-like approach based on swapping the increment of~$\sqrt{L_t}$ in the last~$k$ generations of the tree by that of a GFF, which is the Gaussian process~$h$ defined in \eqref{E:2.10a}. We will write $\wt P$ for the law of~$h$, use~$\wt E$ for the associated expectation  and let $E^\varrho\otimes\wt E$ denote the expectation for the product law~$P^\varrho\otimes\wt P$ under which $L_t$ and~$h$ are independent. Recall that~$\frakm(x)$ denotes, for~$x\ne\varrho$, the nearest vertex to~$x$ on the path to the root. The ``swapping'' argument now comes in:

\begin{proposition}[Swapping the local time for GFF]
\label{prop-3.1}
Given two naturals $n>k\ge1$, a real $t>0$ and a function $f\in\Cloc^+$, define
\begin{equation}
\label{E:3.19w}
U_{n,k}(f):=E^\varrho\biggl(\,\prod_{x\in\BbbL_n}\texte^{-f(\theta_n(x),\sqrt{L_t(x)}-m_n,\,\sqrt{L_t(x)}-\sqrt{L_t(x\cdot)}\,)\,1_{\MM_{n,k}(\sqrt{L_t},\,x)}}\biggr)
\end{equation}
and
\begin{equation}
\label{E:3.16}
V_{n,k}(f):=E^\varrho\otimes \wt E\biggl(\,\prod_{x\in\BbbL_n}\texte^{-f(\theta_{n-k}(z),\sqrt{L_t(z)}+h(x)-h(z)-m_n,\,h(x)-h(x\cdot)\,)\,1_{\MM_{n,k}(h,x)\cap \BB_{n,k}(x)}}\biggr),
\end{equation}
where we put $z:=\frakm^k(x)$ to reduce clutter and with the same convention set
 \begin{equation}
\label{E:3.19x}
\BB_{n,k}(x):=\Bigl\{m_k+k^{1/3}\le m_n-\sqrt{L_t(z)}\le m_k+k^{2/3}\Bigr\}.
\end{equation}
Then
\begin{equation}
\label{E:3.22w}
\lim_{k\to\infty}\limsup_{n\to\infty}\bigl|U_{n,k}(f)-V_{n,k}(f)\bigr|=0
\end{equation}
holds for all $t>0$ and all $f\in\Cloc^{+}$. Moreover, 
\begin{equation}
\label{E:3.19}
\lim_{\ell\to\infty}
\limsup_{k \to \infty} \limsup_{n \to \infty} 
P^\varrho\otimes\wt P\biggl(\,\sum_{x\in\BbbL_n}1_{\MM_{n,k}(h,x)\cap \BB_{n,k}(x)}1_{\{\sqrt{L_t(z)}+h(x)-h(z)>m_n-\lambda\}}>\ell\biggr)=0
\end{equation}
holds for all $t>0$ and all~$\lambda>0$. 
\end{proposition}

To see the relevance of the above for the problem at hand, note that 
\begin{equation}
%\label{}
U_{n,k}(f)=E^\varrho(\texte^{-\langle\zeta_{n,k}^{(t)},f\rangle}).
\end{equation}
The set $\BB_{n,k}(x)$ is in turn a kind of barrier event which are generally events that constrain the whole profile of~$k\mapsto L_t(x_k)$ for~$x_0=\varrho,x_1,\dots,x_n=x$ labeling the vertices on the unique path from the root to~$x$. Here we only impose a restriction on one value on this path. The powers $k^{1/3}$, resp., $k^{2/3}$ can be replaced by $k^\alpha$, resp., $k^{1-\alpha}$ for any~$0 < \alpha < 1/2$. 
The statement \eqref{E:3.19} shows that the number of terms non-trivially contributing to the product in \eqref{E:3.16} is tight.

\begin{comment}
In order to process the quantity $V_{n,k}(f)$ further, let~$\fraka\colon\N\to[0,\infty)$ be defined by
\begin{equation}
%\label{}
\fraka(x):=\slb \,\inf\{k\ge0\colon x\le b^k\}.
\end{equation}
It turns out that
\begin{equation}
%\label{}
\wh C(x,y):=\fraka(x)+\fraka(y)-\fraka\bigl(|y-x|\bigr)
\end{equation}
is a covariance kernel which thus defines a centered Gaussian process indexed by the naturals. We will demonstrate this by showing that this process arises as the $n\to\infty$ weak limit of the GFF on~$\BbbL_n$ conditioned on~$h_0=0$ for~$0\in\BbbL_n$ denoting the vertex represented by the sequence $(0,\dots,0)$. For this reason we call the resulting limit process the \emph{pinned GFF} and write~$\overline P$ for its law and~$\overline E$ for the associated expectation.
\end{comment}

Proposition~\ref{prop-3.1} effectively replaces the last~$k$ generations of the local time by  independent GFFs. In order to extract the local behavior from this, we will condition on the position of the local maxima of~$x\mapsto h(x)-h(\frakm^k(x))$ in subtrees of~$\T_n$ of depth~$k$ rooted at the vertices of~$\BbbL_{n-k}$ and apply the following asymptotic:

\begin{proposition}[``Shape'' of local extrema of GFF]
\label{prop-3.3}
There exists a probability law~$\nu$ on~$\R^\N$ such that the following holds:
Given~$f\in\Cloc^+$, $k\ge1$ and~$u\in\R$, let $g_{k,u}\colon[0,1]\times\R\to[0,\infty)$ be defined by
\begin{equation}
\label{E:3.22}
\texte^{-g_{k,u}(v,s)}:=\wt E\Bigl(\texte^{-f(v,\,s, \, h(\zero)-h(\zero\cdot))}\,\Big|\, h(\zero) = \max_{y\in\BbbL_k}h(y)=m_k+u\Bigr)
\end{equation}
and $g\colon[0,1]\times\R\to[0,\infty)$ defined by
\begin{equation}
\label{E:3.22a}
g(v,s):=-\log\biggl(
\int\nu(\textd\phi)\texte^{-f(v,s,\phi)}\biggr).
\end{equation}
Then
\begin{equation}
\label{E:3.24}
\sup_{v\in[0,1]}\,\,\sup_{s\in\R}\,\,\sup_{k^{ 1/13 }\le u\le k^{ 12/13 }}\,\bigl|g_{k,u}(v,s)-g(v,s)\bigr|\,\,\underset{k\to\infty}\longrightarrow\,\,0.
\end{equation}
The measure $\nu$ is concentrated on $\{\phi_0=0\}\cap\bigcap_{x\in\N\smallsetminus\{0\}}\{\phi_x>0\}$. 
\end{proposition}

We remark that the finite level set property for samples from~$\nu$ is not claimed here but will be proved as part of the proof of Theorem~\ref{thm-3.1}. 

\subsection{Integrating out the local component}
The proofs of Propositions~\ref{prop-3.1} and~\ref{prop-3.3} are technical and so we postpone their execution to Section~\ref{sec-4} and Section~\ref{sec-6}, respectively. The main role of these propositions is to integrate out the third component of the test function. We summarize this in:

\begin{lemma}
\label{lemma-3.4}
Assuming the statements of Propositions~\ref{prop-3.1}--\ref{prop-3.3}, let~$f\in\Cloc^+$ and let~$g$ be related to~$f$ as in \eqref{E:3.22a}. Then for all~$t>0$,
\begin{equation}
\label{E:3.25a}
\lim_{k\to\infty}\limsup_{n\to\infty}\,\bigl|U_{n,k}(f)-U_{n,k}(g)\bigr|=0.
\end{equation}
\end{lemma}

\begin{proofsect}{Proof}
Fix~$t>0$, let~$f\in\Cloc^+$ and define~$g$ as in \eqref{E:3.22a}. Given $n>k\ge1$ and $z\in\BbbL_{n-k}$, let~$\T_k(z)$ denote the subtree of~$\T_n$ rooted at~$z$ and abbreviate~$\BbbL_k(z):=\T_k(z)\cap\BbbL_n$. Abusing our earlier notation, let~$\scrF_k:=\sigma(h(y),L_t(y)\colon y\in\T_{n-k})$. 
As $\MM_{n,k}(h,x)$ depends only on~
$\{h(y)\colon y\in \BbbL_k (z)\}$ 
for~$z$ such that~$x\in \BbbL_k (z)$ while
\begin{equation}
%\label{}
\bigl\{h(x)-h(z)\colon x\in\T_k(z)\bigr\}_{z\in\BbbL_{n-k}}
\end{equation}
 are independent copies of the GFF in~$\T_k$, conditioning the expectation in the definition of~$V_{n,k}(f)$ on~$\scrF_k$ results in the product of conditional expectations 
\begin{equation}
\label{E:3.27}
E^\varrho\otimes \wt E\biggl(\,\prod_{x\in\BbbL_k(z)}\texte^{-f(\theta_{n-k}(z),\sqrt{L_t(z)}+h(x)-h(z)-m_n,\,h(x)-h(x\cdot)\,)\,1_{\MM_{n,k}(h,x)\cap \BB_{n,k}(x)}}\,\bigg|\,\scrF_k\biggr)
\end{equation}
over all~$z\in\BbbL_{n-k}$.  Next observe that, in light of the continuous nature of the law of~$h$, the event $\MM_{n,k}(h,x)$ occurs at exactly one~$x\in\BbbL_k(z)$ a.s.\ 
for each~$z\in\BbbL_{n-k}$. This equates \eqref{E:3.27} with the sum
over all $z \in \BbbL_{n-k}$ of 
\begin{equation}
\label{E:3.28}
\sum_{x\in\BbbL_k(z)}
E^\varrho\otimes \wt E\Bigl(\,1_{\MM_{n,k}(h,x)\cap \BB_{n,k}(x)}
\,\texte^{-f(\theta_{n-k}(z),\,s(x),\,h(x)-h(x\cdot)\,)\,}\,\Big|\,\scrF_k\Bigr),
\end{equation}
where we also abbreviated
\begin{equation}
%\label{}
s(x):=\sqrt{L_t(z)}-m_n+h(x)-h(z).
\end{equation}
Conditioning on~$h(x)-h(z)$, the symmetries of $\T_n$ give 
\begin{equation}
\begin{aligned}
E^\varrho\otimes \wt E\Bigl(\,1_{\MM_{n,k}(h,x)\cap \BB_{n,k}(x)}&
\,\texte^{-f(\theta_{n-k}(z),\,s(x),\,h(x)-h(x\cdot)\,)\,}\,\Big|\,\scrF_k\Bigr)
\\
&=E^\varrho\otimes \wt E\Bigl(\,1_{\MM_{n,k}(h,x)\cap \BB_{n,k}(x)}
\,\texte^{-g_{k,u(x)}(\theta_{n-k}(z),\,s(x))\,}\,\Big|\,\scrF_k\Bigr)
\end{aligned}
\end{equation}
where $g_{k,u}$ is as in \eqref{E:3.22} and where
\begin{equation}
%\label{}
u(x):=h(x)-h(z)-m_k.
\end{equation}
This now has the same form as the  expectation in \eqref{E:3.28} albeit with~$f$ replaced by $g_{k,u(x)}$ and so we readily conclude
\begin{equation}
\label{E:3.26}
V_{n,k}(f)=E^\varrho\otimes \wt E\biggl(\,\prod_{x\in\BbbL_n}\texte^{-g_{k,u(x)}(\theta_{n-k}(z),\,s(x))\,1_{\MM_{n,k}(h,x)\cap \BB_{n,k}(x)}}\biggr)
\end{equation}
by tracing back the above steps.

We will now argue that~$g_{k,u(x)}(\theta_{n-k}(z),s(x))$ can be replaced by~$g(\theta_{n-k}(z),s(x))$ without affecting the $n\to\infty$ and $k\to\infty$ limit. Indeed, the assumptions on~$f$ ensure existence of a~$\lambda>0$ be such that~$f(\cdot,h,\cdot)$ vanishes unless~$|h|\le\lambda$. Then 
$g_{k,u(x)}(\cdot,s(x))$ and~$g(\cdot,s(x))$ vanish unless~$|s(x)|\le\lambda$; i.e., unless 
$|\sqrt{L_t(z)}-m_n+m_k+u(x)|\le\lambda$.
Writing
\begin{equation}
%\label{}
\CalS(x):=\Bigl\{\texte^{-g_{k,u(x)}(\theta_{n-k}(z),\,s(x))\,1_{\MM_{n,k}(h,x)\cap \BB_{n,k}(x)}}\ne\texte^{-g(\theta_{n-k}(z),\,s(x))\,1_{\MM_{n,k}(h,x)\cap \BB_{n,k}(x)}}\Bigr\}
\end{equation}
the restrictions in $\BB_{n,k}(x)$ give 
\begin{equation}
%\label{}
-\lambda+k^{1/3}\le u(x)\le \lambda+ k^{2/3}\quad \text{ on }\CalS(x)
\end{equation}
for all~$x\in\BbbL_n$.
By Proposition~\ref{prop-3.3},  $g_{k,u(x)}(v,s(x))$ is thus uniformly close to~$g(v,s(x))$ whenever~$\CalS(x)$ occurs and~$k$ is so large that $[-\lambda+k^{1/3},\lambda+ k^{2/3}]
\subseteq[k^{ 1/13 },k^{ 12/13 }]$. 

In order to swap one function for the other in the exponent, we have to rule out that the total number of swapped terms explodes with~$n$ and~$k$. Here we observe that if $\CalS(x)$ occurs, the event
\begin{equation}
%\label{}
\MM_{n,k}(h,x)\cap \BB_{n,k}(x)\cap\bigl\{\sqrt{L_t(z)}+h(x)-h(z)\ge m_n-\lambda\bigr\}
\end{equation}
 must occur. Forcing the number of~$x$ where the latter occurs to be at most~$\ell$, the exchange of all $g_{k,u(x)}(v,s(x))$'s for~$g(v,s(x))$'s causes an error at most $\ell\epsilon_k$, where~$\epsilon_k$ is the quantity in \eqref{E:3.24}. Using the Intermediate-Value Theorem to take this error out of the exponential and the expectation then gives 
\begin{equation}
%\label{}
\begin{aligned}
\bigl|V_{n,k}(f)&-\texte^{O(\epsilon_k)\ell}V_{n,k}(g)\bigr|
\\
&\le2\texte^{\epsilon_k \ell}\,P^\varrho\otimes\wt P\biggl(\,\sum_{x\in\BbbL_n}1_{\MM_{n,k}(h,x)\cap \BB_{n,k}(x)}1_{\{\sqrt{L_t(z)}+h(x)-h(z)\ge m_n-\lambda\}}>\ell\biggr),
\end{aligned}
\end{equation}
where we abuse the notation by treating~$g$ as a function of three coordinates, again wrote $z:=\frakm^k(x)$ and where $O(\epsilon_k)$ denotes a deterministic quantity with values in~$[-\epsilon_k,\epsilon_k]$. Thanks to $V_{n,k}(g)\in[0,1]$, the term $\texte^{O(\epsilon_k)\ell}$ can be dropped at the cost of~$\texte^{\epsilon_k\ell}-1$ popping up on the right-hand side. 

Taking $n\to\infty$, $k\to\infty$ and $\ell\to\infty$ with the help of \eqref{E:3.19} now allows us to conclude that the statement holds with~$U_{n,k}$'s replaced by~$V_{n,k}$'s, i.e.,
\begin{equation}
%\label{E:3.25a}
\lim_{k\to\infty}\limsup_{n\to\infty}\,\bigl|V_{n,k}(f)-V_{n,k}(g)\bigr|=0.
\end{equation}
To complete the proof observe that, since $g\in\Cloc^+$, Proposition~\ref{prop-3.1} lets us exchange~$V_{n,k}$'s for $U_{n,k}$'s at no cost under these limits.
\end{proofsect}

Lemma~\ref{lemma-3.4} reduces the proof of Theorem~\ref{thm-3.1} to the convergence of the extremal process associated with local maxima. For that we will also need a slight upgrade of the convergence stated in Proposition~\ref{prop-3}: 

\begin{lemma}
\label{lemma-3.5a}
Let $a_n(t)$ be as in \eqref{E:1.24a} and set $\wt C_\star:=2c_\star\slb$ for~$c_\star$ as in Proposition~\ref{prop-3}. Then for all continuously differentiable $f\colon\R\to[0,\infty)$  with compact support, all natural $n\ge1$ and all real $t>0$ and $u>0$, the quantity $o(1)=o_{n,t,u}(1)$ defined implicitly by
\begin{equation}
\label{E:3.33}
\begin{aligned}
E^\varrho\Bigl(&\texte^{-f(\max_{y\in\BbbL_n}\sqrt{L_t(y)}-a_n(t)-\sqrt t-u)}\Bigr) 
\\
&\quad=\exp\biggl\{-\wt C_\star u \texte^{-2\slb\, u}\Bigl(o(1)+\int\textd h \,\texte^{-2\slb\, h}(1-\texte^{-f(h)})\Bigr)\biggr\}
\end{aligned}
\end{equation}
obeys
\begin{equation}
\label{E:3.29u}
\lim_{m\to\infty}\,\sup_{t,u\ge m}\,\limsup_{n\to\infty}\,\bigl|o_{n,t,u}(1)\bigr|=0.
\end{equation}
\end{lemma}

\begin{proofsect}{Proof}
Given $n\ge1$ and~$t>0$, abbreviate $M_{n,t}:=\max_{y\in\BbbL_n}\sqrt{L_t(y)}$. Assume  $f\colon\R\to[0,\infty)$ to be continuously differentiable with $\supp(f)\subseteq[-\lambda,\lambda]$, for some $\lambda>0$. Writing the expectation as a Stieltjes integral and integrating by parts yields
\begin{equation}
\label{E:3.30}
\begin{aligned}
E^\varrho\Bigl(&\texte^{-f(M_{n,t}-a_n(t)-\sqrt t-u)}\Bigr)
\\
&\qquad=1-\int \textd h\,\texte^{-f(h)}f'(h)\,P^\varrho\Bigl(M_{n,t}>a_n(t)+\sqrt t+u+h\Bigr).
\end{aligned}
\end{equation}
By Lemma~\ref{lemma-3.2a}, there exists~$c\in(0,\infty)$ such that for all~$n\ge1$ and all~$u>0$,
\begin{equation}
\label{E:3.31u}
P^\varrho\bigl(M_{n,t}>a_n(t)+\sqrt t+u\bigr)\le c(1+u)\texte^{-2\slb\, u}.
\end{equation}
Thanks to the restriction on the support of~$f$, the integral in \eqref{E:3.30} is for $u>\lambda$ dominated by the quantity
\begin{equation}
%\label{}
\delta(u):=c(1+u+\lambda)\texte^{-2\slb\, (u-\lambda)}\Vert \texte^{-f}f'\Vert_\infty 2\lambda.
\end{equation}
In particular, once~$u$ is so large that~$\delta(u)\le1/2$,  we get
\begin{equation}
\label{E:3.33u}
\begin{aligned}
\Biggl|\,\log E^\varrho\Bigl(&\texte^{-f(M_{n,t}-a_n(t)-\sqrt t-u)}\Bigr)
\\
&+\int \textd h\,\texte^{-f(h)}f'(h)\,P^\varrho\Bigl(M_{n,t}>a_n(t)+\sqrt t+u+h\Bigr)\biggr|\le \delta(u)^2
\end{aligned}
\end{equation}
relying on the inequality $|\log(1-x) +x|\le x^2$ for~$|x|\le 1/2$.

As is readily checked, $\delta(u)^2=o(u\texte^{-2\slb\, u})$ as~$u\to\infty$ and so for \eqref{E:3.33} it suffices to control the asymptotic of the integral in \eqref{E:3.33u}. For this we call upon the asymptotic stated in Proposition~\ref{prop-3} to get
\begin{equation}
\label{E:3.36u}
\begin{aligned}
\biggl|\int \textd h\,\texte^{-f(h)}f'(h)\,&P^\varrho\Bigl(M_{n,t}>a_n(t)+\sqrt t+u+h\Bigr)
\\
&-c_\star\int \textd h\,\texte^{-f(h)}f'(h)\, (u+h)\texte^{-2\slb\, (u+h)}\biggr|\le\delta(u)\epsilon_{n,t,u},
\end{aligned}
\end{equation}
where~$c_\star$ and $\epsilon_{n,t,u}$ are as in \eqref{E:3.34u} with, we note, $c_\star$ no larger than the constant~$c$ from~\eqref{E:3.31u}.
Since $\delta(u)=O(u\texte^{-2\slb u})$ as~$u\to\infty$, the error is $o(1)u\texte^{-2\slb u}$.  

It remains to find the asymptotic of the second integral in \eqref{E:3.36u}. Here one more integration by parts gives
\begin{equation}
\label{E:3.37u}
\begin{aligned}
\biggl|\int \textd h\,\texte^{-f(h)}&f'(h)\, (u+h)\texte^{-2\slb\, (u+h)}
\\
&\qquad\qquad-2\slb\, u\texte^{-2\slb\, u}
\int \textd h\, (1-\texte^{-f(h)})\, \texte^{-2\slb\, h}
\biggr|
\\
&\le 
\texte^{-2\slb\, u}
\int \textd h\,(1-\texte^{-f(h)})\, \bigl(2\slb\,|h|+1\bigr)\,\texte^{-2\slb\, h}.
\end{aligned}
\end{equation}
Setting $\tilde c:=(2\slb\,\lambda+1)\texte^{2\slb\, \lambda}2\lambda$, the right-hand side is bounded by $\tilde c\texte^{-2\slb\, u}$ which is $u^{-1}$-factor smaller than the order of the exponent in \eqref{E:3.33}. Putting \twoeqref{E:3.36u}{E:3.37u} together then yields the claim.
\end{proofsect}

\subsection{Limit of the structured extremal process}
We  are now ready to give a proof of the main result of this section:

\begin{proofsect}{Proof of Theorem~\ref{thm-3.1} from Propositions~\ref{prop-3.1}-\ref{prop-3.3}}
Let~$f\in\Cloc^+$ and assume that~$n>2k\ge1$ with~$k$ so large that~$(x,h,\phi)\mapsto f(x,h,\phi)$ depends only on $\{\phi_y\colon y=0,\dots,j\}$ for some $j< b^k$. Let~$\lambda>0$ be such that the support of~$f$ in these coordinates is fully contained in $[0,1]\times[-\lambda,\lambda]\times[-\lambda,\lambda]^{\{0,\dots,j\}}$. The proof comes in six steps.

\medskip\noindent
\textsl{Step 1: Reduction to first two coordinates.}
For~$n$ so large that $k_n>k$ we have
\begin{equation}
%\label{}
\bigl|E^\varrho\bigl(\texte^{-\langle\zeta_{n,k_n}^{(t)},f\rangle}\bigr)-U_{n,k}(f)\bigr|\le 
P^\varrho\Bigl(\exists x,y\in\Gamma(\lambda)\colon y\in B_{k_n}(x)\smallsetminus B_k(x)\Bigr),
\end{equation}
which by Lemma~\ref{lemma-3.3a} tends to zero as~$n\to\infty$ followed by~$k\to\infty$. Lemma~\ref{lemma-3.4} then permits us to swap~$U_{n,k}(f)$ for~$U_{n,k}(g)$ in these limits. In light of \eqref{E:3.22a}, for \eqref{E:3.4} it thus suffices to prove the existence of a random measure~$Z_t$ such that
\begin{equation}
\label{E:3.29a}
\lim_{k\to\infty}\limsup_{n\to\infty}\Biggl|\, U_{n,k}(g)-\E\biggl(\exp\Bigl\{-\int Z_t(\textd x)\otimes\texte^{-2\slb\, h}\textd h(1-\texte^{-g(x,h)})\Bigr\}\biggr)\Biggr|=0
\end{equation}
for any~$g\colon[0,1]\times\R\to[0,\infty)$ continuous with compact support. 

Approximating~$g$ by a compactly-supported, continuous~$\tilde g$ which is~$C^1$ in the second variable, for~$\lambda'>0$ such that~$\supp(g)\cup\supp(\tilde g)\subseteq[0,1]\times[-\lambda',\lambda']$ we have
\begin{equation}
%\label{}
\bigl|U_{n,k}(\tilde g) - U_{n,k}(g)\bigr|\le(\texte^{\Vert g-\tilde g\Vert_\infty \ell}-1)+P^\varrho\bigl(|\Gamma(\lambda')|>\ell\bigr).
\end{equation}
Thanks to Corollary~\ref{cor-3.4}, for~$\ell:=\Vert g-\tilde g\Vert_\infty^{-1/2}$ the \emph{limes superior} as~$n\to\infty$ of the left-hand side tends to zero with $\Vert g-\tilde g\Vert_\infty\to0$. A similar approximation argument applies inside the integral \eqref{E:3.29a}. It thus suffices to prove \eqref{E:3.29a} for~$g$ continuous with compact support that is~$C^1$ in the second variable. We will assume~$g$ to be such in what follows.

\medskip\noindent
\textsl{Step 2: Reduction to absolute maxima in subtrees.}
Write~$\T_{n-k}(z)$ for the subtree rooted at~$z\in\BbbL_k$, let~$M(z)$ denote the absolute maximum of~$\sqrt{L_t}$ on the leaves of~$\T_{n-k}(z)$ and, relying on the lexicographic ordering of~$\BbbL_n$, let $X(z)\in\BbbL_n\cap\T_{n-k}(z)$ be the minimal leaf-vertex where that maximum is achieved, i.e., $\sqrt{L_t(X(z))}=M(z)$. (The maximum can be degenerate when~$L_t$ vanishes on~$\BbbL_n\cap\T_{n-k}(z)$.) Observe $z\in\BbbL_k$ contributes to the product defining $U_{n,k}(g)$ only if $M(z)-m_n\in[-\lambda,\lambda]$. Moreover, the events $\MM_{n,k}(\sqrt{L_t},\cdot)$ force that~$X(z)$ is then the only point in~$\BbbL_n\cap\T_{n-k}(z)$ that contributes unless there exists another $k$-local maximum in 
$B_{n-k}(X(z))\smallsetminus B_k(X(z))$. 
Hence we get
\begin{equation}
\label{E:3.30a}
\begin{aligned}
\Biggl|\,U_{n,k}(g) - E^\varrho\biggl(\,\prod_{z\in\BbbL_k}&\texte^{-g(\theta_{n}(X(z)),M(z)-m_n)}\biggr)\Biggr|
\\
&\le P^\varrho\Bigl(\exists x,y\in\Gamma(\lambda)\colon y\in 
B_{n-k}(x)\smallsetminus B_k(x)\Bigr),
\end{aligned}
\end{equation}
where,  by Lemma~\ref{lemma-3.3a}, the right-hand side vanishes as~$n\to\infty$ and~$k\to\infty$.

Next observe that, by our assumptions on~$g$ we have
\begin{equation}
\label{E:3.53a}
\omega_g(r):=\sup_{s\in\R}\sup_{|v-v'|<r}\bigl|g(v,s)-g(v',s)\bigr|\,\,\underset{r\downarrow0}\longrightarrow\,\,0.
\end{equation}
Proceeding as in the proof of Lemma~\ref{lemma-3.4}, we curb the number of terms in the product that potentially contribute in order to prevent explosions of exponential factors and then use $|\theta_{n}(X(z))-\theta_k(z)|\le b^{-k}$ 
to swap $\theta_{n}(X(z))$ for~$\theta_k(z)$ with the result
 \begin{equation}
\label{E:3.31a}
\begin{aligned}
\Biggl|\,E^\varrho\biggl(\,\prod_{z\in\BbbL_k}\texte^{-g(\theta_{n}(X(z)),M(z)-m_n)}\biggr)-&E^\varrho\biggl(\,\prod_{z\in\BbbL_k}\texte^{-g(\theta_k(z),M(z)-m_n)}\biggr)\Biggr|
\\
&\le \texte^{\omega_g(b^{-k})\ell}-1+P^\varrho\bigl(|\Gamma(\lambda)|>\ell\bigr).
\end{aligned}
\end{equation}
This 
tends to zero as~$n\to\infty$, $k\to\infty$ and~$\ell\to\infty$ by \eqref{E:3.53a} and Corollary~\ref{cor-3.4}.

\medskip\noindent
\textsl{Step 3: Representation via a random measure.}
We have so far reduced the asymptotic computation of~$U_{n,k}(g)$ to the expectation of $\prod_{z\in\BbbL_k}\texte^{-g(\theta_k(z),M(z)-m_n)}$. We will address that by conditioning on~$\scrF_k:=\sigma(L_t(z)\colon z\in\BbbL_k)$ which results in the factorization
\begin{equation}
%\label{}
E^\varrho\Bigl(\prod_{z\in\BbbL_k}\texte^{-g(\theta_k(z),M(z)-m_n)}\,\Big|\,\scrF_k\Bigr) 
=\prod_{z\in\BbbL_k}E^\varrho\bigl(\texte^{-g(\theta_k(z),M(z)-m_n)}\,\big|\,\scrF_k\bigr)
\end{equation}
implied by the Markov property of the local time. This permits us to control the asymptotic of the expressions term by term.

Assume containment in the event
\begin{equation}
\label{E:3.56w}
\AA_k(t):=\Bigl\{\max_{z\in\BbbL_k}\sqrt{L_t(z)}\le m_k+\log\log k\Bigr\}.
\end{equation}
Denoting~$s:=L_t(z)$, the argument of~$g$ in the subtree rooted at~$z\in\BbbL_k$ then becomes
\begin{equation}
%\label{}
M(z) - m_n = M(z) - a_{n-k}(s\vee1)-\sqrt s - u,
\end{equation}
where
\begin{equation}
\label{E:3.50}
\begin{aligned}
u&:=m_n-a_{n-k}(s\vee1)-\sqrt s
\\
&=\slb\, k+\frac3{4\slb}\log\frac{n-k}n+\frac1{4\slb}\log\frac{n-k+\sqrt{s\vee1}}{n\sqrt{s\vee1}}-\sqrt s
\\
&=\slb\, k-\sqrt s - \frac1{4\slb}\log\sqrt{s\vee1}+O(k/n)+O(\sqrt s/n).
\end{aligned}
\end{equation}
For~$\sqrt s\le m_k+\log\log k$ as enforced by~$\AA_k(t)$, we get
\begin{equation}
%\label{}
u\ge \frac3{4\slb}\log k -\tilde c\log\log k +o(1)
\end{equation}
for some constant~$\tilde c>0$, 
implying that~$u$ is large uniformly once~$k$ is large (with~$2k<n$). We will now use this to extract the asymptotic of the conditional expectation in the terms where $L_t(z)$ is large and effectively bound it away for the other terms.

Suppose first that~$z\in\BbbL_k$ is such that $\AA_k(t)\cap\{L_t(z)\ge\log\log k\}$ occurs. For~$k$ large, Lemma~\ref{lemma-3.5a} shows
\begin{equation}
\label{E:3.58}
\begin{aligned}
E^\varrho\bigl(&\texte^{-g(\theta_k(z),M(z)-m_n)}\,\big|\,\scrF_k\bigr) 
\\
&\quad=\exp\biggl\{-\wt C_\star w_k(z)\Bigl(o(1)+\int\textd h \,\texte^{-2\slb\, h}(1-\texte^{-g(\theta_k(z),h)})\Bigr)\biggr\},
\end{aligned}
\end{equation}
where $o(1)\to0$ uniformly as~$k\to\infty$ and
\begin{equation}
\label{E:3.53}
\begin{aligned}
w_k(z): = b^{-2k}\biggl(\slb\, k-\sqrt{L_t(z)}-\frac1{8\slb}&\log(L_t(z)\vee 1)\biggr)^+ 
\\
&\times\,(L_t(z)\vee 1)^{1/4}\,\texte^{2\slb\sqrt{L_t(z)}}
\end{aligned}
\end{equation}
captures the relevant portion of the expression~$u\texte^{-2\slb\, u}$ for~$u$ as in \eqref{E:3.50}.

On $\AA_k(t)\cap\{L_t(z)<\log\log k\}$ the asymptotic \eqref{E:3.50} in turn gives $u-\lambda\ge\slb\, k-\log\log k$ once~$k$ is sufficiently large. We then  invoke the restriction on the support along with conditional Jensen's inequality and the uniform bound on the upper tail of the maximum from Lemma~\ref{lemma-3.2a} to get
\begin{equation}
\label{E:3.60}
\begin{aligned}
1\ge E^\varrho\bigl(&\texte^{-g(\theta_k(z),M(z)-m_n)}\,\big|\,\scrF_k\bigr) 
\\
&\ge \exp\bigl\{- E^\varrho(g(\theta_k(z),M(z)-m_n)\,\big|\,\scrF_k)\bigr\}
\\
&\ge \exp\bigl\{-\Vert g\Vert P^\varrho(M(z)\ge m_n-\lambda\,|\,\scrF_k)\bigr\}
\\
&\ge\exp\bigl\{-c k^2b^{-2k}\Vert g\Vert\bigr\},
\end{aligned}
\end{equation}
where in the last line we first used that~$(1+u-\lambda)\texte^{-2\slb(u-\lambda)}\le k^2 b^{-2k}$ once~$k$ is sufficiently large and~$c$ is the constant from \eqref{E:3.6x}. 

Note that the exponent in \eqref{E:3.60} is $O(k^2 b^{-k})$ even after summation over~$z\in\BbbL_k$. Since $w_k(z)\le k^2b^{-2k}$ when $L_t(z)<\log\log k$, the same applies even to the summation of the corresponding terms in the exponent in \eqref{E:3.58}. The terms with~$L_t(z)<\log\log k$ are thus negligible on both sides and so, denoting
\begin{equation}
%\label{}
Z_t^{(k)}:= \wt C_\star \sum_{z\in\BbbL_k}w_k(z)\delta_{\theta_k(z)},
\end{equation}
on~$\AA_k(t)$ we then get
\begin{equation}
%\label{E:3.33}
\begin{aligned}
%\biggl| 
E^\varrho\Bigl(&\prod_{z\in\BbbL_k}\texte^{-g(\theta_k(z),M(z)-m_n)}\,\Big|\,\scrF_k\Bigr) 
\\
&
=\exp\biggl\{o(1)\bigl(1+Z_t^{(k)}([0,1])\bigr)-\int \,Z_t^{(k)}(\textd x)\otimes\texte^{-2\slb\, h}\textd h\,
(1-\texte^{-g(x,h)})\Bigr)\biggr\},
\end{aligned}
\end{equation}
where~$o(1)\to0$ as~$k\to\infty$ uniformly on~$\AA_k(t)$.

\medskip\noindent
\textsl{Step 4: Proof of convergence.}
Note that the tightness of~$\max_{z\in\BbbL_k}\sqrt{L_t(z)}-m_k$ from Lemma~\ref{lemma-3.2a} gives
\begin{equation}
%\label{}
\lim_{k\to\infty} P^\varrho\bigl(\AA_k(t)\bigr)=1.
\end{equation}
Summarizing the above arguments, we thus get 
\begin{equation}
\label{E:3.66ii}
\lim_{k\to\infty}\limsup_{n\to\infty}\biggl|\,E^\varrho\bigl(\texte^{-\langle\zeta_{n,k_n}^{(t)},g\rangle}\bigr)-E^\varrho\Bigl(\texte^{-\int Z_t^{(k)}(\textd x)\otimes\texte^{-2\slb\, h}\textd h(1-\texte^{-g(x,h)})}\Bigr)\biggr|=0.
\end{equation}
But the first term does not depend on~$k$ while second term does not depend on~$n$, which is only possible if they both converge in their respective limits. 

Since \cite[Corollary~3.2 and Lemma~3.3]{BL4} gives tightness of $\{Z_t^{(k)}([0,1])\colon k\ge1\}$, we are permitted to extract a subsequential weak limit~$Z_t$ relative to the topology of vague convergence. The $k\to\infty$ limit of the second term in \eqref{E:3.66ii} is then realized by substituting~$Z_t$ for~$Z_t^{(k)}$. But the resulting quantity equals the limit of the first term in \eqref{E:3.66ii} and so is the same for all convergent subsequences. This means that the weak limit of $Z_t^{(k)}$ as~$k\to\infty$ actually exists and we get \eqref{E:3.29a} as desired.

Next we will prove that~$Z_t$ is given by the weak limit \eqref{E:1.18}. Here the proof of \cite[Lemma~3.3]{BL4} 
shows that the measure $Z_t^{(k)}$ above receives asymptotically vanishing contribution from $z\in\BbbL_k$ at which $\slb\, k-\sqrt{L_t(z)}<k^\delta$, for some fixed $\delta\in(0,1/4)$. This means that we can drop the third term in the positive part in \eqref{E:3.53} as well as the truncation by one in $(L_t(z)\vee 1)^{1/4}$ without affecting the convergence and/or the limit measure~$Z_t$. This proves \eqref{E:1.18}. The limit \eqref{E:1.13} concerns the total mass of~$Z_t$ and can thus be referred to~\cite[Theorem~1.2]{BL4}.

\medskip\noindent
\textsl{Step 5: Finite level sets in samples of~$\nu$.}
The next item to show is that
\begin{equation}
%\label{}
N_\lambda:=\bigl|\{j\in\N\colon \phi_j\le \lambda\}\bigr|
\end{equation}
 is finite for all $\lambda>0$ and $\nu$-a.e.~$\phi$.  Given~$a>0$ and $g\colon\R\to[0,1]$ continuous with~$g=1$ on~$[-1,1]$ and~$g=0$ outside~$(-2,2)$, let
 \begin{equation}
%\label{}
f(x,h,\phi):= a g(\lambda^{-1}h) \sum_{j\le r}g\bigl((2\lambda)^{-1}(h-\phi_j)\bigr).
\end{equation}
The properties of~$g$ imply that $f\in\Cloc^+$ and that $f(x,h,\phi)$ exceeds $a|\{j\le r\colon 0\le \phi_j\le\lambda\}|$ whenever $h\in[-\lambda,0]$. Noting also  that, for $k_n>r$,
\begin{equation}
%\label{}
\langle \zeta_{n,k_n}^t,f\rangle \le a|\Gamma(4\lambda)|,
\end{equation}
the limit result in \eqref{E:3.4} shows
\begin{equation}
\label{E:3.66}
\begin{aligned}
&\liminf_{n\to\infty}\,E^\varrho\bigl(\texte^{-a|\Gamma(4\lambda)|}\bigr) 
\\
&\quad\le \E\Bigg(\exp\biggl\{-Z_t([0,1])\int_{[-\lambda,0]} \!\!\textd h\,\texte^{-2\slb\,h}\int\nu(\textd\phi)
\bigl(1-\texte^{-a|\{j\le r\colon 0\le \phi_j\le\lambda\}|}
\bigr)\biggr\}\Biggr).
\end{aligned}
\end{equation}
Next observe that the left-hand side is independent of~$r$ and, by Corollary~\ref{cor-3.4}, tends to one in the limit as~$a\downarrow0$.
Taking $r\to\infty$ followed by $a\downarrow0$ with the help of the Bounded Convergence Theorem then gives
\begin{equation}
%\label{}
1\le\E\Bigl(\exp\bigl\{-c(\lambda)Z_t([0,1])\nu(N_\lambda=\infty)\bigr\}\Bigr)
\end{equation}
for $c(\lambda):=(2\slb)^{-1}(\texte^{2\slb\,\lambda}-1)$.
As~$Z_t([0,1])>0$ with positive probability, it follows that $\nu(N_\lambda=\infty)=0$ as desired.

\medskip\noindent
\textsl{Step 6: No atoms in~$Z_t$.}
It remains to check that~$Z_t$ is a.s.~diffuse. The intuitive argument is simple: If~$Z_t$ had an atom at some (random) point, say~$X$, then with positive probability the process of limiting extreme local maxima would have two points with non-negative second coordinates at~$X$. Rolling back the $n\to\infty$ limit we infer that,  for sufficiently large~$n$, there would have to be two local maxima of $\sqrt{L_t}$ that are more than~$k_n$ but less than, say, $n/2$ in graph-theoretical distance on~$\T_n$ from each other. This is impossible unless the event in \eqref{E:3.8a} occurs. 

A formal argument unfortunately requires work. We start by noting that, by linearity of $f\mapsto\langle\zeta_{n,k_n}^t,f\rangle$, the above shows that for all $f_1,\dots,f_N\in\Cloc^+$ and all $\lambda_1,\dots,\lambda_N\ge0$,
\begin{equation}
\label{E:3.71}
\lim_{n\to\infty}
E^\varrho\Bigl(\,\prod_{i=1}^N\texte^{-\lambda_k\langle \zeta^t_{n,k_n},f_i\rangle}\Bigr)
=\E\bigl(\,\texte^{-\Phi(\lambda_1 f_1+\dots+\lambda_N f_N)}\bigr)
\end{equation}
holds with
\begin{equation}
\label{E:3.70}
\Phi(f):=\int Z_t(\textd x)\otimes \texte^{-2\slb\,h}\textd h\otimes\nu(\textd\phi)\bigl(1-\texte^{-f(x,h,\phi)}\bigr).
\end{equation}
Next observe that, while~$Z_t$ may \emph{a priori} have atoms, the fact that~$Z_t(\R)<\infty$ $\BbbP$-a.s.\ implies that 
\begin{equation}
%\label{}
S:=\bigl\{x\in\R\colon \BbbP(Z_t(\{x\})>0)>0\bigr\}
\end{equation}
is at most countable. (Indeed,~$S$ is the set where $x\mapsto \E(Z_t((-\infty,x])1_{\{Z_t(\R)\le M\}})$ has a discontinuity for some~$M\in\N$.) Writing~$\II$ for the set of finite intervals with both endpoints in~$\R\smallsetminus S$, we claim that \eqref{E:3.71} applies even to functions $f_1,\dots,f_N$ of the form
\begin{equation}
\label{E:3.72}
f_i(x,h,\phi):=1_{I_i}(x)1_{(0,1)}(h),
\end{equation}
where~$I_1,\dots,I_N\in\II$. Indeed, let~$\hat f_1,\dots,\hat f_N$ be defined using the interiors and~$\tilde f_1,\dots,\tilde f_N$ using the closures of some~$I_1,\dots,I_N\in\II$, respectively, where we also use $1_{[0,1]}(h)$ in the second variable for~$\tilde f_i$'s. Note that then $\hat f_i\le f_i  \le  \tilde f_i$ for all $i=1,\dots,N$. Now use that each~$\hat f_i$ can be written as an increasing limit of functions from $\Cloc^+$ and~$\tilde f_i$ as a decreasing limit of functions from~$\Cloc^+$ to get
\begin{equation}
%\label{}
\begin{aligned}
\E\bigl(\,\texte^{-\Phi(\lambda_1 \tilde f_1+\dots+\lambda_N \tilde f_N)}&\bigr)
\le \liminf_{n\to\infty} E^\varrho\Bigl(\,\prod_{i=1}^N\texte^{-\lambda_k\langle \zeta^t_{n,k_n},f_i\rangle}\Bigr)
\\&\le \limsup_{n\to\infty}
E^\varrho\Bigl(\,\prod_{i=1}^N\texte^{-\lambda_k\langle \zeta^t_{n,k_n},f_i\rangle}\Bigr)
=\E\bigl(\,\texte^{-\Phi(\lambda_1 \hat f_1+\dots+\lambda_N \hat f_N)}\bigr).
\end{aligned}
\end{equation}
The a.s.-equality of~$\Phi(\lambda_1 f_1+\dots+\lambda_N f_N)$ to $\Phi(\lambda_1 \tilde f_1+\dots+\lambda_N \tilde f_N)$ implied by the restriction on the interval endpoints then proves  \eqref{E:3.71} for~$f_1,\dots,f_N$.

As a consequence of \eqref{E:3.72}, the Curtiss theorem implies joint convergence in law of random variables $\{\langle \zeta^t_{n,k_n},f_i\rangle\colon i=1,\dots,N\}$ for any~$f_1,\dots,f_N$ of the form \eqref{E:3.72} with intervals in~$\II$. If these functions have also disjoint supports, then
\begin{equation}
%\label{}
\Phi(\lambda_1 f_1+\dots+\lambda_N f_N) = \sum_{i=1}^N\alpha Z_t(I_i)(1-\texte^{-\lambda_i})
\end{equation}
with $\alpha:=(2\slb)^{-1}(1-\texte^{-2\slb})$ implying that, conditional on~$Z_t$, the limit law of the $N$-tuple $\{\langle \zeta^t_{n,k_n},f_i\rangle\colon i=1,\dots,N\}$ is Poisson with parameters $(\alpha Z_t(I_1),\dots,\alpha Z_t(I_N))$. The Portmanteau theorem then gives
\begin{equation}
\label{E:3.73}
\limsup_{n\to\infty}
P^\varrho\Bigl(\,\max_{i=1,\dots,N}\langle \zeta^t_{n,k_n},f_i\rangle\le1\Bigr)
\le\E\biggl(\,\prod_{i=1}^N\texte^{-\alpha Z_t(I_i)}\bigl(1+\alpha Z_t(I_i)\bigr)\biggr)
\end{equation}
whenever~$f_1,\dots,f_N$ are as in \eqref{E:3.72} for disjoint~$I_1,\dots,I_N\in\II$.

We are now ready to give a formal version of our intuitive argument. Denote by~$\GG_{n,k}$ the complement of the event in \eqref{E:3.8a} with~$\lambda:=0$ and, for~$k\ge1$ and $i\in\Z$, let~$I_i\in\II$ be an open subinterval of~$((i-1)b^{-k},ib^{-k})$. Define~$f_i$ by \eqref{E:3.72} using these intervals. Assuming~$\GG_{n,k}$ along with the full probability event that no two positive values of~$L_t$ coincide, $L_t$ has only one local maximum in $\theta_n^{-1}([(i-1)b^{-k},ib^{-k}))$ unless it vanishes there. For~$k_n\ge k$ this implies $\langle  \zeta^t_{n,k_n},f_i\rangle\le1$ a.s.\ on~$\GG_{n,k}$ for all~$i\in\Z$ and so
\begin{equation}
\label{E:3.74}
P^\varrho(\GG_{n,k})\le P^\varrho\Bigl(\,\max_{i\in\Z}\,\langle \zeta^t_{n,k_n},f_i\rangle\le1\Bigr),
\end{equation}
where the use of ``max'' is justified by noting that $\langle \zeta^t_{n,k_n},f_i\rangle\ne0$ for only a finite number~$i\in\Z$.
Using that $u\mapsto \texte^{-u}(1+u)$ is decreasing on~$\R_+$ we in turn get
\begin{equation}
\label{E:3.75}
\E\biggl(\,\prod_{i\in\Z}\texte^{-\alpha Z_t(I_i)}\bigl(1+\alpha Z_t(I_i)\bigr)\biggr)
\le 1-\bigl(1-\texte^{-\alpha\delta}(1+\alpha\delta)\bigr)\BbbP\bigl(\max_{i\in\Z}Z_t(I_i)>\delta\bigr)
\end{equation}
for any~$\delta>0$, where we again use that only a finite number of intervals can possibly contribute on each side. Invoking \eqref{E:3.73} and, relying on~$S$ being countable, we now increase each~$I_i$ to fill all of $((i-1)b^{-k},ib^{-k})$ to conclude
\begin{equation}
\label{E:3.81}
\BbbP\Bigl(\,\max_{i\in\Z} Z_t\bigl(((i-1)b^{-k},ib^{-k})\bigr)>\delta\Bigr)\le\frac{1-\limsup_{n\to\infty}P^\varrho(\GG_{n,k})}{1-\texte^{-\alpha\delta}(1+\alpha\delta)}
\end{equation}
for any~$k\ge1$ and $\delta>0$. But the symmetries of~$\T_n$ observed by~$L_t$ imply
\begin{equation}
%\label{}
\Bigl\{Z_t\bigl([ib^{-k},b^{-k-1}+ib^{-k})\bigr)\colon i\in\Z\Bigr\}\laweq
\Bigl\{Z_t\bigl([b^{-k-1}+ib^{-k},2b^{-k-1}+ib^{-k})\bigr)\colon i\in\Z\Bigr\}
\end{equation}
and the inclusions
\begin{equation}
%\label{}
\{ib^{-k}\}\subseteq[ib^{-k},b^{-k-1}+ib^{-k})
\end{equation}
and
\begin{equation}
\bigl[b^{-k-1}+ib^{-k},2b^{-k-1}+ib^{-k}\bigr)\subseteq\bigl(ib^{-k},(i+1)b^{-k}\bigr)
\end{equation}
give
\begin{equation}
\label{E:3.85}
\BbbP\Bigl(\,\max_{i\in\Z} Z_t\bigl(\{ib^{-k}\}\bigr)>\delta\Bigr)
\le \BbbP\Bigl(\,\max_{i\in\Z} Z_t\bigl(((i-1)b^{-k},ib^{-k})\bigr)>\delta\Bigr).
\end{equation}
Taking $k\to\infty$ in \eqref{E:3.81} and \eqref{E:3.85} with the help of Lemma~\ref{lemma-3.3a} rules out atoms of~$Z_t$ of size in excess of~$\delta$. Taking~$\delta\downarrow0$ then shows that~$Z_t$ has no atoms a.s.
\end{proofsect}

%%%%%%
%\newpage
\section{Swapping the local time for GFF}
\label{sec-4}
\noindent
We now move to the proof of Proposition~\ref{prop-3.1} dealing with the swap of the local time for GFF in the last couple of generations of the tree. This proposition has served as one of the key inputs for our proof of Theorem~\ref{thm-3.1} and thus also Theorem~\ref{thm-2}. The other key input, Proposition~\ref{prop-3.3}, will be proved in Section~\ref{sec-6}. 

\subsection{Useful lemmas}
We start by proving a sequence of useful lemmas. The first of these deals with the construction of a coupling that drives the rest of the argument.

\begin{lemma}[Coupling to GFF]
\label{lemma-4.1}
Fix~$n>k\ge1$ and~$t>0$. There exists a coupling of the local time $L_t$ and two GFFs~$h$ and~$\tilde h$ on~$\T_n$ such that
\begin{equation}
\label{E:3.29}
h\,\independent L_t\quad\wedge\quad\tilde h \,\,\independent\,\,\bigl\{L_t(z)\colon z\in\T_{n-k}\bigr\}
\end{equation}
and such that the following holds almost surely:
\begin{equation}
\label{E:4.2}
\forall x\in\T_{n-k}\colon\,\,\tilde h(x) = h(x)
\end{equation}
and
\begin{equation}
\label{E:3.25}
\forall x\in\T_n\smallsetminus\T_{n-k}\colon\,\,L_t(x)+\bigl[h(x)-h(z)\bigr]^2 = \Bigl(\tilde h(x)-\tilde h(z)+\sqrt{L_t(z)}\Bigr)^2,
\end{equation}
where we used the shorthand~$z:=\frakm^k(x)$
 to make \eqref{E:3.25} easier to parse.
\end{lemma}

\begin{proofsect}{Proof}
The idea of the proof is to apply Lemma~\ref{lemma-D} or, more precisely, arguments from its proof, separately in each~$\T_k(z)$ for~$z$ ranging through~$\BbbL_{n-k}$. Recall that~$\wt P$ denotes the law of GFF. Define $\frakp_k\colon \{-1,+1\}^{\T_k}\times[0,\infty)^{\T_k}\times[0,\infty)\to[0,1]$ to be the conditional probability mass function of the signs of $\{h_x+\sqrt s\colon x\in\T_k\}$ given  their  absolute values $\{|h_x+\sqrt s|\colon x\in\T_k\}$, i.e.,
\begin{equation}
%\label{}
\frakp_k\bigl(\sigma,y, s) := \wt P\biggl(\,\bigcap_{x\in\T_k}\bigl\{h_x =\sigma_x|h_x + \sqrt{s}| - \sqrt{s}
\Bigr\}\,\bigg|\, |h_x + \sqrt{s}| = y(x), x \in \T_k \biggr).
\end{equation}
Given a sample of $(L_t,h)$ from~$P^\varrho\otimes\wt P$, use these to sample $\sigma \in \{-1, +1 \}^{\mathbb{T}_n \backslash \mathbb{T}_{n-k}}$
with probability 
\begin{equation}
\label{E:4.5}
\prod_{z\in\BbbL_{n-k}}
\frakp_k\Bigl(\{\sigma_x\colon x\in\T_k(z)\},\bigl\{\textstyle\sqrt{L_t(x)+[h(x)-h(z)]^2}\colon x\in\T_k(z)\bigr\}, L_t(z) \Bigr).
\end{equation}
Now set~$\tilde h=h$ on~$\T_{n-k}$ as required by \eqref{E:4.2} and, for each~$z\in\BbbL_{n-k}$, let
\begin{equation}
%\label{}
\tilde h(x):=h(z)-\sqrt{L_t(z)}+\sigma_x \sqrt{L_t(x)+[h(x)-h(z)]^2},\qquad x\in\T_k(z).
\end{equation}
This identifies~$\tilde h$ on all of~$\T_n$ and thus defines a coupling of~$h$, $L_t$ and~$\tilde h$.

The definition ensures the validity of \twoeqref{E:4.2}{E:3.25} and~$h\independent L_t$ was assumed from the beginning. What remains to be shown is the second part of \eqref{E:3.29} and that~$\tilde h$ is a GFF. This is where Lemmas~\ref{lemma-M} and~\ref{lemma-D} come handy: Condition on
\begin{equation}
%\label{}
\scrF_k:=\sigma\bigl(\{L_t(z)\colon z\in\T_{n-k}\}\bigr)
\end{equation}
and note that, by the Markov property of~$L_t$, 
the conditional law of $\{L_t(x)\colon x\in\T_k(z)\}$, for $z\in\BbbL_{n-k}$, is that of $\{L_{u(z)}(x)\colon x\in\T_k\}$ for~$u(z):=L_t(z)$. In addition, the Markov property implies that $\{\{L_t(x)\colon x\in\T_k(z)\}\colon z\in\BbbL_{n-k}\}$ are independent conditionally on~$\scrF_k$. Noting also that $\{h(x)-h(z)\colon x\in\T_k(z)\}$ are independent samples of GFF in~$\T_k$ for each~$z$, Lemma~\ref{lemma-D} along with the product structure of \eqref{E:4.5} gives that, conditionally on $\scrF_k$, the fields $\{\{\tilde h(x)-\tilde h(z)\colon\T_k(z)\}\colon z\in\BbbL_{n-k}\}$, are independent copies of GFF. 

Thanks to \eqref{E:4.2} we now conclude the second part of \eqref{E:3.29}. The tree-indexed random walk structure of GFF combined with \eqref{E:4.2} in turn shows that~$\tilde h$ is a GFF in~$\T_n$.\end{proofsect}

We will write~$\overline P^\varrho$ for the coupling measure and prove the claim with~$\tilde h$ instead of~$h$. We will repeatedly use the following explicit version of \eqref{E:2.14}:
\begin{equation}
\label{E:4.8w}
\sqrt{L_t(x)}\le \sqrt{L_t(x)+[h(x)-h(z)]^2}\le\sqrt{L_t(x)}+\frac{[h(x)-h(z)]^2}{\sqrt{L_t(x)}},
\end{equation}
where~$z:=\frakm^k(x)$. Next we observe that we can always assume the containment in the ``barrier'' event~$\BB_{n,k}(x)$ whenever~$L_t(x)$ is large:

\begin{lemma}
\label{lemma-3.5}
For all~$\lambda\ge0$,
\begin{equation}
\label{E:3.5}
\lim_{k\to\infty}\limsup_{n\to\infty} P^\varrho\biggl(\,\bigcup_{x\in \BbbL_n}\bigl\{\sqrt{L_t(x)}\ge m_n-\lambda\bigr\}\cap\BB_{n, k}(x)^\cc\biggr) = 0,
\end{equation}
where~$\BB_{n, k}(x)$ is the barrier event from \eqref{E:3.19x}.
\end{lemma}

\begin{proofsect}{Proof}
The proof builds on a barrier estimate proved as part of the proof of~\cite[Proposition~3.1]{A18}, which we will cite from heavily. The argument there works with the local time of a Brownian motion on an associated metric tree, which is constructed by replacing the edges of~$\T_n$ by line-segments 
of length $1/2$. We will write~$d(\cdot,\cdot)$ for the associated metric and, abusing the notation both here and in~\cite{A18}, will keep writing $\T_n$ for the metric tree and~$L_t$ for the underlying local time there.

We start by recalling the event $G_u^n (t)$ from \cite[Eq.~(3.3)]{A18} with $\kappa:=\frac58(\log b)^{-1/2}$ and~$y$ replaced by~$u>0$. 
For each $z \in\mathbb{L}_{n-k}$ and $s \in [0, n-k]$,
let $z_s$ be the point on the unique path from $\varrho$ to $z$ with $d(\varrho, z_s) = s/2$.
Then, the probability on the left-hand side of (\ref{E:3.5})
is bounded from above by $P^{\varrho} (G_u^n (t))$ plus the sum over~$z\in \mathbb{L}_{n-k}$ of
\begin{equation}
\label{E:4.10w}
\begin{aligned}
P^{\varrho}&
\Biggl(\left\{\max_{x \in \mathbb{L}_k(z)} \sqrt{L_t (x)} \ge m_n-\lambda \right\}
\cap \left\{\sqrt{L_t (z)} -(m_n-m_k) \notin [-k^{2/3}, -k^{1/3}] \right\} \\
&\cap \left\{\sqrt{L_t (z_s)} \le \sqrt{t} + \frac{a_n(t)}{n} s + \kappa (\log (s \wedge (n-s)))^+ + u + 1,
\forall s \in [0, n-k] \right\}
 \Biggr),
 \end{aligned}
 \end{equation}
where the last event arises from the complement of~$G_u^n(t)$. Thanks to \cite[Lemma 3.2]{A18} we have
 $\lim_{u \to \infty} \limsup_{n \to \infty} P^{\varrho} (G_u^n (t)) = 0$ and so we only need to focus on \eqref{E:4.10w}. 
 
Abbreviate 
 \begin{equation}
%\label{}
\psi (s) := P^{\varrho} \left(\max_{x \in \mathbb{L}_k} \sqrt{L_s (x)} \ge m_n - \lambda \right).
\end{equation}
Since the second and third event in \eqref{E:4.10w} depend only on~$z_s$ for~$s\le n-k$ (note also that $z_{n-k}=z$), 
the Markov property (Lemma \ref{lemma-M}) allows us to condition on~$L_t(z)$ and rewrite the probability in~\eqref{E:4.10w} as 
\begin{equation}
\begin{aligned}
 E^{\varrho} \biggl[\psi \bigl(L_t(z)\bigr)
 &1_{\{\sqrt{L_t (z)} - (m_n-m_k) \notin [-k^{2/3}, -k^{1/3}],~L_t (z) > 0 \}} \\
 &\times 1_{\left\{\sqrt{L_t (z_s)} \le \sqrt{t} + \frac{a_n(t)}{n} s + \kappa (\log (s \wedge (n-s)))^+ + u + 1,
\forall s \in [0, n-k] \right\}} \biggr],
\end{aligned}
\end{equation}
where we also assumed that~$n$ is so large that~$m_n>\lambda$. Next we recall that $\{\sqrt{L_s}\colon s\ge0\}$ has the law of a zero-dimensional Bessel process (see Belius, Rosen and Zeitouni~\cite[Lemma 3.1(e)]{BRZ}). Using the well-known connection between this process and Brownian motion (see, e.g.,~\cite[Lemma 2.6]{BL4}) the above expectation is bounded by 
 \begin{equation}
 \label{E:4.13w}
 \begin{aligned}
\mathbb{E} \Biggl[
\sqrt{\frac{\sqrt{t}}{\sqrt{t} + B_{n-k}}} \psi \bigl((\sqrt{t} + B_{n-k})^2\bigr)
 &1_{\{\sqrt{t} + B_{n-k} - (m_n-m_k) \notin [-k^{2/3}, -k^{1/3}],~\sqrt{t} + B_{n-k} > 0 \}} \\
 &\times 1_{\left\{B_s \le \frac{a_n(t)}{n} s + \kappa (\log (s \wedge (n-s)))^+ + u + 1,
\forall s \in [0, n-k] \right\}} \Biggr],
\end{aligned}
\end{equation}
where the expectation is with respect to the law of Brownian motion $(B_s)_{s \ge 0}$ on $\mathbb{R}$ with $\E(B_s)=0$ and $\text{Var}(B_s) = s/2$ for $s \ge 0$. 

We will now estimate \eqref{E:4.13w} by methods of stochastic calculus. 
Let~$\widetilde{\mathbb{P}}$ be the probability measure defined by
\begin{equation}
%\label{}
\widetilde{\mathbb{P}} (A)
:= \mathbb{E} \left(1_A \texte^{\frac{2a_n(t)}{n} B_{n-k} - \frac{a_n(t)^2}{n^2} (n-k)} \right),
\quad A \in \sigma (B_s \colon s \le n-k).
\end{equation}
By the Girsanov theorem, under $\widetilde{\mathbb{P}}$, the process
\begin{equation}
%\label{}
\widetilde{B}_s := B_s - \frac{a_n(t)}{n} s, s \in [0, n-k]
\end{equation}
is a Brownian motion on $\mathbb{R}$ with $\text{Var}(B_s) = s/2$ for $s \ge 0$
started at $0$.
Then, for sufficiently large $n$ and $k$,
the expectation is bounded from above by
\begin{equation}
 \begin{aligned}
\widetilde{\mathbb{E}} \Biggl[&\texte^{-\frac{2a_n(t)}{n} \widetilde{B}_{n-k} - \frac{a_n(t)^2}{n^2} (n-k)}
\sqrt{\frac{\sqrt{t}}{\sqrt{t} + \frac{a_n(t)}{n} (n-k)
+ \widetilde{B}_{n-k}}} \,\,
\psi \biggl(\Bigl(\sqrt{t} + \frac{a_n(t)}{n} (n-k) + \widetilde{B}_{n-k}\Bigr)^2\biggr)
  \\
 &\qquad\qquad\times 1_{\{\widetilde{B}_{n-k} \in I_{n, k} \}}
 1_{\left\{\widetilde{B}_s \le \kappa \log (k+1) + \kappa (\log (s \wedge (n-k-s)))^+ + u + 1,
\forall s \in [0, n-k] \right\}} \Biggr],
\end{aligned}
\end{equation}
where
\begin{equation}
%\label{}
I_{n, k} := \left(-\sqrt{t} - \frac{a_n(t)}{n} (n-k), ~- \frac{k^{2/3}}{2} \right)
\cup \left(-k^{1/3}, ~\kappa \log k + u + 1 \right)
\end{equation}
and we have used the inequality
\begin{equation}
\label{E:4.18}
(\log (s \wedge (n-s)))^+
\le \log (k+1) + (\log (s \wedge (n-k-s)))^+,\quad s \in [0, n-k]
\end{equation}
for $n \gg k \gg 1$.
Note that for any $s \in I_{n, k}$, we have
\begin{equation}
 \begin{aligned}
m_n - \lambda 
&\ge \left(\sqrt{t} + \frac{a_n(t)}{n} (n-k) + s \right)
+ a_k \left(\left(\sqrt{t} + \frac{a_n(t)}{n} (n-k) + s \right)^2 \right)\\
&\qquad\quad+ \frac{3}{4\sqrt{\log b}} \log k - s - \lambda - \sqrt{t} + O\left(\frac{\log n}{n} k \right).
\end{aligned}
\end{equation}
The second line is positive for sufficiently large $k$
in light of $s \le \kappa \log k + u + 1$ and $\kappa < \frac34(\log b)^{-1/2}$.
Lemma~\ref{lemma-3.2a} implies that for any $s \in I_{n, k}$,
\begin{equation}
\label{E:4.20}
\begin{aligned}
\psi \biggl(\Bigl(\sqrt{t} + &\frac{a_n(t)}{n} (n-k) + s \Bigr)^2 \biggr)
\\
&\le c \Bigl(\frac{3}{4\sqrt{\log b}} \log k - s \Bigr)
k^{-3/2} \texte^{2s \sqrt{\log b}} \texte^{-c^{\prime} \frac{(\frac{3}{4\sqrt{\log b}} \log k - s)^2}{k}}
\end{aligned}
\end{equation}
for some constants $c,c^{\prime}>0$
that  depend on neither~$n$ nor~$k$.
This estimate along with~\cite[Lemma 2.4(i)]{A18} show that \eqref{E:4.10w} is
bounded from above by
\begin{equation}
\label{E:4.21ww}
 \begin{aligned}
c_1 \bigl(\kappa \log (k+1) + u + 1\bigr) k^{-3/2} 
&\int_{I_{n, k}} (\kappa \log (k+1) + u + 1 - s)
\left(\frac{3}{4\sqrt{\log b}} \log k - s \right)\\
&\times \sqrt{\frac{\sqrt{t} + n}{\sqrt{t} + \frac{a_n(t)}{n} (n-k) + s}}
\,\,\texte^{c_2 \frac{\log n}{n} s -c_3 \frac{(\frac{3}{4\sqrt{\log b}} \log k - s)^2}{k}} \textd s
\end{aligned}
\end{equation}
for some positive constants $c_1, c_2, c_3$
that depend on neither $n$ nor $k$.

It remains to bound the integral in \eqref{E:4.21ww}. Partitioning the integration domain into the intervals 
$(-\sqrt{t} - \frac{a_n(t)}{n} (n-k),~-\frac{1}{2} (\sqrt{t} + \frac{a_n(t)}{n} (n-k)))$,
$(-\frac{1}{2} (\sqrt{t} + \frac{a_n(t)}{n} (n-k)),~-\frac{k^{2/3}}{2})$,
and $(-k^{1/3}, \kappa \log k + u + 1)$
and estimating the resulting three integrals,
\eqref{E:4.21ww} is bounded from above by 
\begin{equation}
\label{E:4.22}
c_4 \log k~ k^{-3/2} \left(\int_{k^{2/3}/2}^{\infty} s^2 \texte^{-c_5 \frac{s^2}{k}} ds
+ n^{c_6} \texte^{-c_7 n} + k \right)
\end{equation}
for some positive constants $c_4, \dots, c_7$ that depend on neither~$n$ nor~$k$.
This quantity vanishes in the limits as $n \to \infty$ and $k \to \infty$. 
\end{proofsect}

Next we note that, whenever~$\BB_{n,k}(x)$ occurs, we can bound $h(x)-h(z)$ by a quantity proportional to~$k$:

\begin{lemma}
\label{lemma-3.6}
There exists~$\hat a>0$ such that
\begin{equation}
\label{E:3.6}
\lim_{k\to\infty}\limsup_{n\to\infty} P^\varrho\otimes \wt P\biggl(\,\bigcup_{ x\in\BbbL_n}\Bigl\{\bigl|h(x)-h(\frakm^k(x))\bigr|> \hat a k\Bigr\}\cap \BB_{n,k}(x)\biggr) = 0.
\end{equation}
\end{lemma}

\begin{proofsect}{Proof}
Set~$\hat a:=a+\slb$ for~$a>0$ to be determined and write~$z:=\frakm^k(x)$ whenever $x\in\BbbL_n$ is clear from context. By symmetry of the Gaussian distribution, the probability in the statement is bounded by twice the probability without the absolute value around the term $h(x)-h(\frakm^k(x))$. It suffices to prove the claim without the absolute value.

Note that on~$\BB_{n,k}(x)$ we have $\sqrt{L_t(z)}\ge m_n - m_k - k^{2/3}$ and so $h(x)-h(z)>(a+\slb)k$ along with $m_k\le\slb\,k$ show
\begin{equation}
%\label{}
h(x)-h(z)+\sqrt{L_t(z)}> (a+\slb)k+m_n-m_k - k^{2/3} \ge m_n+ (a/2) k.
\end{equation}
Interpreting this under the coupling measure~$\overline P^\varrho$, the independence in \eqref{E:3.29} allows us to swap~$h$ for~$\tilde h$ conditional on~$\{L_t(z')\colon z'\in\T_{n-k}\}$. This bounds the probability in \eqref{E:3.6} (without absolute value) by
\begin{equation}
\label{E:4.12w}
\begin{aligned}
\overline P^\varrho\Bigl(\exists x\in\BbbL_n\colon &\tilde h(x)-\tilde h(z)+\sqrt{L_t(z)}> m_n+
(a/2)k\Bigr)
\\
&\le
\overline P^\varrho\Bigl(\exists x\in\BbbL_n\colon \sqrt{L_t(x)+[h(x)-h(z)]^2}> m_n+
(a/2)k\Bigr)
\end{aligned}
\end{equation}
where the inequality follows from the coupling identity \eqref{E:3.25}. 

Next notice that $|h(x)-h(z)|\le [1+\log b]^{1/2}\sqrt{k}\sqrt n$ and $L_t(x)+ |h(x)-h(z)|^2> m_n^2$ force $L_t(x)>\frac12m_n^2$ whenever $k\le[2(1+\log b)]^{-1}m_n^2/n$. Under these circumstances the inequality in the event on the right of \eqref{E:4.12w} along with \eqref{E:4.8w} give
\begin{equation}
%\label{}
\sqrt{L_t(x)}\ge m_n+(a/2)k - 2\frac{1+\log b}{\slb}\, k
\end{equation}
whenever~$n$ is so large that $m_n\ge\frac1{\sqrt2}\slb\, n$.
For $a:=4\frac{1+\log b}{\slb}+2$, the probability in~\eqref{E:4.12w} is thus bounded by
\begin{equation}
\begin{aligned}
\wt P\Bigl(\,\max_{x\in\BbbL_n}\,\bigl|h(x)-h(z)\bigr|&>[1+\log b]^{1/2}\sqrt{k}\sqrt n\Bigr)
+P^\varrho\Bigl(\,\max_{x\in\BbbL_n}\sqrt{L_t(x)}\ge m_n+k\Bigr).
\end{aligned}
\end{equation}
The union bound along with a standard Gaussian tail estimate dominate the first probability by $b^n\texte^{-(1+\log b)n}=\texte^{-n}$, which tends to zero as~$n\to\infty$. The second probability vanishes as~$n\to\infty$ and~$k\to\infty$ by Lemma~\ref{lemma-3.2a}.
\end{proofsect}

Finally, we observe some consequences of the previous proof for the objects entering the coupling identity \eqref{E:3.25}:

\begin{lemma}
\label{lemma-3.7}
We have
\begin{equation}
\label{E:4.15w}
\lim_{k\to\infty}\limsup_{n\to\infty} P^\varrho\biggl(\,\bigcup_{x\in\BbbL_n}\Bigl\{\sqrt{L_t(x)}\le\frac12 m_n\Bigr\}\cap \BB_{n,k}(x)\biggr)=0
\end{equation}
and, recalling that $ \overline P^\varrho$ denotes the coupling measure from Lemma~\ref{lemma-4.1}, 
\begin{equation}
\label{E:4.16w}
\lim_{k\to\infty}\limsup_{n\to\infty} \overline P^\varrho\biggl(\,\bigcup_{x\in\BbbL_n}\bigl\{\tilde h(x)- \tilde h(z)+\sqrt{L_t(z)}<0\bigr\}\cap\BB_{n,k}(x)\biggr)=0,
\end{equation}
where we again invoked the shorthand $z:=\frakm^k(x)$.
\end{lemma}

\begin{proofsect}{Proof}
Interpreting \eqref{E:3.6} using the coupling measure,  \eqref{E:3.29} shows that the statement of \eqref{E:3.6} holds with~$\tilde h$ in place of~$h$. Hence we may assume that $|h(x)- h(z)|\le\hat a k$ and $|\tilde h(x)-\tilde h(z)|\le\hat a k$ whenever $\BB_{n,k}(x)$ occurs. But then 
$\sqrt{L_t(z)}\ge m_n - m_k - k^{2/3}$ imposed by~$\BB_{n,k}(x)$ forces $\tilde h(x)-\tilde h(z)+\sqrt{L_t(z)}\ge\frac34 m_n$ once~$n\gg k$, proving \eqref{E:4.16w}. 
The identity \eqref{E:3.25} turns that bound into $L_t(x)+[h(x)-h(z)]^2\ge\frac9{16}m_n^2$ which then gives~$\sqrt{L_t(x)}>\frac12 m_n$ once~$n\gg k$, proving \eqref{E:4.15w} as well.
\end{proofsect}

\subsection{Key lemma and proof of Proposition~\ref{prop-3.1}}
The proof of Proposition~\ref{prop-3.1} requires one additional lemma that we will state and prove next. First, given~$n\ge k>1$ and $\lambda>0$, let~$\GG_{n,k}^\lambda$ denote the 
intersection of the complements of the events in \eqref{E:3.5}, \eqref{E:4.15w} and \eqref{E:4.16w} along with the complement of the event in \eqref{E:3.6} for both~$h$ and~$\tilde h$. The above lemmas then show
\begin{equation}
\label{E:4.17w}
\lim_{k\to\infty}\liminf_{n\to\infty}\, \overline P^\varrho\bigl(\GG_{n,k}^\lambda\bigr)=1.
\end{equation}
Next let $\EE_{n,k}^\lambda$ be the event that the set
\begin{equation}
\label{E:4.18w}
\Bigl\{x\in\BbbL_n\colon \MM_{n, k}(\sqrt{L_t},x)\cap\BB_{n,k}(x)\text{ occurs }\wedge\,\,\sqrt{L_t(x)+[h(x)-h(z)]^2}\ge m_n-\lambda\Bigr\}
\end{equation}
equals the set
\begin{equation}
\label{E:4.19w}
\Bigl\{x\in\BbbL_n\colon \MM_{n, k}(\tilde h,x)\cap\BB_{n,k}(x)\text{ occurs }\wedge\,\tilde h(x)-\tilde h(z)+\sqrt{L_t(z)}\ge m_n-\lambda\Bigr\},
\end{equation}
where, as before, we are using the shorthand~$z:=\frakm^k(x)$ throughout. This event addresses the most difficult part of the ``swap'' of the local time for GFF because being a local maximum is generally not preserved by perturbations, however small they may be. Notwithstanding, we still get:

\begin{lemma}
\label{lemma-4.5}
For all $\lambda>0$,
\begin{equation}
%\label{}
\lim_{k\to\infty}\liminf_{n\to\infty} \,\overline P^\varrho\bigl(\EE_{n,k}^\lambda\bigr)=1.
\end{equation}
\end{lemma}

\begin{proofsect}{Proof}
For each~$x\in\BbbL_n$, let $\text{gap}(\tilde h,x)$ be the difference between the largest and second largest value of~$\tilde h$ in~$B_k(x)$. If~$\GG_{n,k}^{\lambda}\cap\BB_{n,k}(x)$ occurs, then the inequality \eqref{E:4.8w} gives
\begin{equation}
\label{E:4.21w}
\begin{aligned}
\Bigl|\sqrt{L_t(x)}&-\sqrt{L_t(z)} - \bigl(\tilde h(x)-\tilde h(z)\bigr)\Bigr|
\\
&=\Bigl|\sqrt{L_t(x)}-\sqrt{L_t(x)+[h(x)-h(z)]^2}\Bigr|
\le\frac{[h(x)-h(z)]^2}{\sqrt{L_t(x)}}\le 2\hat a^2\frac{k^2}{m_n}
\end{aligned}
\end{equation}
for $z:=\frakm^k(x)$. Thanks to the coupling identity \eqref{E:3.25}, the conditions for the fields on the right of \eqref{E:4.18w} and \eqref{E:4.19w} are identical so the only difference is the status of the events $\MM_{n,k}(\sqrt{L_t},x)$ and $\MM_{n,k}(\tilde h,x)$.
A routine use of the triangle inequality then shows that, on $\GG_{n,k}^\lambda\smallsetminus \EE_{n,k}^\lambda$, there exists~$x$ --- namely, one that lies in one but not both of the sets \eqref{E:4.18w} and \eqref{E:4.19w} --- for which $\BB_{n,k}(x)$ occurs and yet 
$\text{gap}(\tilde h,x)\le 4\hat a^2 k^2/m_n$. Hence we get
\begin{equation}
\label{E:4.22w}
\overline P^\varrho\bigl(\GG_{n,k}^\lambda\smallsetminus \EE_{n,k}^\lambda\bigr)\le \overline P^\varrho\Biggl(\,\bigcup_{x\in\BbbL_n}\biggl\{\text{gap}(\tilde h,x)\le 4\hat a^2\frac{k^2}{m_n}\biggr\}\cap \BB_{n,k}(x)\Biggr).
\end{equation}
Since~$\BB_{n,k}(x)$ is determined by~$L_t(z)$, it is independent of~$\text{gap}(\tilde h,x)$ by \eqref{E:3.29}. Using~$N_{n,k}$ to denote the number of~$z\in\BbbL_{n-k}$ such that~$\BB_{n,k}(x)$ occurs for some (and thus all) $x\in\BbbL_n\cap\T_k(z)$ and denoting
\begin{equation}
%\label{}
F_k(u):=\overline P^\varrho\Bigl(\text{gap}(\tilde h,x)\le u\Bigr),
\end{equation}
which by the symmetries of the tree does not depend on~$x$ or~$n$, the probability on the right of \eqref{E:4.22w} equals
\begin{equation}
\label{E:4.24w}
1-E^\varrho\Biggl( \biggl[1-F_{k}\Bigl(4\hat a^2\frac{k^2}{m_n}\Bigr)\biggr]^{N_{n,k}}\Biggr).
\end{equation}
As $m_n-m_k = m_{n-k}+O(\log k)$, Corollary~\ref{cor-3.4}
shows that $\{N_{n,k}\colon n\ge1\}$ is tight for each $k\ge1$. The continuity of the law of GFF in turn implies that $F_k(u)\to0$ as~$u\downarrow0$. It follows that \eqref{E:4.24w} tends to zero as~$n\to\infty$. The proof is completed by invoking \eqref{E:4.17w}.
\end{proofsect}

We are now ready for:

\begin{proofsect}{Proof of Proposition~\ref{prop-3.1}}
Let~$f\in\Cloc^+$ and assume $n>k\gg1$. Let~$\lambda>0$ be such that~$(x,h,\phi)\mapsto f(x,h,\phi)$ depends only on $\{\phi_y\colon y=0,\dots,j\}$ for some~$j<b^k$ and, when restricted to these coordinates, is supported in~$[0,1]\times[-\lambda,\lambda]\times[-\lambda,\lambda]^{\{0,\dots,j\}}$. Let us write $\text{osc}_f(r)$ for the supremum of $|f(x,t,\phi)-f(x',t',\phi')|$ over all triplets of coordinates satisfying $|x-x'|<r$, $|t-t'|<r$ and~$\Vert \phi-\phi'\Vert_\infty<r$.

Thanks to \eqref{E:4.17w} and Lemma~\ref{lemma-4.5}, the event
\begin{equation}
%\label{}
\HH_{n,k}:=\GG_{n,k}^\lambda\cap \EE_{n,k}^{\lambda}
\end{equation}
occurs with probability tending to one as~$n\to\infty$ and~$k\to\infty$. Suppose that~$\HH_{n,k}$ occurs along with the event 
\begin{equation}
%\label{}
\JJ_{n,k}^q:=\biggr\{\sum_{x\in\BbbL_n}1_{\MM_{n,k}(\tilde h,x)\cap\BB_{n,k}(x)}1_{\{\tilde h(x)-\tilde h(z)+\sqrt{L_t(z)}\ge m_n-\lambda\}}\le q\biggr\}
\end{equation}
for some~$q\ge1$.
Noting that~$f$ in the term corresponding to~$x$ in $U_{n,k}(f)$ vanishes unless~$\sqrt{L_t(x)}\ge m_n-\lambda$, on $\HH_{n,k}\cap  \JJ_{n,k}^q$ we  can insert the indicator of~$\BB_{n,k}(x)$ in front of that term without changing the result. Then we use \eqref{E:4.21w} along with $|\theta_n(x)-\theta_{n-k}(z)|\le\frac b{b-1}b^{-(n-k)}$ to replace the arguments of~$f$ in~$U_{n,k}(f)$ by those in~$V_{n,k}(f)$ causing an error at most $q\,\text{osc}_f(r_{n,k})$ in the exponent, where
\begin{equation}
%\label{}
r_{n,k}:=\max\Bigl\{\frac b{b-1}b^{-(n-k)},\,4\hat a^2\frac{k^2}{m_n}
 \Bigr\}.
\end{equation}
Since~$f$ in the term corresponding to~$x$ in~$U_{n,k}(f)$ vanishes unless $\sqrt{L_t(x)}\ge m_n-\lambda$ while~$f$ in the term corresponding to~$x$ in~$V_{n,k}(f)$ vanishes unless $\tilde h(x)-\tilde h(z)
+\sqrt{L_t(z)}\ge m_n-\lambda$, the containment in~$\EE_{n,k}^{\lambda}$ permits us to swap the event $\MM_{n,k}(\sqrt{L_t},x)$ for~$\MM_{n,k}(\tilde h,x)$ at all~$x\in\BbbL_n$. 

As a result of these manipulations, we get
\begin{equation}
%\label{}
\bigl|U_{n,k}(f)-V_{n,k}(f)\bigr|\le 2\bigl[1-\overline P^\varrho(\HH_{n,k}\cap\JJ_{n,k}^q)\bigr]+\texte^{q\,\text{osc}_f(r_{n,k})}-1.
\end{equation}
The continuity of~$f$ now ensures that $\text{osc}_f(r_{n,k})\to0$ as~$n\to\infty$, regardless of~$k\ge1$. As to the event~$\HH_{n,k}\cap\JJ_{n,k}^q$, here the bound \eqref{E:4.21w} gives
\begin{equation}
%\label{}
\HH_{n,k}\smallsetminus\JJ_{n,k}^q\subseteq\bigl\{|\Gamma(\lambda+1)|>q\bigr\}
\end{equation}
as soon as $2\hat a^2\frac{k^2}{m_n}\le 1$. Corollary~\ref{cor-3.4} shows that the probability of the event on the right tends to zero as~$n\to\infty$ and~$q\to\infty$. As $\overline P^\varrho(\HH_{n,k})\to1$ in the limits~$n\to\infty$ and~$k\to\infty$, we get both \eqref{E:3.22w} and \eqref{E:3.19}.
\end{proofsect}

%%%%%
%\newpage
\section{Cluster process and walk started from a leaf}
\label{sec-5}
\noindent
We will now address the proof of Corollary~\ref{cor-1.5}, which links the cluster process of the local time to that of GFF/BRW, and that of Theorem~\ref{thm-1}, which deals with the random walk started from the leaves. As part of the latter, we also prove Corollary~\ref{cor-1.4} which identifies the law of the position of the local-time maximizer.

\subsection{Connection to BRW cluster process}
Our proof of Corollary~\ref{cor-1.5} follows that of Theorem~\ref{thm-3.1} with the local time replaced by the~BRW/GFF. A key point is that Proposition~\ref{prop-3.3}, which is where the cluster process is extracted for the local time, plugs in seamlessly for the GFF as well.

We start by recalling some needed facts about the extremal properties of BRW/GFF. Let us write~$M_n$ for the maximum of the GFF on~$\BbbL_n$. Then, as shown in Addario-Berry and Reed~\cite[Theorem~3]{ABR} for fairly general BRW (and for the particular case at hand in~\cite[Lecture~7]{B-notes}), 
\begin{equation}
\label{E:4.73}
\wt P\bigl(|M_n-\wt m_n|>\lambda\bigr)\le c_1\texte^{-c_2\lambda}
\end{equation}
holds for all~$n\ge1$ and all~$\lambda>0$, where
\begin{equation}
%\label{}
\wt m_n:=\slb\, n-\frac{3}{4\slb}\log n,
\end{equation}
and~$c_1$ and~$c_2$ are positive constants. For the extremal level set,
\begin{equation}
\wt\Gamma(\lambda):=\bigl\{x\in\BbbL_n\colon h(x)\ge\wt m_n-\lambda\bigr\},
\end{equation}
Mallein~\cite[Theorem 4.5]{M-review} proves the clustering property
\begin{equation}
\label{E:4.75}
\lim_{k\to\infty}\,\limsup_{n\to\infty}\,\wt P\Bigl(\exists x,y\in\wt\Gamma(\lambda)\colon y\in B_{n-k}(x)\smallsetminus B_k(x)\Bigr) = 0.
\end{equation}
This also readily implies tightness
\begin{equation}
\label{E:4.76}
\lim_{\ell\to\infty}\limsup_{n\to\infty}\wt P\Bigl(\bigl|\wt\Gamma(\lambda)\bigr|>\ell\Bigr)=0
\end{equation}
of the cardinality of~$\wt\Gamma(\lambda)$.
Note that the above are the analogues of Lemmas~\ref{lemma-3.2a}-\ref{lemma-3.3a} and Corollary~\ref{cor-3.4} for the case at hand. 

The extraction of the limit process will also require the asymptotic statement
\begin{equation}
\label{E:4.77}
\wt P\bigl(M_n-\wt m_n>u\bigr)=\tilde c_\star\bigl(1+\tilde\epsilon_{n,u}
\bigr)\,u\texte^{-2u\slb},
\end{equation}
where~$\tilde c_\star$ is a positive constant and $\tilde\epsilon_{n,u}$ defined by this expression obeys
\begin{equation}
\label{E:4.78a}
\lim_{u\to\infty}\limsup_{n\to\infty}\,\bigl|\,\tilde\epsilon_{n,u}\bigr|=0.
\end{equation}
This statement, which plays the role of Proposition~\ref{prop-3} for the case at hand, can be found in Bramson, Ding and Zeitouni~\cite[Proposition 3.1]{BDZ-BRW}.

In addition to the above, we will need an analogue of Lemma~\ref{lemma-3.5}:

\begin{lemma}
\label{lemma-4.7}
For all~$\lambda\ge0$,
\begin{equation}
\label{E:4.79}
\lim_{k\to\infty}\limsup_{n\to\infty} \wt P\biggl(\,\bigcup_{x\in \BbbL_n}\bigl\{h(x)\ge \wt m_n-\lambda\bigr\}\cap\wt\BB_{n, k}(x)^\cc\biggr) = 0,
\end{equation}
where for $\tilde\sigma:=\frac1{12}$, 
\begin{equation}
\label{E:4.78}
\wt \BB_{n,k}(x):=\Bigl\{\wt m_k+k^{\tilde\sigma}\le \wt m_n-h(z)\le \wt m_k+k^{1-\tilde\sigma}\Bigr\}
\end{equation}
with our convention $z:=\frakm^k(x)$.
\end{lemma}

The proof of this lemma requires ideas similar to those entering the proof of Proposition~\ref{prop-3.3} and so we relegate it to Section~\ref{sec-6}. Taking this lemma for granted in what follows, we are now ready for:

\begin{proofsect}{Proof of Corollary~\ref{cor-1.5} from Proposition~\ref{prop-3.3} and Lemma~\ref{lemma-4.7}}
We follow the proof of Theorem~\ref{thm-3.1} and its reduction to Theorem~\ref{thm-2}. We start by the definition of the structured extremal process associated with the BRW/GFF. This is a random measure on $[0,1]\times\R\times\R^\N$ defined by
\begin{equation}
%\label{}
\wt\eta_{n,k} :=\sum_{x\in\BbbL_n}1_{\MM_{n,k}(h,\,x)}\,\delta_{\theta_n(x)}\otimes\delta_{h(x)-\wt m_n}\otimes\delta_{h(x)-h(x\cdot)}.
\end{equation}
Given~$f\in\Cloc^+$,
set
\begin{equation}
%\label{E:3.16}
\wt V_{n,k}(f):=\wt E\biggl(\,\prod_{x\in\BbbL_n}\texte^{-f(\theta_{n-k}(z),\,
h(x)-\wt m_n,\,h(x)-h(x\cdot)\,)\,1_{\MM_{n,k}(h,x)\cap\wt\BB_{n,k}(x)}}\biggr),
\end{equation}
where we used our standard shorthand $z:=\frakm^k(x)$. Lemma~\ref{lemma-4.7} along with uniform continuity of~$f$ then show
\begin{equation}
%\label{}
\lim_{k\to\infty}\limsup_{n\to\infty}\bigl|\wt E(\texte^{-\langle\wt\eta_{n,k},f\rangle})-\wt V_{n,k}(f)\bigr|=0.
\end{equation}
We thus need to compute the limit of $\wt V_{n,k}(f)$ as~$n\to\infty$ followed by~$k\to\infty$.

Assume that~$n>2k\ge1$ with~$k$ so large that~$(x,h,\phi)\mapsto f(x,h,\phi)$ depends only on $\{\phi_y\colon y=0,\dots,j\}$ for some $j< b^k$. Then, similarly as in \eqref{E:3.26}, the Markov property of the BRW/GFF along with a simple calculation show
\begin{equation}
\label{E:4.84}
\wt V_{n,k}(f)=\wt E\biggl(\,\prod_{x\in\BbbL_n}\texte^{-g_{k,\tilde u(x)}(\theta_{n-k}(z),\,h(x)-\wt m_n)\,1_{\MM_{n,k}(h,x)\cap\wt\BB_{n,k}(x)}}\biggr),
\end{equation}
where
\begin{equation}
%\label{}
\tilde u(x):=h(x)-\wt m_n +\wt m_n-m_k-h(z).
\end{equation}
For~$\lambda>0$ such that $\text{supp}(f)$ is contained in $[0,1]\times[-\lambda,\lambda]\times[-\lambda,\lambda]^{\{0,\dots,j\}}$, the term corresponding to~$x$ in \eqref{E:4.84} is non-trivial only if
\begin{equation}
%\label{}
-\lambda+k^{\tilde\sigma}\le \tilde u(x) \le \lambda+k^{1-\tilde\sigma}+\frac1{4\slb}\log k
\end{equation}
where $\tilde\sigma:=\frac1{12}$. 
Under these conditions Proposition~\ref{prop-3.3} permits us to replace $g_{k,\tilde u(x)}$ on the right of \eqref{E:4.84} by~$g$ from \eqref{E:3.22a}; Lemma~\ref{lemma-4.7} then also allows us to drop the barrier event~$\wt\BB_{n,k}(x)$ at the cost of an error term that vanishes as~$n\to\infty$ and~$k\to\infty$. 

The previous manipulations reduce the computation to the limit of the process of local maxima that, with the help of the tightness of the extreme level sets \twoeqref{E:4.75}{E:4.76} is reduced to the asymptotic \twoeqref{E:4.77}{E:4.78a}. Leaving the details to the reader, for any~$k_n\to\infty$ with~$n-k_n\to\infty$ this gives
\begin{equation}
%\label{E:3.4}
\begin{aligned}
E^\varrho\bigl(&\texte^{-\langle\wt\eta_{n,k_n},f\rangle}\bigr)
\\
&\,\,\,\underset{n\to\infty}\longrightarrow\,\, \E\biggl(\exp\Bigl\{-\int W(\textd x)\otimes\texte^{-2\slb\, h}\textd h\otimes\nu(\textd\phi)(1-\texte^{-f(x,h,\phi)})\Bigr\}\biggr),
\end{aligned}
\end{equation}
for a random measure~$W$ on~$[0,1]$ and~$\nu$ as in Proposition~\ref{prop-3.3}. Proceeding along the same argument as in the reduction of Theorem~\ref{thm-3.1} to Theorem~\ref{thm-2} then shows that the cluster process of GFF is that defined in \eqref{E:3.13ui}, thus identifying it with the cluster process of the local time of the random walk on~$\T_n$.
\end{proofsect}

\subsection{Random walk started from the leaves}
We now move to the proof of Theorem~\ref{thm-1} which will be deduced from Theorem~\ref{thm-2}. Here we note that the corresponding reduction in~\cite{BL4} for the maximum of the local time  relied on a convenient trick: The law of $\max_{x\in\BbbL_n}\ell_{\tau_\varrho}(x)$ is that of~$\max_{x\in\BbbL_n}L_t(x)$ conditioned on being positive, in the limit as~$t\downarrow0$. This capitalized on the observation that the law of $\max_{x\in\BbbL_n}L_t(x)$ does not depend on where the walk started from the root first hits the leaves. Unfortunately, this symmetry no longer helpful once we aim to include information on where the maximum is achieved, let alone how the other nearly-maximal values are distributed, and so we proceed along different lines. 

We will rely on the Markovian nature of the local time that yields the following decomposition of~$\ell_{\tau_\varrho}$: 

\begin{lemma}
\label{lemma-5.1}
Writing $x_0:=\varrho,x_1,\dots,x_n:=\zero$ for the vertices on the unique path in~$\T_n$ from the root to vertex~$\zero$, for any Borel sets $E_k\subseteq\R^{\T_k'}$ indexed by $k=1,\dots,n$ we have
\begin{equation}
\label{E:5.19i}
\begin{aligned}
P^\zero\biggl(\,\bigcap_{k=1}^{n}\Bigl\{\bigl\{\ell_{\tau_\varrho}(x)&\colon 
x\in \T_{n-k}^{\prime} (x_k)  \bigr\}\in E_{n-k}\Bigr\}\,\bigg|\, \sigma\bigl(\ell_{\tau_{\varrho}}(x_k)\colon k=1,\dots,n\bigr)\biggr)
\\
&=\prod_{k=1}^n P^\varrho\Bigl(\bigl\{L_{t_k}(x)\colon x\in
\T_{n-k}^{\prime}\bigr\}
\in E_{n-k}\Bigr)\Big|_{t_k:=\ell_{\tau_\varrho}(x_k)}\quad\text{\rm a.s.}
\end{aligned}
\end{equation}
where~$\T_{n-k}'(x_k)$ is the connected component of~$\T_n$ containing~$x_k$ when the edges on the path $(x_0,\dots,x_n)$ are removed from~$\T_n$. (We regard $\T_{n-k}'(x_k)$ as isomorphic to~$\T_{n-k}'$ defined earlier.)
\end{lemma}

\begin{proofsect}{Proof}
This follows using the exponential memoryless property of the exponential distribution similarly as Lemma~\ref{lemma-M}. 
\end{proofsect}

In order to understand the extremal process associated with $\ell_{\tau_\varrho}$ we thus need to understand the extremal behavior of independent processes $L_{t_1},L_{t_2},\dots, L_{t_n}$, for which we now have Theorem~\ref{thm-2} at our disposal, with the sequence $\{t_k\}_{k=1}^n$ set via $t_k:=\ell_{\tau_{\varrho}}(x_k)$. Here are the facts we need to know about this sequence:

\begin{lemma}
\label{lemma-5.2}
Let~$\{T_k\}_{k\ge0}$ be as defined in Corollary~\ref{cor-1.3}. Then for all $n\ge1$,
\begin{equation}
\label{E:5.2}
\bigl(\ell_{\tau_\varrho}(x_0),\dots,\ell_{\tau_\varrho}(x_n)\bigr)\text{\rm\ under }P^\zero
\,\,\,\laweq\,\,(T_0,\dots,T_n).
\end{equation}
In addition,
\begin{equation}
\label{E:5.3}
\lim_{k \to \infty} \inf_{n \ge k}\,P\biggl(\,\bigcap_{j=k}^n\bigl\{T_j\le\alpha j\log(j)\bigr\}\biggr)=1
\end{equation}
holds true for each~$\alpha>1$. 
\end{lemma}

\begin{proofsect}{Proof}
The law of $(\ell_{\tau_\varrho}(x_0),\dots,\ell_{\tau_\varrho}(x_n))$ was identified in Zhai~\cite[Corollary 5.3]{Zhai} (it appears that a factor $1/2$ is missing there). The formula \eqref{E:5.3} then follows from a union bound and a standard tail estimate for Gaussian random variables.
\end{proofsect}

The tail bound in the previous lemma shows that the relevant contribution to the extremal process of~$\ell_{\tau_\varrho}$ arrive only in the trees $\T_{n-k}(x_k)$ with~$k$ small:

\begin{corollary}
\label{cor-5.3}
For each~$\lambda>0$,
\begin{equation}
\label{E:5.5}
\sup_{n\ge k}P^\zero\biggl(\exists x\in\BbbL_n\cap\bigcup_{j=k}^n 
\T_{n-j}^{\prime}(x_j)
\colon\,
\sqrt{\ell_{\tau_\varrho}(x)}\ge m_n-\lambda\biggr)\,\,\underset{k\to\infty}\longrightarrow\,\,0.
\end{equation}
\end{corollary}

\begin{proofsect}{Proof}
Suppose~$k\ge1$ is such that the event in  \eqref{E:5.3} with~$\alpha:=2$ occurs. Then \eqref{E:3.50} gives that, for each~$j=k,\dots, n$,
\begin{equation}
%\label{}
m_n-a_{n-j}(T_j) - \sqrt{T_j}\ge\slb\,j-\sqrt{2j\log j}-\frac1{4\slb}\log(2j\log j)+O(1),
\end{equation}
where we also used that $t\mapsto a_{n-j}(t)+\sqrt t$ is non-decreasing. As this grows linearly with~$j$, Lemma~\ref{lemma-3.2a} along with the observation that the maximum of the local time on~$\T_{n-j}$ dominates that on~$\T_{n-j}'$ show that, on the event in \eqref{E:5.3}, the probability that the maximum of~$L_{T_j}$ exceeds~$m_n-\lambda$ is exponentially small in~$j$. 
Along with \eqref{E:5.3}, a routine union bound then proves the claim.
\end{proofsect}

We are now ready for:

\begin{proofsect}{Proof of Theorem~\ref{thm-1} except for \eqref{E:1.7}}
Given~$n\ge k>0$ and a sample of~$\ell_{\tau_\varrho}$, write $\eta_n$ for the process on the left-hand side of \eqref{E:1.5} and denote its truncated version by
\begin{equation}
%\label{E:1.5}
\eta_{n,k}:=\sum_{x\in\BbbL_n}1_{\{\theta_n(x)\ge b^{-k}\}}\delta_{\theta_n(x)}\otimes\delta_{\sqrt{\ell_{\tau_\varrho}(x)}-m_n}.
\end{equation}
Given $t\ge0$ and a sample of~$L_t$ on~$\T_n$, let also
\begin{equation}
%\label{}
\zeta_n^t:=\sum_{x\in\BbbL_n}\delta_{\theta_n(x)}\otimes\delta_{\sqrt{L_t(x)}-m_n}.
\end{equation}
Pick~$f\colon[0,1]\times\R\to[0,\infty)$ continuous with support in $[0,1]\times[-\lambda,\lambda]$, for some~$\lambda>0$, and observe that
\begin{equation}
\label{E:5.9}
\lim_{k\to\infty}\,\limsup_{n\to\infty}\Bigl|\,E^\zero\bigl(\texte^{-\langle\eta_{n,k},f\rangle}\bigr)-E^\zero\bigl(\texte^{-\langle\eta_n,f\rangle}\bigr)\Bigr|=0
\end{equation}
by Corollary~\ref{cor-5.3}.

Next let us call upon Lemma~\ref{lemma-5.1} to produce a coupling of~$(\ell_{\tau_\varrho}(x_0),\dots,\ell_{\tau_\varrho}(x_n))$ with conditionally independent processes $(L_{t_0},\dots,L_{t_n})$ on~$\T_n$ for the choices~$t_j:=\ell_{\tau_\varrho}(x_j)$ such that~$\ell_{\tau_\varrho}$ coincides with appropriately parametrized~$L_{t_j}$ on~$\T_{n-j}'(x_j)$, for each~$j=0,\dots,n$. Under this coupling we have
\begin{equation}
\begin{aligned}
\langle\eta_{n,k},f\rangle
&=\sum_{j=0}^{k-1}\sum_{x\in\BbbL_n} 1_{[b^{-j-1},b^{-j})}\bigl(\theta_n(x)\bigr) f\Bigl(\theta_n(x),\sqrt{\ell_{\tau_\varrho}(x)}-m_n\Bigr)
\\
&=\sum_{j=0}^{k-1}\sum_{x\in\BbbL_{n-j}} 1_{[b^{-1},1)}\bigl(\theta_{n-j}(x)\bigr) f\Bigl(b^{-j}\theta_{n-j}(x),\sqrt{L_{t_j}(x)}-m_n\Bigr)\Big|_{t_j=\ell_{\tau_\varrho}(x_j)}.
\end{aligned}
\end{equation}
Abbreviating
\begin{equation}
%\label{}
f_{n,j}(x,h):=1_{[b^{-1},1)}(x) f\bigl(b^{-j}x,h+m_{n-j}-m_n\bigr),
\end{equation}
Lemmas~\ref{lemma-5.1}--\ref{lemma-5.2} (along with the fact that the local time process~$L_t$ on~$\T_{n-i}'$ is a restriction of the process on~$\T_{n-i}$) give
\begin{equation}
\label{E:5.28i}
E^\zero\bigl(\texte^{-\langle\eta_{n,k},f\rangle}\bigr) = E\biggl(\,
\prod_{i=1}^{k-1} E^\varrho\bigl(\texte^{-\langle\zeta_{n-i}^{t_i},f_{n,i}\rangle}\bigr)\big|_{t_i:=T_i}\biggr),
\end{equation}
where the outer expectation is with respect to the law of $(T_0,T_1,\dots)$. Here the term corresponding to~$j=0$ has been dropped because~$T_0=0$ a.s.\ and so the corresponding process $\zeta^{t_0}_n$ is thus zero when~$t_0=T_0$. 

Note that $m_{n-j}-m_n=-j\slb+o(1)$ as~$n\to\infty$. Using the uniform continuity of~$f$ to absorb the error and the lack of atoms of~$Z_u$ to deal with the discontinuity of~$f_{n,j}$ at~$x:=b^{-1}$, 
Theorem~\ref{thm-2} shows that $E^\varrho(\texte^{-\langle\zeta_{n-j}^{u},f_{n,j}\rangle})$ converges, for each~$u\ge0$, to
\begin{equation}
%\label{}
\E\biggl(\exp\Bigl\{-\int 1_{[b^{-1},1)}(x) Z_{u}(\textd x)\otimes\texte^{-2h\slb}\textd h\otimes\DD(\textd\xi)\bigl(1-\texte^{-\langle\xi, f_j(x,h+\cdot)\rangle}\bigr)\Bigr\}\biggr),
\end{equation}
where the expectation is with respect to the law of~$Z_u$ and~$f_j(x,h):=f(b^{-j}x,h-j\slb)$. A change of variables casts this in the form
\begin{equation}
%\label{}
\E\biggl(\exp\Bigl\{-\int Z^{(j)}_{u}(\textd x)\otimes\texte^{-2h\slb}\textd h\otimes\DD(\textd\xi)\bigl(1-\texte^{-\langle\xi, f(x,h+\cdot)\rangle}\bigr)\Bigr\}\biggr),
\end{equation}
where
\begin{equation}
%\label{}
Z^{(j)}_u(\textd x) := b^{-2j} 1_{[b^{-j-1},b^{-j})}(x)
Z_u\bigl(b^j\textd x\bigr).
\end{equation}
Plugging the above convergence statements into \eqref{E:5.28i} it follows that if $\{Z_{T_j}\}_{j\ge0}$ is the family of random measures 
associated with~$\{T_j\}_{j\ge0}$ as specified in \eqref{E:1.21a} and 
\begin{equation}
\label{E:5.15}
\ZZ_k(\textd x):=\sum_{j=1}^{k-1}b^{-2j} 1_{[b^{-j-1},b^{-j})}(x)
Z_{T_j}\bigl(b^j\textd x\bigr)
\end{equation}
then the Bounded Convergence Theorem gives
\begin{equation}
\label{E:5.33i}
\begin{aligned}
&E^\zero\bigl(\texte^{-\langle\eta_{n,k},f\rangle}\bigr) 
\\
&\underset{n\to\infty}\longrightarrow\,\,
\E\biggl(\exp\Bigl\{-\int \ZZ_k(\textd x)\otimes\texte^{-2h\slb}\textd h\otimes\DD(\textd\xi)\bigl(1-\texte^{-\langle\xi, f(x,h+\cdot)\rangle}\bigr)\Bigr\}\biggr)
\end{aligned}
\end{equation}
for each~$k\ge1$. 

To finish the proof of convergence \eqref{E:1.5} observe that $\ZZ_k$ increases to the measure~$\ZZ$ defined by the right-hand side of \eqref{E:1.23a} and so, by the Bounded Convergence Theorem, the right-hand side of \eqref{E:5.33i} tends to that with~$\ZZ_k$ replaced by~$\ZZ$ as~$k\to\infty$. Thanks to \eqref{E:5.9}, the limiting expression is then also the limit of $E^\zero(\texte^{-\langle\eta_n,f\rangle})$ as~$n\to\infty$ (and the limit thus exists) for each~$f$ as above. 

The properties of~$\DD$ were proved along with Theorem~\ref{thm-2}.
To address the properties of~$\ZZ$, note that ~$\ZZ$ is finite a.s.\ by the tightness of~$\max_{x\in\BbbL_n}\ell_{\tau_\varrho}(x)$ proved in~\cite[Theorem~1.1]{BL4}. While the measures~$Z_t$ may vanish with positive probability, by \eqref{E:1.13} this probability decreases with~$t$. Since~$T_k\to\infty$ a.s.\ it follows that a positive fraction of $Z_{T_i}$-measures are non-vanishing almost surely, thus showing $\ZZ([0,\epsilon))>0$ for each~$\epsilon>0$~a.s. That $\ZZ(\{0\})=0$~a.s.\ follows from the construction of~$\ZZ$ as an increasing limit of~$\ZZ_k$ and the fact that, since~$\ZZ_k$ is supported in~$[b^{-k},1)$ and $\ZZ([0,1])<\infty$, we have $\ZZ([0,b^{-k}))\le\ZZ([0,1])-\ZZ_k([0,1])\to0$ as~$k\to\infty$. 
\end{proofsect}

From the above we immediately get:

\begin{proofsect}{Proof of Corollary~\ref{cor-1.3}}
The above proof shows that~$\ZZ_k$ defined in \eqref{E:5.15} increases in law to measure~$\ZZ$ governing the limit extremal process. Hence we get \eqref{E:1.23a}.
\end{proofsect}

\begin{proofsect}{Proof of Corollary~\ref{cor-1.4}}
For $n\ge k\ge1$ abbreviate $\BbbL_{n,k}:=\BbbL_n\cap \T_{n-k}'(x_k)$. Corollary~\ref{cor-5.3} tells us that the maximum of~$\ell_{\tau_\varrho}$ occurs in $\bigcup_{j=1}^{k} \BbbL_{n,j}$ with probability tending to one as~$n\to\infty$ and~$k\to\infty$. 
Writing $M_{n,j}:= \max_{x\in\BbbL_{n,j}} \sqrt{\ell_{\tau_\varrho}(x)}- m_n$, the joint law of $(M_{n,1},\dots,M_{n,k})$ has a limit described by
\begin{equation}
%\label{}
P^\zero\biggl(\,\bigcap_{j=1}^k\{M_{n,j}\le u_j\}\biggr)\,\,\underset{n\to\infty}\longrightarrow\,\,
\E\biggl(\,\,\prod_{j=1}^k\texte^{-\frac1{2\slb}b^{-2j}Z_{T_j}([b^{-1},1))\texte^{-2\slb\,u_j}}\biggr).
\end{equation}
Another way to write this is as 
\begin{equation}
\label{E:5.18}
\begin{aligned}
\bigl\{M_{n,j}\colon j=1,&\dots,k\bigr\}
\\
&\underset{n\to\infty}\lawarrow\,\,
\Bigl\{\frac1{2\slb}\bigl(\log Z_{T_j}([b^{-1},1))-2j\log b+G_j'\bigr)\colon j=1,\dots,k\Bigr\},
\end{aligned}
\end{equation}
where~$G_1',\dots,G_k'$ are i.i.d.\ standard Gumbel random variables shifted by $\log(2\slb)$. 

As the $G_i$'s are continuously distributed, the largest term in the sequence on the right of \eqref{E:5.18} is unique a.s. Thanks to Corollary~\ref{cor-5.3} again, so must be the maximizer of the infinite sequence. This yields the claim by the fact that, within~$\BbbL_{n,k}$, the maximizer of~$x\mapsto\ell_{\tau_\varrho}(x)$ is uniformly distributed by the symmetries of the tree.
\end{proofsect}

\subsection{Limit characterization of~$\ZZ$-measure}
To complete the proof of Theorem~\ref{thm-1} it remains to prove the limit characterization \eqref{E:1.7}. Our argument relies on the convergence of the total mass of the measure on the left of~\eqref{E:1.7} to the total mass of~$\ZZ$. This was claimed in \cite[Theorem~1.5]{BL4}; unfortunately, the proof of this part appears to be missing. As this fact enters rather delicately our proof of~\eqref{E:1.7}, we state and prove the result here.

\begin{proposition}
\label{prop-5.5}
Let $\wt C_\star:=2c_\star\slb$ for~$c_\star$ as in Proposition~\ref{prop-3}. Then
\begin{equation}
%\label{}
\wt C_\star\, b^{-2n} \sum_{x \in \BbbL_n} \Bigl(n\slb - 
\sqrt{\ell_{\tau_\varrho}(x)}\,\Bigr)^+ \ell_{\tau_\varrho}(x)^{1/4}\,
	\texte^{2\slb\,\sqrt{\ell_{\tau_\varrho}(x)}}\,\,\,\underset{n\to\infty}\lawarrow\,\,\,\ZZ\bigl([0,1]\bigr),
\end{equation}
where~$\ZZ$ is the measure from Theorem~\ref{thm-1}. 
\end{proposition}

The proof of this proposition will rely on one useful fact:

\begin{lemma}
\label{lemma-5.6}
Given~$n\ge1$, let~$\ell^{(n)}_{\tau_\varrho}$ be the local time on~$\T_n$ for the walk started from~$\zero\in\BbbL_n$ and stopped upon first hit of~$\varrho$. Then for all~$1\le k\le n$, 
\begin{equation}
%\label{}
\bigl\{\ell^{(n)}_{\tau_\varrho}(x)\colon x\in\T_k\bigr\}\laweq \bigl\{\ell^{(k)}_{\tau_\varrho}(x)\colon x\in\T_k\bigr\}.
\end{equation}
In particular, the laws of $\{\ell^{(n)}_{\tau_\varrho}\colon n\ge1\}$ are Kolmogorov-consistent.
\end{lemma}

\begin{proofsect}{Proof}
This follows by combining \eqref{E:5.19i} and \eqref{E:5.2}, along with the fact that the law of~$\{T_k\}_{k\ge0}$ does not depend on~$n$. (A direct argument can be based on the memoryless property of the exponential distribution and the fact that the walk started from~$\zero$ on~$\T_n$ enters~$\T_k$ through~$\zero$ on~$\T_k$.)
\end{proofsect}

Using this we now give:

\begin{proofsect}{Proof of Proposition~\ref{prop-5.5}}
Given $n\ge1$,  write $\ZZ^{(n)}$ for the measure from \eqref{E:1.7}. Our goal is to show that $\ZZ^{(n)}([0,1])$ tends in law to~$\ZZ([0,1])$ as~$n\to\infty$. Let~$s\in\R$ and let us continue writing $\ell^{(n)}_{\tau_\varrho}$ for the local time~$\ell_{\tau_\varrho}$ on~$\T_n$ under~$P^\zero$. Conditioning on the values in~$\T_k$ for some~$1\le k\le n$, Lemmas~\ref{lemma-5.1} and~\ref{lemma-5.6} imply
\begin{equation}
\label{E:5.37}
P^\zero\Bigl(\,\max_{x\in\BbbL_n} \sqrt{\ell^{(n)}_{\tau_\varrho}(x)}  \le m_n+s\Bigr)
=E^\zero\biggl(\,\prod_{z\in\BbbL_k}P^\varrho\Bigl(\max_{x\in\BbbL_{n-k}}  \sqrt{L_t(x)}  \le m_n+s\Bigr)\Bigr|_{t:=\ell^{(k)}_{\tau_\varrho}(z)}\biggr).
\end{equation}
We will represent the probabilities under the product using the asymptotic form in Proposition~\ref{prop-3}. For this we assume that
\begin{equation}
\label{E:5.38i}
\AA_k':=\Bigl\{\max_{z\in\BbbL_k}\sqrt{\ell^{(k)}_{\tau_\varrho}(z)}\le m_k+\log\log k\Bigr\}
\end{equation}
occurs and note that, for $t:= \ell^{(k)}_{\tau_\varrho}(z)$, we have $m_n+s = a_{n-k}(t)+\sqrt t+u$ with
\begin{equation}
\label{E:5.37a}
\begin{aligned}
u&:=m_n+s-a_{n-k}(t)-\sqrt{t} 
\\
&\,= k\slb -\frac1{8\slb}\log t+s+\frac1{4\slb}\log\Bigl(\frac{n-k+\sqrt t}{n}\Bigr) - \sqrt{t}
 + \frac{3}{4\sqrt{\log b}} \log \left(1-\frac{k}{n} \right)
\\
&\,\ge k\slb -\frac1{8\slb}\log t+O(1) - m_k-\log\log k
\\
&\,\ge \frac1{\slb}\log k - 2\log\log k+O(1),
\end{aligned}
\end{equation}
where the third inequality follows by plugging the constraint from \eqref{E:5.38i}. 
Also note that for such~$u$ and~$t$ we get 
\begin{equation}
\label{E:5.40iu}
\begin{aligned}
u\texte^{-2u\slb}=
\biggl(1+O\Bigl(\frac{\log\log k}{\log k}\Bigr)&\biggr)\Bigl(k\slb -\sqrt{\ell_{\tau_\varrho}(x)}\Bigr)^{ + }
\\
&\times b^{-2k}\ell_{\tau_\varrho}(x)^{1/4}\texte^{2\slb\sqrt{\ell_{\tau_\varrho}(x)}}\texte^{-2s\slb},
\end{aligned}
\end{equation}
where the error term is random but bounded by a deterministic constant times the stated ratio.
For~$k$ large, $u$ is large and the quantity \eqref{E:5.40iu} is small uniformly for all terms in the product in \eqref{E:5.37}. This allows us to apply Proposition~\ref{prop-3} to the terms with $t:= \ell^{(k)}_{\tau_\varrho}(z)\ge\log\log k$ while handling the remaining terms using  the bound \eqref{E:3.6x} in Lemma~\ref{lemma-3.2a}. As is seen from \eqref{E:5.40iu}, each of the latter terms contributes at most $\exp\{O(1)b^{-2k}k^2\}$ to the product and contributes at most $O(b^{-2k}k^2)$ to~$\ZZ^{(k)}$. These terms are thus negligible and we may rewrite the expectation in \eqref{E:5.37} as
\begin{equation}
%\label{}
O\bigl(P^\zero(\AA_k'^\cc)\bigr)+E^{\zero}\biggl(\exp\Bigl\{\tilde o(1) -\wt C_\star^{-1}(c_\star+o(1)) \texte^{-2s\slb}\ZZ^{(k)}\bigl([0,1]\bigr)\Bigr\}\biggr),
\end{equation}
where the~$o(1)$-terms vanish in the limit as $n\to\infty$ followed by~$k\to\infty$. 

Appealing to the structured-process convergence in Theorem~\ref{thm-3.1}, a simple approximation argument shows that $P^\zero(\AA_k'^\cc)\to0$ as~$k\to\infty$ while the left-hand side of \eqref{E:5.37} converges to
\begin{equation}
%\label{}
\E\biggl(\exp\Bigl\{-(2\slb)^{-1}\texte^{-2s\slb}\ZZ([0,1])\Bigr\}\biggr)
\end{equation}
as~$n\to\infty$. Since~$s$ is arbitrary, taking $k\to\infty$ we get convergence of the Laplace transforms of the laws of~$\ZZ^{(k)}([0,1])$ to that of~$\ZZ([0,1])$. This is enough to imply the claim.
\end{proofsect}

With the above in hand, we are ready to give:

\begin{proofsect}{Proof of \eqref{E:1.7} in Theorem~\ref{thm-1}}
Set~$\wt C_\star:=2c_\star\slb$ for~$c_\star$ as in Proposition~\ref{prop-3}. Given $n\ge k\ge1$ and a sample $\ell_{\tau_\varrho}$ of the local time on~$\T_n$ from~$P^{\zero}$, denote
\begin{equation}
%\label{}
\begin{aligned}
\ZZ^{(n)}_k:=\wt C_\star\, b^{-2n} &\sum_{x \in \BbbL_n}\sum_{j=0}^{k-1}1_{[b^{-j-1},b^{-j})}(\theta_n(x)) 
\\&\times\biggl((n-j)\slb - 
\sqrt{\ell_{\tau_\varrho}(x)}\,\biggr)^+ \ell_{\tau_\varrho}(x)^{1/4}\,
	\texte^{2\slb\,\sqrt{\ell_{\tau_\varrho}(x)}}\delta_{\theta_n(x)}.
\end{aligned}
\end{equation}
Observe that, if the factor $n-j$ were replaced by~$n$, then this would simply be the restriction of the measure in \eqref{E:1.7} to~$[b^{-k},1)$.

Next let~$x_0:=\varrho,x_1,\dots,x_n:=\zero$ label the vertices on the unique path from the root to~$\zero$. Lemma~\ref{lemma-M} implies that, conditional on $(\ell_{\tau_\varrho}(x_0),\dots, \ell_{\tau_\varrho}(x_n))$, the local time process~$\ell_{\tau_\varrho}$ restricted to $\{x\in\BbbL_n\colon \theta_n(x)\in[b^{-j-1},b^{-j})\}$ has the law of~$L_{t_j}$ on~$\{x\in\BbbL_{n-j}\colon \theta_{n-j}(x)  \in[b^{-1},1)\}$ for $t_j:=\ell_{\tau_\varrho}(x_j)$. Using this we get
\begin{equation}
\label{E:5.36}
\ZZ^{(n)}_k(\textd x)\laweq\sum_{j=1}^{k-1} b^{-2j}\,1_{[b^{-j-1},b^{-j})}(x) \wt Z_{t_j}^{(n-j)}(b^j\textd x)\bigl|_{t_j:=T_j},
\end{equation}
where $(T_0,\dots,T_n)$ are independent with law given in Corollary~\ref{cor-1.3} and $(\wt Z^{(n-1)}_{t_1},\dots,\wt Z^{(1)}_{t_{n-1}})$ are independent samples of measures from \eqref{E:1.18}; i.e., measures of the form
\begin{equation}
%\label{}
\wt Z^{(k)}_t:=\wt C_\star\,b^{-2k} \sum_{x \in \BbbL_k} \Bigl(k\slb - 
\sqrt{L_t^{(k)}(x)}\,\Bigr)^+ L_t^{(k)}(x)^{1/4}\,
	\texte^{2\slb\,\sqrt{L_t^{(k)}(x)}}\delta_{\theta_k(x)}
\end{equation}
with~$\{t\mapsto L^{(k)}_t\colon k\ge1\}$ independent for different~$k$. We again dropped the $j:=0$ term in \eqref{E:5.36} due to the fact that~$T_0=0$ a.s. Invoking the convergence \eqref{E:1.18} we then get
\begin{equation}
\label{E:5.38}
\ZZ^{(n)}_k\,\,\underset{n\to\infty}\lawarrow\,\,\ZZ_k
\end{equation}
where~$\ZZ_k$ is as in \eqref{E:5.15}. Recall also that $\ZZ_k\,\lawarrow\,\ZZ$ as~$k\to\infty$.

Let~$\ZZ^{(n)}$ denote the measure on the left of \eqref{E:1.7}. Clearly, $\ZZ^{(n)}_k\le \ZZ^{(n)}$ pointwise as measures for each~$k=1,\dots,n$.  Moreover, by Proposition~\ref{prop-5.5}, 
\begin{equation}
\label{E:5.45}
\ZZ^{(n)}\bigl([0,1]\bigr)\,\,\underset{n\to\infty}\lawarrow\,\,\ZZ\bigl([0,1]\bigr).
\end{equation}
In particular, $\{\ZZ^{(n)}\colon n\ge1\}$ is a tight family of random Borel measures on~$[0,1]$ which permits consideration of weak subsequential limits. But the inequality $\ZZ^{(n)}_k\le \ZZ^{(n)}$ along with \eqref{E:5.38} and the weak convergence of~$\ZZ_k$ to~$\ZZ$ show that any subsequential weak limit~$\ZZ'$ of $\{\ZZ^{(n)}\colon n\ge1\}$  dominates~$\ZZ$ in law and yet, by \eqref{E:5.45}, has the same total mass. This forces~$\ZZ'$ to coincide with~$\ZZ$, proving the claim.
\end{proofsect}

%\newpage
\section{Barrier estimates for BRW}
\label{sec-6}\noindent
To settle all aspects of the proof of Theorem~\ref{thm-3.1}, it remains to prove the convergence in Proposition~\ref{prop-3.3} and the barrier estimate in Lemma~\ref{lemma-4.7} for the Branching Random Walk with normal step distribution.

\subsection{Gaussian random walk above a barrier}
We begin by a statement that belongs to the theory of inhomogeneous Ballot Theorems. This subject has been treated systematically in many sources; e.g., Bramson~\cite{Bramson1},  Biskup and Louidor~\cite{BL3}, Mallein~\cite{M-review} or Cortines, Hartung and Louidor~\cite{CHL17Sup}; unfortunately, none of these treatments seem to give precisely what we need due to either a different setting or non-uniformity in the relevant parameters (specifically, allowing that~$u$ in \eqref{E:3.22} scales as a power of~$k$). We will therefore work out the needed details here focusing solely on the setting relevant for the above claims. 

\begin{proposition}
\label{prop-6.1}
Fix any $\sigma \in(0,1/10)$ and let~$X_0,X_1,\dots$ be a random walk with step distribution~$\NN(0,1/2)$. For each~$\epsilon\in(0,1)$ and each $\delta\in(0,1)$ there exists~$a_0=a_0(\epsilon,\delta)\ge1$ 
such that for all $a\ge a_0$, all naturals $k\ge2$ and all reals~$r,u$ satisfying
\begin{equation}
\label{E:6.1a}
r,u\ge a^{1/\sigma}\quad\text{\rm and}\quad ru\le k^{1-\delta}
\end{equation}
and all~$\gamma\colon\{0, \dots,k\}\to\R$ obeying
\begin{equation}
\label{E:4.31}
\bigl|\gamma(i)\bigr|\le a+\bigr(\min\{i,k-i\}\bigr)^{\sigma},\quad i=0,\dots,k, 
\end{equation}
we have
\begin{equation}
\label{E:4.32}
(1-\epsilon)\frac{4ru}k\le
E\biggl(\,\prod_{j=1}^{k-1}1_{\{X_j\ge\gamma(j)\}}\,\bigg|\, X_0=r,\,X_k=u\biggr)\le(1+\epsilon)\frac{4ru}k.
\end{equation}
\end{proposition}

A key step of the proof is to show that, conditional on $X\ge-|\gamma|$, the walk actually obeys~$X\ge|\gamma|$. We prove this first as:

\begin{lemma}
\label{lemma-4.6b}
Let~$X$ be as in Proposition~\ref{prop-6.1}. For each~$\epsilon>0$ and $\sigma\in(0,1/\10)$ there exists~$a_1\ge1$ such that for all~$a\ge a_1$ and all non-negative~$\gamma$ satisfying \eqref{E:4.31},
\begin{equation}
\label{E:4.46b}
E\biggl(\,1_{\{\min_{1\le j\le k-1}
[X_j-\gamma(j)]<0\}}\prod_{j=1}^{k-1}1_{\{X_j\ge -\gamma(j)\}}\,\bigg|\, X_0=r,\,X_k=u\biggr)\le \epsilon\frac{ru}k
\end{equation}
holds whenever~$r,u\ge a^{1/\sigma}$, uniformly in~$k\ge2$. 
\end{lemma}

\begin{proofsect}{Proof}
Abbreviate the law of~$X$ conditioned on~$X_i=r$ and~$X_j=u$ as~$P_{i,j}^{r,u}$. We first make one useful observation. Given a path~$X$ of the random walk, note that if $i\mapsto X_i-\gamma(i)$ is minimized at some~$\ell$, then $X_i\ge X_\ell-\gamma(\ell)+\gamma(i)\ge X_\ell-\gamma(\ell)$ for all~$i=1,\dots,k-1$. Hence we get
\begin{equation}
\label{E:4.46w}
\begin{aligned}
P_{0,k}^{r,u}\Bigl(\,\min_{1\le i<k}[X_i-&\gamma(i)]=X_\ell-\gamma(\ell)\,\Big|\,X_\ell = s\Bigr)
\\
&\le P_{0,\ell}^{r,s}\biggl(\,\bigcap_{i=1}^{\ell-1}\{X_i\ge s-\gamma(\ell)\}\biggr)
P_{\ell,k}^{s,u}\biggl(\,\bigcap_{i=\ell+1}^{k-1}\{X_i\ge s-\gamma(\ell)\}\biggr)
\end{aligned}
\end{equation}
for all~$\ell=1,\dots,k-1$.
For~$s\in[-\gamma(\ell),\gamma(\ell)]$ the homogeneous Ballot Theorem (or the argument in the last part of the proof of Proposition~\ref{prop-6.1}) bounds the first probability by a constant that depends only on the distribution of~$X$ times $\ell^{-1}[r+\gamma(\ell)]\gamma(\ell)$. A similar argument applies to the second probability as well. 
 
Next observe that, by a calculation with Gaussian densities, the probability density~$f$ of~$X_\ell$ under $P_{0,k}^{r,u}$ equals
\begin{equation}
\label{E:6.6w}
f(s)=\frac1{\sqrt{\pi}}\sqrt{\frac{k}{\ell(k-\ell)}}\exp\biggl\{-\Bigl(\frac1\ell+\frac1{k-\ell}\Bigr)^{-1}\Bigl[\frac1\ell(r-s)+\frac1{k-\ell}(u-s)\Bigr]^2\biggr\}
\end{equation}
which is readily bounded as
\begin{equation}
\label{E:4.48b}
f(s)\le \frac1{\sqrt{\pi}}\sqrt{\frac{k}{\ell(k-\ell)}}\,\texte^{-[(\min\{r,u\}-s)_+]^2/\min\{\ell,k-\ell\}}.
\end{equation}
Let $E_{0,k}^{r,u}$ denote the expectation with respect to~$P_{0,k}^{r,u}$. Partitioning according to the first maximizer of $i\mapsto X_i-\gamma(i)$ we then dominate the quantity of interest as 
\begin{equation}
\label{E:4.48a}
\begin{aligned}
E_{0,k}^{r,u}\biggl(\,&1_{\{\min_{1\le j\le k-1}
[X_j-\gamma(j)]<0\}}\prod_{j=1}^{k-1}1_{\{X_j\ge -\gamma(j)\}}\biggr)
\\
&\le \tilde c\sum_{\ell=1}^{k-1} \gamma(\ell)^3 [r+\gamma(\ell)][u+\gamma(\ell)]\frac{\sqrt k}{[\ell(k-\ell)]^{3/2}}\,\texte^{-\frac{[(\min\{r,u\}-\gamma(\ell))_+]^2}{\min\{\ell,k-\ell\}}},
\end{aligned}
\end{equation}
where the terms $\ell^{-1} \gamma (\ell) [r+\gamma(\ell)]$ and $( k  -\ell)^{-1}\gamma(\ell)[u+\gamma(\ell)]$ arise from the estimates of the conditional probability as discussed after \eqref{E:4.46w}, the exponential on the right dominates that in \eqref{E:4.48b} uniformly in 
$s \in [-\gamma (\ell), \gamma (\ell)] $ and another factor~$2\gamma(\ell)$ arises from the integral over~$s$ subject to the aforementioned restrictions. 

Assuming $a\ge1$ and $r,u\ge a^{1/\sigma}\ge1$ we have
\begin{equation}
\label{E:6.9i}
\bigl[r+\gamma(\ell)\bigr]\bigl[u+\gamma(\ell)\bigr]\le ru[1+\gamma(\ell)]^2.
\end{equation}
We now finally call upon the assumed bound \eqref{E:4.31} on~$\gamma$ which permits us to dominate the sum in \eqref{E:4.48a} by $\frac{ru}k$ times
\begin{equation}
\label{E:4.50}
 \sum_{\ell=1}^{k-1}\,\bigl[1+a+(\min\{\ell,k-\ell\})^{\sigma}\bigr]^5\frac{k^{3/2}}{[\ell(k-\ell)]^{3/2}}\,\texte^{-\frac{[( a^{1/\sigma}-a-(\min\{\ell,k-\ell\})^{\sigma} )_+]^2}{\min\{\ell,k-\ell\}}}.
\end{equation}
This sum is bounded by twice that for~$\ell\le \lceil k/2\rceil$, with $k^{3/2}$ subsequently cancelling against the lower bound on $(k-\ell)^{3/2}$ arising in the denominator. We then split the resulting sum according to whether~$\ell^\sigma\le a/2$ or not. Under $\sigma<1/10$, the former part is checked to be exponentially small in~$a^{1/\sigma}$ while the latter part is bounded by a constant times $a^{(-1/2+5\sigma)/\sigma}$. As $5\sigma<1/2$, \eqref{E:4.50} is thus bounded uniformly in~$k\ge1$ and tends to zero as~$a\to\infty$. We conclude that \eqref{E:4.46b} holds for~$a$ sufficiently large. 
\end{proofsect}

With \eqref{E:4.46b} established, we now use it to prove the bounds in \eqref{E:4.32}.

\begin{proofsect}{Proof of Proposition~\ref{prop-6.1}}
Denote~$\gamma_\star(i):=a+\bigr(\min\{i,k-i\}\bigr)^{\sigma}$ and assume~$a\ge a_1$, for~$a_1$ related to~$\epsilon$ as in Lemma~\ref{lemma-4.6b}.
Abbreviate the expectation in \eqref{E:4.32} as~$F_k(\gamma)$. If~$|\gamma|\le\gamma_\star$ and~$r,u\ge a^{1/\sigma}$, then the above tells us
\begin{equation}
\label{E:4.51}
F_k(\gamma)\ge F_k(\gamma_\star)\ge F_k(-\gamma_\star)-\frac{ur}k\epsilon\ge F_k(0)-\frac{ur}k\epsilon.
\end{equation}
To estimate the right-hand side, let $B=\{B_s\colon 0\le s\le k/2\}$ be the standard Brownian bridge of time-length $k/2$ and endpoints fixed to~$B_0=r$ and~$B_{k/2}=u$ a.s. Thanks to the specific step distribution, under $P^{r,u}_{0,k}$ the family $\{X_j\colon j=0,\dots ,k\}$ is equidistributed to $\{B_{j/2}\colon j=0,\dots,k\}$. Using $ P^{r,u}$ to denote the law of~$B$, the Reflection Principle  gives
\begin{equation}
\label{E:4.53i}
F_k(0)\ge  P^{r,u}\Bigl(\,\inf_{0\le s\le k/2}B_s\ge0\Bigr)=1-\exp\Bigl\{-4\frac{ur}k\Bigr\}.
\end{equation}
The right-hand side exceeds $(1-\epsilon/2)\frac{4ur}k$ when $ur\le k^{1-\delta}$ and both~$u$ and~$r$ (and thus~$k$) are sufficiently large. Plugging this in \eqref{E:4.51} proves the lower bound in \eqref{E:4.32}. 

For the bound on the right of \eqref{E:4.32}, a similar argument as one used for \eqref{E:4.51} shows
\begin{equation}
\label{E:4.37}
F_k(\gamma)\le F_k(-\gamma_\star)\le F_k(\gamma_\star)+\frac{ur}k\epsilon.
\end{equation}
Proceeding as in the proof of \cite[Lemma 4.16]{BL3}, we now bound $F_k(\gamma_\star)$ as follows. First note that, for $B$ the above Brownian bridge  and each~$j=1,\dots,k$,
\begin{equation}
%\label{}
W^{(j)}_s:=B_{(1-s)\frac{j-1}2+s\frac j2}-(1-s)B_{\frac{j-1}2}-sB_{\frac j2},\quad 0\le s\le 1,
\end{equation}
defines a family $\{W^{(j)}\colon j=1,\dots,k\}$ of i.i.d.\ Brownian bridges of time-length~$1$ and endpoints fixed to zero. Abbreviating $W_\star^{(j)}:= \inf_{0\le s\le1}W_s^{(j)}$, we readily check
\begin{equation}
%\label{}
\bigcap_{j=1}^{k-1}\bigl\{B_{\frac j2}\ge\gamma_\star(j)\bigr\}\cap\bigcap_{j=1}^k\Bigl\{W_\star^{(j)}\ge-\min\{\gamma_\star(j-1),\gamma_\star(j)\}\Bigr\}
\\
\subseteq\bigcap_{0\le s\le k/2}\bigl\{B_s\ge0\bigr\}.
\end{equation}
The properties of Brownian motion imply that the two giant intersections on the left are independent. This yields
\begin{equation}
\label{E:4.57}
F_k(\gamma_\star)\prod_{j=1}^k P^{0,0}\Bigl(W_\star^{(j)}\ge-\min\{\gamma_\star(j-1),\gamma_\star(j)\}\Bigr)\le 1-\exp\Bigl\{-4\frac{ur}k\Bigr\}.
\end{equation}
The Reflection Principle shows that $W_\star^{(j)}$ has a Gaussian lower tail and our choice of~$\gamma_\star$ then ensures that the product on the left is at least $(1+\epsilon/2)^{-1}$ when~$a$ is sufficiently large, uniformly in~$k\ge1$. The inequality $1-\texte^{-s}\le s$ bounds the right-hand side by $4\frac{ur}k$. Using this in \eqref{E:4.37} proves the bound on the right of \eqref{E:4.32}. 
\end{proofsect}

If we do not insist on precise asymptotic and just aim for an upper bound, we can drop most of the restrictions between the parameters in Proposition~\ref{prop-6.1}:

\begin{corollary}
\label{cor-6.3}
Assume~$X$ is as Proposition~\ref{prop-6.1}. For each $\sigma\in(0,1/10)$ and~$a\ge0$ there exists $c>0$ such that
\begin{equation}
\label{E:6.17i}
E\biggl(\,\prod_{j=1}^{k-1}1_{\{X_j\ge\gamma(j)\}}\,\bigg|\, X_0=r,\,X_k=u\biggr)\le c\,\frac{(1+u)(1+r)}k
\end{equation}
holds for all~$r, u\ge 0$, all~$k\ge 2$ and all~$\gamma\colon\{0,\dots,k\}\to\R$ satisfying \eqref{E:4.31}.
\end{corollary}

\begin{proofsect}{Proof}
Fix~$\sigma\in(0,1/10)$ and let~$a_1$ be the constant from Lemma~\ref{lemma-4.6b} with $\epsilon:=1$. Given~$a\ge0$ let~$\gamma$ satisfy \eqref{E:4.31}. Set $\tilde a:=\max\{a_1,a,1\}$, abbreviate $\wt\gamma(i):=\tilde a+(\min\{i,k-i\})^\sigma$ and note that $\gamma\ge -\wt\gamma$. Using that the expectation increases when~$r$ and~$u$ increase, the argument \twoeqref{E:4.37}{E:4.57} along with the inequality $1-\texte^{-s}\le s$ show
\begin{equation}
\label{E:6.18i}
\begin{aligned}
E\biggl(\,&\prod_{j=1}^{k-1}1_{\{X_j\ge\gamma(j)\}}\,\bigg|\, X_0=r,\,X_k=u\biggr)
\\
&\le E\biggl(\,\prod_{j=1}^{k-1}1_{\{X_j\ge -\wt\gamma(j)\}}\,\bigg|\, X_0=\tilde a^{1/\sigma}+r,\,X_k=\tilde a^{1/\sigma}+u\biggr)
\\
&\le \biggl(\epsilon+4\prod_{j=1}^k P^{0,0}\Bigl(W_\star^{(j)}\ge-\min\{\wt\gamma(j-1),\wt\gamma(j)\}\Bigr)^{-1}\biggr)\frac{(\tilde a^{1/\sigma}+r)(\tilde a^{1/\sigma}+u)}k.
\end{aligned}
\end{equation}
To get the claim, note that the product is bounded uniformly in~$k\ge2$ and that $\tilde a^{1/\sigma}\ge1$ implies $(\tilde a^{1/\sigma}+r)(\tilde a^{1/\sigma}+u)\le \tilde a^{2/\sigma}(1+r)(1+u)$.
\end{proofsect}

\subsection{Reduction to a barrier estimate for random walks}
We now move attention to the Branching Random Walk with step distribution \eqref{E:1-step} whose first~$n$ generations can alternatively be viewed as a $\T_n$-indexed Markov chain with steps distributed as~$\NN(0,1/2)$. Our next goal will be to show that, along a path from the root to a near-maximal leaf, the Markov chain stays above a barrier of the kind studied earlier, modulo a linear tilt of the whole path.

We start with some notation. Given a sample~$h$ of the BRW on~$\T_n$ with step distribution $\NN(0,1/2)$,  for each~$k=1,\dots,n$ let 
\begin{equation}
%\label{}
M_k':=\max\bigl\{h(x)\colon x\in\BbbL_k,\,\theta_k(x)\ge b^{-1}\bigr\}.
\end{equation}
Departing from our previous labeling convention, let~$x_0:=\zero, x_1,\dots,x_n:=\varrho$ be the path from~$\zero$ to the root on~$\T_n$. 
For $u\in\R$ and~$A\in\sigma(h(x_0),\dots,h(x_n))$, define the probability measure
\begin{equation}
%\label{}
P_{n,u}(A):=\wt P\bigl(A\,\big|\,\,h(x_0)=\wt m_n+u\bigr)
\end{equation}
and, for~$s\ge u$, another measure
\begin{equation}
\label{E:6.18}
Q_{n,s,u}(A):=\wt E\biggl(1_A\prod_{k=1}^{n} \wt P\bigl(M_{k}'\le \wt m_{k}+t\bigr)\big|_{t:=\wt m_n+s-\wt m_{k}-h(x_k)}\,\bigg|\,h(x_0)=\wt m_n+u\biggr).
\end{equation}
In both cases the conditioning is well defined for all~$u\in\R$ thanks to the  probability density of~$h$ being a continuous function.

The measure~$P_{n,u}$ gives us access to the law of the sequence $(h(x_0),\dots,h(x_n))$ conditioned on~$h(\zero)=\wt m_n+u$, which by the structure of the Branching Random Walk reduces to the law of a random walk with step distribution~$\NN(0,1/2)$ conditioned to reach~$m_n+u$ at time~$n$. For~$Q_{n,s,u}$ we in turn get:

\begin{lemma}
\label{lemma-6.4a}
For each~$s\ge u$ and~$A\in\sigma(h(x_0),\dots,h(x_n))$,
\begin{equation}
%\label{}
Q_{n,s,u}(A) = \wt P\Bigl(A\cap\bigl\{\,\max_{x\in\BbbL_n}h(x)\le \wt m_n+s\bigr\}\,\Big|\,h(\zero)=\wt m_n+u\Bigr).
\end{equation}
\end{lemma}

\begin{proofsect}{Proof}
Given a sample~$h$ of the BRW on~$\T_n$, for each~$k=1,\dots,n$ set
\begin{equation}
%\label{}
M_k'':=\max\Bigl\{h(x)-h(x_k)\colon x\in\BbbL_n,\,b^{-n+k-1}\le\theta_n(x)<b^{-n+k}\Bigr\}.
\end{equation}
By the structure of~$\T_n$ and the Markov property of the BRW, the random variables $M_1'',\dots,M_n''$ are independent of each other and of $(h(x_0),\dots,h(x_n))$ with~$M_k''\laweq M_k'$ for each~$k=1,\dots,n$. Using these variables, the event $\{\max_{x\in\BbbL_n}h(x)\le \wt m_n+s\}$ becomes
\begin{equation}
%\label{}
\{h(x_0)\le \wt m_n+s\bigr\}\cap\bigcap_{k=1}^{n}\bigl\{h(x_k)+M_{k}''\le \wt m_n+s\bigr\}.
\end{equation}
For~$s\ge u$ the first event occurs automatically under the conditioning on $h(\zero)=m_n+u$. Using the independence of $M_1'',\dots,M_n''$, the conditional probability of the second event given $(h(x_0),\dots,h(x_n))$ turns into the product in \eqref{E:6.18}.
\end{proofsect}

We now proceed to prove three lemmas about $Q_{n,s,u}$. In the first lemma we observe that the product of the probabilities ``inside'' $Q_{n,s,u}$ effectively pushes the path $(h(x_0),\dots,h(x_n))$ below a barrier:

\begin{lemma}
\label{lemma-6.6}
For all~$\lambda>0$, $\sigma\in(0,1/10)$ and~$\delta\in(0,1)$ there exist~$c>0$ and~$a_2\ge0$ such that for all $s\ge 1$, all~$u$ satisfying $s\ge u\ge s-\lambda$ and all~$n$ with~$s\le n^{1-\delta}$,
\begin{equation}
\label{E:6.31a}
\begin{aligned}
Q_{n,s,u}\Biggl(\biggl(\,\bigcap_{k=1}^n\Bigl\{h(x_k)\le&\frac {n-k}n\wt m_n+s+\gamma_{n,a}(k)\Bigr\}\biggr)^\cc\,\Biggr)
\\
&\le
c\frac{(1+s-u)s}n\sum_{1\le k\le n/2}\bigl[k+\gamma_{n,a}(k)^2\bigr]\texte^{-\frac18 c_2\gamma_{n,a}(k)}
\end{aligned}
\end{equation}
holds for all~$a\ge a_2$ with
\begin{equation}
\label{E:6.34ii}
\gamma_{n,a}(k):=a+\bigl(\min\{k,n-k\}\bigr)^\sigma.
\end{equation}
Here~$c_2$ is as in \eqref{E:4.73}.
\end{lemma}

\begin{proofsect}{Proof}
Abbreviate
\begin{equation}
\label{E:6.34i}
A_j:=\biggl\{h(x_j)\le\frac {n-j}n\wt m_n+s+\gamma_{n,a}(j)\biggr\}.
\end{equation}
On~$A_k^\cc$ we have
\begin{equation}
\label{E:6.27}
\begin{aligned}
\wt m_n+s-\wt m_{k}-h(x_k)
&\le \wt m_n-\wt m_{k}-\frac {n-k}n\wt m_n -\gamma_{n,a}(k)
\\
&= \frac3{4\slb}\Bigl(\log k-\frac {k}n\log n\Bigr)-\gamma_{n,a}(k)\le-\frac12\gamma_{n,a}(k)
\end{aligned}
\end{equation}
once~$a$ is sufficiently large. Using \eqref{E:4.73} along with the fact that, by the FKG inequality and the symmetries of the tree, $\wt P(M_k'\le \wt m_k-\lambda)^2\le \wt P(M_k\le \wt m_k-\lambda)$, we then get
\begin{equation}
\label{E:6.35}
Q_{n,s,u}\biggl( (A_k\cap A_{n-k})^\cc \cap  \bigcap_{j=k+1}^{n-k-1}A_j\biggr)
\le 2\sqrt{c_1}\,\texte^{-\frac14 c_2\gamma_{n,a}(k)}P_{n,u}\biggl(\,\bigcap_{j=k+1}^{n-k-1}A_j\biggr).
\end{equation}
Since the probability in \eqref{E:6.31a} is dominated by the sum of the probabilities on the left of \eqref{E:6.35} and $\gamma_{n,a}(k)-\gamma_{n,a}(k-1)$ is uniformly bounded, it suffices to show
\begin{equation}
\label{E:6.38i}
P_{n,u}\biggl(\,\bigcap_{j=k}^{n-k}A_j\biggr)\le
c \frac{(1+s-u)s}n\bigl[k+\gamma_{n,a}(k)^2\bigr]\texte^{\frac18 c_2\gamma_{n,a}(k)}
\end{equation}
for some constant~$c>0$, uniformly in $1\le k\le n/2$ and~$u$ and~$s$ as above. 

As the quantity on the right of \eqref{E:6.38i} exceeds~$1$ once $n^\delta<k\le n/2$ and~$n$ is large, it suffices to focus on~$k\le n^\delta$. Abbreviate
\begin{equation}
\label{E:sub}
X_k:=\frac {n-k}n\wt m_n+s-h(x_k)
\end{equation}
 and note that $h(x_0)=\wt m_n+u$ translates into~$X_0 = s-u$ and~$h(x_n)=0$ into $X_n=s$. Using this we get
\begin{equation}
\label{E:6.31b}
P_{n,u}\biggl(\,\bigcap_{j=k}^{n-k}A_j\biggr)
=P\biggl(\,\bigcap_{j=k}^{n-k}\bigl\{X_j\ge-\gamma_{n,a}(j)\bigr\}\,\bigg| X_0=s-u,\,X_n=s\biggr),
\end{equation}
where, capitalizing on the fact that conditioning i.i.d.\ normals on their sum makes their mean irrelevant, ~$X$ is a random walk with step distribution~$\NN(0,1/2)$ under the law on the right. Shifting the whole path~$X$ by $\gamma_{n,a}(k)$ and abbreviating $\gamma(i):=\gamma_{n,a}(k)-\gamma_{n,a}(i)$ dominates the probability in \eqref{E:6.31b} by
\begin{equation}
\label{E:6.39i}
P\biggl(\,\bigcap_{j=k}^{n-k}\bigl\{X_j\ge\gamma(j)\bigr\}\,\bigg|\, X_0=s-u+\gamma_{n,a}(k),\,X_n=s+\gamma_{n,a}(k)\biggr).
\end{equation}
We now invoke $(\alpha+\beta)^\sigma\le\alpha^\sigma+\beta^\sigma$ to get $|\gamma(k+j)|\le j^\sigma$ for all $j=0,\dots,n/2-k$. This allows us to call upon Corollary~\ref{cor-6.3} under the conditional measure given $X_{k}$ and~$X_{n-k}$ to bound the probability \eqref{E:6.39i} by a constant times
\begin{equation}
%\label{}
E\biggl(\frac{|1+X_{k}||1+X_{n-k}|}{n-2k}\bigg|\, X_0=s-u+\gamma_{n,a}(k),\,X_n=s+\gamma_{n,a}(k)\biggr).
\end{equation}
Shifting $X$ by~$k\mapsto\frac kns+\frac{n-k}n(s-u)+\gamma_{n,a}(k)$ while using the Gaussian nature of~$X$, the expectation is written as $(n-2k)^{-1}$ times
\begin{equation}
\label{E:6.43i}
E\Bigl(\bigl|X_{k}+\varphi_n(k)\bigr|\,\bigl|X_{n-k}+\varphi_n(n-k)\bigr|\Big|\, X_0=0,\,X_n=0\Bigr),
\end{equation}
where~$\varphi_n(k):=1+\gamma_{n,a}(k)+\frac {k}ns+\frac{n-k}n(s-u)$. 

In order to bound \eqref{E:6.43i}, we first separate terms using the Cauchy-Schwarz inequality and the inequality $(\alpha+\beta)^2\le 2\alpha^2+2\beta^2$. Then we invoke the observation that 
\begin{equation}
%\label{}
E\bigl(X_k^2\,\big|\,X_0=0, X_n=0\bigr)\le c'\min\{k,n-k\},
\end{equation}
holds with some constant~$c'>0$ for all $k=1,\dots,n-1$ and note that
\begin{equation}
%\label{}
\varphi_n(k)\le\gamma_{n,a}(k)+s+1
\end{equation}
and, relying on $k\le n^{\delta}$ and $u\le s\le n^{1-\delta}$, also
\begin{equation}
%\label{}
 \varphi_n(k) \le\gamma_{n,a}(k)+2+s-u.
\end{equation}
Putting these together bounds \eqref{E:6.43i} by a constant times $(1+s-u)s[k+\gamma_{n,a}(k)^2]$ whenever $k\le n^\delta$, thus proving \eqref{E:6.38i} in this case.
\end{proofsect}

We now invoke the same argument as in the proof of Lemma~\ref{lemma-4.6b} to get:

\begin{lemma}
\label{lemma-6.7}
There exists~$\tilde c>0$ such that for all $n\ge1$, all $\gamma\colon\{0,\dots,n\}\to[0,\infty)$ and all $s\ge u\ge0$,
\begin{equation}
\label{E:6.31}
\begin{aligned}
Q_{n,s,u}\Biggl(\biggl(\,\bigcap_{j=k}^{n-k'}\Bigl\{h(x_j)\le&\frac {n-j}n\wt m_n+s-\gamma(j)\Bigr\}\biggr)^\cc\,\cap\,\bigcap_{i=1}^n\Bigl\{h(x_i)\le\frac {n-i}n\wt m_n+s+\gamma(i)\Bigr\}\Biggr)
\\
&\le
\tilde c\frac{(1+s-u)(1+s)}n\sum_{\ell=k}^{n-k'}\frac{n^{3/2}}{[\ell(n-\ell)]^{3/2}}\,\gamma(\ell)^5
\end{aligned}
\end{equation}
for any~$k,k'=1,\dots,\lfloor n/2\rfloor$ such that $\gamma(i)\ge1$ for all~$i=k,\dots,n-k'$.
\end{lemma}

\begin{proofsect}{Proof}
Abbreviate the second intersection in \eqref{E:6.31} as~$A$ while noting that the $i=n$ term can be dropped due to the fact that $s\ge0$. Denote
\begin{equation}
\label{E:6.44i}
B_j:=\Bigl\{h(x_j)\le\frac {n-j}n\wt m_n+s-\gamma(j)\Bigr\}.
\end{equation}
Relying on the substitution \eqref{E:sub} with~$\gamma$ instead of~$\gamma_{n,a}$ we then get
\begin{equation}
%\label{}
\begin{aligned}
Q_{n,s,u}(B_\ell^\cc\cap A)&\le P_{n,u}(B_\ell^\cc\cap A)
\\
&=E\biggl(1_{\{X_\ell< \gamma(\ell) \}}\prod_{j=1}^{n-1}1_{\{X_j\ge -\gamma(j)\}}\,\bigg|\, X_0=s-u,\,X_n=s\biggr).
\end{aligned}
\end{equation}
Proceeding as in \eqref{E:4.48a}, this is bounded by a term proportional to
\begin{equation}
%\label{}
\gamma(\ell)^3 [s-u+\gamma(\ell)][s+\gamma(\ell)]  \frac{\sqrt n}{[\ell(n-\ell)]^{3/2}}.
\end{equation}
Using the bounds $s-u+\gamma(\ell)\le\gamma(\ell)(1+s-u)$
and $s+\gamma(\ell)\le\gamma(\ell)(1+s)$
once~$\gamma(\ell)\ge1$ instead of \eqref{E:6.9i}, the claim follows by summing the bound over $\ell=k,\dots,  n-k'$.
\end{proofsect}

The conclusion of Lemma~\ref{lemma-6.7} will be used both in the proof of  Lemma~\ref{lemma-4.7} and Proposition~\ref{prop-3.3}. For Lemma~\ref{lemma-4.7} we will also need:

\begin{lemma}
\label{lemma-6.8}
Let~$\sigma\in(0,1/10)$ and let~$\gamma_{n,a}$ be as in \eqref{E:6.34ii}. For each~$a\ge0$ there exists a constant~$c'>0$ such that for all $n,k\ge1$ with $2k\le n$, all $s\ge 0$, all $u\le s$ and all $v>0$,
\begin{equation}
%\label{E:6.31}
\begin{aligned}
Q_{n,s,u}\Biggl(\Bigl\{h(x_k)<&\frac {n-k}n\wt m_n+s-v\Bigr\}\,\cap\,\bigcap_{i=1}^n\Bigl\{h(x_i)\le\frac {n-i}n\wt m_n+s+\gamma_{n,a}(i)\Bigr\}\Biggr)
\\
&\le
c'\frac{(1+s-u)(1+s)}nE\biggl((1+X_k)^2 1_{\{X_k>v\}}\,\bigg|X_0=s-u,\,X_n=s\biggr).
\end{aligned}
\end{equation}
\end{lemma}

\begin{proofsect}{Proof}
Proceeding via the substitution \eqref{E:sub}, the probability is dominated by
\begin{equation}
\label{E:6.43}
E\biggl(1_{\{X_k>v\}}\prod_{j=1}^{n-1}1_{\{X_j>-\gamma_{n,a}(j)\}}\,\bigg|\, X_0=s-u,\,X_n=s\biggr).
\end{equation}
We will bound this by conditioning on~$X_k$ and separately estimating the conditional probability of the product for indices less than~$k$ and larger than~$k$. 

Noting that~$\gamma_{n,a}(j)\le\gamma_{n,a}(k)$ for all~$j\le k$, the homogeneous Ballot Theorem (or Corollary~\ref{cor-6.3}) dominates the former probability as
\begin{equation}
%\label{}
\begin{aligned}
E\biggl(\,\prod_{j=1}^{k-1}1_{\{X_j>-\gamma_{n,a}(j)\}}\,\bigg|&\, X_0=s-u,\,X_k=\cdot\biggr)
\\
&\le c\frac{(1+s-u+\gamma_{n,a}(k))(1+X_k+\gamma_{n,a}(k))}k.
\end{aligned}
\end{equation}
For the latter probability we note that $\gamma_{n,a}(j)\le\gamma_{n-k,a}(j-k)+\gamma_{n,a}(k)$ which then gives
\begin{equation}
%\label{}
E\biggl(\,\prod_{j=k+1}^{n-1}1_{\{X_j>-\gamma_{n,a}(j)\}}\,\bigg|\, X_n=s,\,X_k=\cdot\biggr)
\le c\frac{(1+s+\gamma_{n,a}(k))(1+X_k+\gamma_{n,a}(k))}{n-k}.
\end{equation}
Separating terms as before, this dominates \eqref{E:6.43} by $(1+s-u)(1+s)$ times
\begin{equation}
%\label{}
c^2\frac{[1+\gamma_{n,a}(k)]^4}{k(n-k)}E\biggl((1+X_k)^2 1_{\{X_k>v\}}\,\big|X_0=s-u,\,X_n=s\biggr).
\end{equation}
Since $\sigma<1/4$, the prefactor is bounded by $c'/n$ independently of~$k\le n/2$.
\end{proofsect}

\subsection{Proofs of Proposition~\ref{prop-3.3} and Lemma~\ref{lemma-4.7}}
For the proofs of our claims, we first summarize the above lemmas as: 

\begin{lemma}
\label{lemma-6.8a}
Given any~$\sigma\in(0,1/10)$, $a\ge0$ and with $\gamma_{n,a}$ as in \eqref{E:6.34ii}, for $k\le n$ let
\begin{equation}
%\label{}
G_{n,k}:=\bigl\{X_k\in[k^{\sigma},k^{1-\sigma}]\bigr\}\cap\bigcap_{j=k}^{n-k}\bigl\{X_j\ge\gamma_{n,a}(j)\bigr\},
\end{equation}
where~$X$ is related to~$h$ via \eqref{E:sub}. For each~$\epsilon>0$ and $\lambda>0$ there exist~$a\in\R$ and $k_0\ge1$ such that for all~$s\ge 1$ and all $u\le s$ with $s-u\le\lambda$,
\begin{equation}
%\label{}
Q_{n,s,u}(G_{n,k}^\cc)\le\epsilon\frac {1+u_+} n
\end{equation}
holds once $n^\sigma\ge k\ge k_0$ and $n^{1-\sigma}\ge s$.
\end{lemma}

\begin{proofsect}{Proof}
Let~$A_j$ be defined by \eqref{E:6.34i} and set $A:=\bigcap_{j=1}^n A_j$. Next let~$B_j$ be as in \eqref{E:6.44i} with~$\gamma:=\gamma_{n,a}$ and set $\wt B_k:=\bigcap_{j=k}^{n-k} B_j$. Finally, set
\begin{equation}
%\label{}
C_k:=\Bigl\{h(x_k)\ge\frac{n-k}n\wt m_n+s-k^{1-\sigma}\Bigr\}.
\end{equation}
Then 
\begin{equation}
%\label{}
 G_{n, k}^\cc  =(\wt B_k\cap C_k)^\cc\subseteq A^\cc\cup (A\cap\wt B_k^\cc)\cup 
 (A \cap C_k^\cc).
\end{equation}
Invoking the union bound, it suffices to show that the $Q_{n,s,u}$-measure of each event on the right is less than~$\frac13\frac{1+u_+}{n} \epsilon$ once~$k$ is large and the other restrictions hold.

As $s\le u+\lambda\le (1+\lambda)(1+u_+)$, Lemma~\ref{lemma-6.6} with~$\delta:=\sigma$ gives $Q_{n,s,u}(A^\cc)\le\frac13\epsilon\frac{1+u_+}n $ once~$a$ is sufficiently large, uniformly in~$1\le k\le n/2$. With this~$a$ fixed, Lemma~\ref{lemma-6.7} does the same for $Q_{n,s,u}(A\cap\wt B_k^\cc)$ once~$k$ is sufficiently large (and~$2k<n$). Finally, Lemma~\ref{lemma-6.8} with $v:=k^{1-\sigma}$  bounds $Q_{n,s,u}( A \cap C_k^\cc )$  by a constant times $\frac {1+u_+}n$ times
\begin{equation}
\label{E:6.52}
E\biggl((1+X_k)^2 1_{\{X_k>k^{1-\sigma}\}}\,\bigg|\,X_0=s-u,\,X_n=s\biggr).
\end{equation}
A calculation (or an inspection of \eqref{E:6.6w}) shows that the probability density~$f$ of~$X_k$ under the conditional measure equals
\begin{equation}
\label{E:6.53i}
f(x) = \frac1{\sqrt\pi}\sqrt{\frac n{k(n-k)}}\exp\biggl\{- \frac{n}{k(n-k)} \Bigl(x-(s-u)-\frac kn u\Bigr)^2\biggr\}.
\end{equation}
Invoking the constraints $s-\lambda\le u\le s\le n^{1-\sigma}$ and $k\le n^\sigma$ (and~$2k<n$), this is bounded by a constant times the probability density of~$\NN(1+\lambda,k)$ whenever~$x\ge1+\lambda$. Since $1-\sigma>1/2$, it follows that the expectation \eqref{E:6.52} can be made as small as desired by taking~$k$ large, uniformly in~$n$ subject to $k\le n^\sigma$.
\end{proofsect}

We are now in a position to complete the proof of our first desired claim:

\begin{proofsect}{Proof of Lemma~\ref{lemma-4.7}}
Invoking the union bound and the symmetries of the tree, the desired probability is for any~$s\ge-\lambda$ bounded~by
\begin{equation}
\label{E:6.47i}
\begin{aligned}
\wt P\Bigl(\,\max_{x\in\BbbL_n}&\,h(x)>\wt m_n+s\Bigr)
\\
&+b^n\wt P\biggl(\,\bigl\{h(\zero)\ge \wt m_n-\lambda\bigr\}\cap\wt\BB_{n, k}(\zero)^\cc\cap\Bigl\{\max_{x\in\BbbL_n}h(x)\le\wt m_n+s\Bigr\}\biggr).
\end{aligned}
\end{equation}
Given any $\epsilon>0$, by \eqref{E:4.73} the first probability can be made smaller than~$\epsilon$ by taking~$s$ sufficiently large. (We assume $s\ge1$ in what follows.) Writing~$f_n$ for the probability density of~$h(\zero) - \widetilde{m}_n$, Lemma~\ref{lemma-6.4a} casts the second probability as
\begin{equation}
\label{E:6.47}
\int_{-\lambda}^s 
Q_{n,s,u}\biggl(\,\Bigl\{\wt m_n-\wt m_k-h(x_k)\in[k^{\tilde\sigma},k^{1-\tilde\sigma}]\Bigr\}^\cc\biggr)\,f_n(u)\,\textd u,
\end{equation}
where~$\tilde\sigma:=\frac1{12}$. 
As is checked by a calculation, there exists a constant~$c>0$ such that, uniformly in~$n\ge1$ over the interval of integration,
\begin{equation}
\label{E:6.49}
f_n(u)\le c n b^{-n}\texte^{-2u\slb}.
\end{equation}
We set $\epsilon':=\epsilon(2\slb)\texte^{-2\lambda\slb}c^{-1}$ in what follows. 
 
Next let~$\sigma:=\frac1{11}$ and use Lemma~\ref{lemma-6.8a} to find~$a\ge0$ and $k_0\ge1$ such that
\begin{equation}
%\label{}
Q_{n,s,u}(G_{n,k}^\cc)\le \epsilon'n^{-1}
\end{equation}
 for the given~$s$, all $u\in[-\lambda,s]$ and all $n>2k$ obeying $n^\sigma\ge k\ge k_0$ and~$n^{1-\sigma}\ge s$. As
\begin{equation}
%\label{}
\Bigl|\wt m_n -\wt m_k-\frac{n-k}n\wt m_n\Bigr|\le k^{\tilde\sigma}
\end{equation}
once~$k$ is large and $n>2k$, we get the inclusion
\begin{equation}
%\label{}
\Bigl\{\wt m_n-\wt m_k-h(x_k)\in[k^{\tilde\sigma},k^{1-\tilde\sigma}]\Bigr\}
\supseteq G_{n,k}
\end{equation}
once~$k$ is so large that also $2k^{\tilde\sigma}+s\le k^\sigma$ and 
$k^{1-\tilde\sigma}-k^{\tilde\sigma} +s \ge k^{1-\sigma}$ hold. Hence the probability under the integral is also at most $\epsilon' n^{-1}$ and the integral is thus bounded by
\begin{equation}
%\label{}
\epsilon'n^{-1}c n b^{-n}(2\slb)^{-1}\texte^{2\lambda\slb} =\epsilon b^{-n}.
\end{equation}
Plugging this in \eqref{E:6.47i}, the desired probability is less than~$2\epsilon$ once~$k$ is large and $n\gg k$. Since~$\epsilon$ is arbitrary, this implies the claim.
\end{proofsect}

The proof of Proposition~\ref{prop-3.3} requires additional lemmas. We start with control of the the total mass of~$Q_{n,u,u}$:

\begin{lemma}
\label{lemma-6.9}
For all~$\delta\in(0,1)$ there exist~$c>c'>0$ and $u_0>0$ such that for all~$n\ge1$ and all~$u\ge u_0$ with~$u \le n^{1-\delta}$,
\begin{equation}
%\label{}
c'\frac un\le Q_{n,u,u}(\R^{n+1})\le c\frac un.
\end{equation}
\end{lemma}

\begin{proofsect}{Proof}
We start with the lower bound.
Given~$\sigma\in(0,\ffrac1{10})$ and~$a\ge0$, let
\begin{equation}
%\label{}
A':=\bigcap_{k=1}^n\Bigl\{h(x_k)\le\frac{n-k}n\wt m_n+u + a^{1/\sigma} -\gamma_{n,a}(k)\Bigr\}.
\end{equation}
For $h(x_0)=\wt m_n+u$, the substitution \eqref{E:sub} gives
\begin{equation}
\label{E:6.63}
h(x_0)-h(x_j) =\wt m_n+u-h(x_j) = X_j+\frac jn\wt m_n.
\end{equation}
On~$A'$ we can bound the product in the definition of~$Q_{n,u,u}$ by
\begin{equation}
\label{E:6.64a}
\tilde c_n(a):=\prod_{k=1}^n P(M_k'\le \wt m_k+t)|_{t:=\wt m_n-\wt m_k-\frac{n-k}n\wt m_n+\gamma_{n,a}(k)- a^{1/\sigma}}
\end{equation}
from below and hereby get
\begin{equation}
\label{E:6.64}
Q_{n,u,u}( A')\ge \tilde c_n(a) E\biggl(\,\prod_{k =1}^{n} 1_{\{X_k\ge \gamma_{n,a} (k) - a^{1/\sigma}\}}\,\bigg|\,X_0=0,\,X_n=u\biggr).
\end{equation}
Using $M_k\ge M_k'$, the bound \eqref{E:4.73} and that $\inf_{n\ge1}P(M_n\le \wt m_n+t)>0$ for all $t\in\R$ we check $\inf_{n\ge1}\tilde c_n(a)>0$ for each~$a\ge0$, so it remains to find a lower bound on the expectation in \eqref{E:6.64}. Here a shift of the whole path by~$a^{1/\sigma}$ turns the expectation into the form in \eqref{E:4.32} with~$(r,u)$ given by $(a^{1/\sigma},u+a^{1/\sigma})$. Proposition~\ref{prop-6.1} with~$\epsilon:=1/2$ and~$\delta$ equal to half of that above  bounds the expectation from below by $2a^{1/\sigma}(u+a^{1/\sigma}) n^{-1}$ provided that $a\ge a_0$ and $a^{1/\sigma}(u+a^{1/\sigma})\le n^{1-\delta/2}$. This gives the claim once~$u$ is sufficiently large with~$u\le n^{1-\delta}$. 

For the upper bound let~$A_j$ be defined by \eqref{E:6.34i} and set $A:=\bigcap_{j=1}^n A_j$.  Lemma~\ref{lemma-6.6} gives that $Q_{n,u,u}(A^\cc)$ is at most a constant times~$u n^{-1}$ once~$u$ is sufficiently large. Next we invoke $Q_{n,u,u}(A)\le P_{n,u}(A)$ and bound the right-hand side by the expectation in \eqref{E:6.64} albeit with $X_k\ge \gamma_{n,a} (k) - a^{1/\sigma}$ replaced by~$X_k\ge -\gamma_{n,a} (k)$.  Corollary~\ref{cor-6.3} then bounds the expectation by a constant times~$u n^{-1}$ once~$u$ is sufficiently large, proving that also~$Q_{n,u,u}(A)$ is at most a constant times $un^{-1}$. 
\end{proofsect}

Next we give a representation of expectations of functions that depend only on a few initial values from $h(x_0),\dots,h(x_n)$:

\begin{lemma}
\label{lemma-6.10a}
Given $k\ge \ell\ge 1$, $\sigma\in(0,1/10)$ and $f\in C_\cc(\R^{\ell})$ let
\begin{equation}
%\label{}
\Xi_k(f):=E\Biggl(f\Bigl(\{X_j+j\slb\}_{j=1}^{\ell}\Bigr)\biggl(\,\prod_{j=1}^{k} \wt P(M_{j}'\le t)|_{t:=X_j+j\slb}\biggr)\,1_{\{X_k\in[k^{\sigma},k^{1-\sigma}]\}}X_k\Biggr),
\end{equation}
where~$X$ is the random walk with step distribution $\NN(0,1/2)$. Let~$\delta\in(0, 1/2)$. Then $\epsilon_{n,k}(u,f)$ defined for each~$u>0$ and $n\ge k$ by
\begin{equation}
%\label{}
E_{Q_{n,u,u}}\Bigl(f\bigl(\{h(x_0)-h(x_j)\}_{j=1}^\ell\bigr)\Bigr)=\frac{4u}n\Xi_k(f)+\frac un\epsilon_{n,k}(u,f)
\end{equation}
obeys
\begin{equation}
\label{E:6.57}
\lim_{k\to\infty}\limsup_{n\to\infty}\sup_{n^{\delta}\le u\le n^{1-\delta}}\bigl|\epsilon_{n,k}(u,f)\bigr|=0.
\end{equation}
Moreover, writing $1$ for the constant function equal to~$1$, we have
\begin{equation}
\label{E:6.58}
0<\inf_{k\ge1}\Xi_k(1)\le\sup_{k\ge1}\Xi_k(1)<\infty.
\end{equation}
\end{lemma}

\begin{proofsect}{Proof}
Fix $\sigma\in(0,1/10)$ and let~$\hat\sigma\in(0,\sigma)$.  Let~$\epsilon\in(0,1/2)$ and let~$a\ge0$ and $k_0\ge1$ be as in Lemma~\ref{lemma-6.8a} for $s=u$ (and so, e.g., $\lambda:=1$). Define the event
\begin{equation}
%\label{}
\wh G_{n,k}:=\bigl\{X_k\in[k^{\sigma},k^{1-\sigma}]\bigr\}\cap\bigcap_{j=k}^{n-k}\bigl\{X_j\ge\wh\gamma_{n,a}(j)\bigr\},
\end{equation}
where~$\wh\gamma_{n,a}$ is defined using~$\hat\sigma$ instead of~$\sigma$. Note that $\wh\gamma_{n,a}\le\gamma_{n,a}$ gives $\wh G_{n,k}\supseteq G_{n,k}$ and so the bound in Lemma~\ref{lemma-6.8a} applies to~$\wh G_{n,k}$ as well. Denoting
\begin{equation}
%\label{}
F_{n,k}:=\bigcap_{j=1}^k\bigl\{ (1-\epsilon)u\le X_{n-j}\le (1+\epsilon)u\bigr\}
\end{equation}
a calculation based on the conditional probability density \eqref{E:6.53i} shows that
\begin{equation}
\label{E:6.72i}
Q_{n,u,u}(F_{n,k}^\cc)\le P_{n,u}(F_{n,k}^\cc)\le c k\texte^{-\frac14\epsilon^2 u^2/k}
\end{equation}
for some constant~$c>0$, uniformly in $k,n\ge1$ subject to~$k/n<\epsilon/2$.
For~$u\ge n^\delta$, this decays faster than polynomially in~$n$.

In light of these observations, it suffices to prove the claim for the expectation restricted to the event $\wh G_{n,k}\cap F_{n,k}$. Using  \eqref{E:6.63} this can be written as 
\begin{equation}
%\label{}
\begin{aligned}
&E_{Q_{n,u,u}}\Bigl(1_{\wh G_{n,k}\cap F_{n,k}}\,f\bigl(\{h(x_0)-h(x_j)\}_{j=1}^\ell\bigr)\Bigr)
\\
&\,=E\Biggl(f\bigl(\{X_j+\tfrac jn\wt m_n\}_{j=1}^\ell\bigr)\biggl(\,\prod_{j=1}^k P(M_j'\le t)|_{t:=X_j+\frac jn\wt m_n}\biggr) Y_{n,k}(X) 1_{\wh G_{n,k}\cap F_{n,k}}\,\Bigg|\, X_n=u\biggr),
\end{aligned}
\end{equation}
where for brevity we dropped the explicit conditioning on~$X_0=0$ and where
\begin{equation}
%\label{}
Y_{n,k}(X):=\prod_{j=k+1}^n P(M_j'\le m_j+t)|_{t:=X_j+\frac jn\wt m_n - m_j}.
\end{equation}
On $\wh G_{n,k} \cap F_{n,k}$ we have
\begin{equation}
%\label{}
\begin{aligned}
Y_{n,k}(X)\ge \biggl(\,\prod_{j=k+1}^{n-k} P(M_j'\le m_j+&t)|_{t:= \wh \gamma_{n,a}(j)  +\frac jn\wt m_n - m_j}\biggr)
\\
&\times\prod_{j=n-k}^{n}P(M_j'\le m_j+t)|_{t:=\frac12 u+\frac jn\wt m_n - m_j}.
\end{aligned}
\end{equation}
Since $\frac jn\wt m_n-m_j \ge - \frac12  \wh \gamma_{n,a}(j)  $ once~$j$ is sufficiently large (with~$j \le n-k$), while $u\ge  n^{\delta}$ implies $\frac12 u+\frac jn\wt m_n - m_j\ge\frac14 u$ for all~$j=n-k,\dots,n$ once~$k$ is large (with~$2k<n$), the bound \eqref{E:4.73} shows that $1\ge Y_{n,k}\ge 1-\epsilon$ on~$\wh G_{n,k} \cap F_{n,k}$ once~$k$ is large and~$2k<n$.

With $Y_{n,k}(X)$ effectively removed, conditioning on~$X_k$ separates the indicator of the event $\wh G_{n,k}\cap F_{n,k}$ from the remaining terms inside the expectation. Conditioning on~$X_{n-k}$ in turn separates the events $\wh G_{n,k}$ and~$F_{n,k}$ from each other. As $\hat\sigma<\sigma$ gives $X_k-\wh\gamma_{n,a}(k)\ge \frac12 k^{\sigma}$ for~$k$ large, on $F_{n,k}$ Proposition~\ref{prop-6.1} gives
\begin{equation}
%\label{}
 (1-\epsilon)^4 \frac{4 uX_k}{n} 1_{\{X_k\in[k^\sigma,k^{1-\sigma}]\}}
\le P\bigl(\wh G_{n,k}\,\big|\, \sigma(X_k,X_{n-k})\bigr)
\le (1+\epsilon)^3\frac{4 uX_k}{n}1_{\{X_k\in[k^\sigma,k^{1-\sigma}]\}},
\end{equation}
where we also assumed that $n-2k\ge n(1+\epsilon)^{-1}$ on the right-hand side
 and that $X_k - k^{\wh \sigma} \ge (1-\epsilon) X_k$
and $u - k^{\wh \sigma} \ge (1-\epsilon) u$ on the left-hand side. Since \eqref{E:6.72i} implies $|P_{n,u}(F_{n,k})-1|\le\epsilon$ once~$n\gg k\gg1$ and~$u\ge n^\delta$, the Intermediate Value Theorem allows us to summarize the above bounds as
\begin{equation}
\begin{aligned}
E_{Q_{n,u,u}}\Bigl(1_{\wh G_{n,k}\cap F_{n,k}}\,f\bigl(\{h(x_0)-h(x_j)\}_{j=1}^\ell\bigr)\Bigr)=(1+O(\epsilon))\frac{4u}n\Xi_{n,k}(f),
\end{aligned}
\end{equation}
where $\Xi_{n,k}(f)$ abbreviates
\begin{equation}
%\label{}
E\Biggl(f\bigl(\{X_j+\tfrac jn\wt m_n\}_{j=1}^\ell\bigr)\biggl(\,\prod_{j=1}^k P(M_j'\le t)|_{t:=X_j+\frac jn\wt m_n}\biggr)1_{\{X_k\in[k^\sigma,k^{1-\sigma}]\}} X_k \,\Bigg|\, X_n=u\biggr).
\end{equation}
The claim now reduces to showing that
\begin{equation}
\label{E:6.72}
\lim_{k\to\infty}\limsup_{n\to\infty}\sup_{n^{\delta}\le u\le n^{1-\delta}}\bigl|\,\Xi_{n,k}(f)-\Xi_k(f)\bigr|=0.
\end{equation}
Note that this and Lemma~\ref{lemma-6.9} then give the bounds \eqref{E:6.58}.

To get \eqref{E:6.72} note that, thanks to~$u\le n^{1-\delta}$, the probability density of~$X_k$ under the conditional measure converges to that of the unconditional variable as~$n\to\infty$. Since also $\frac1n\tilde m_n\to\slb$, the convergence \eqref{E:6.72} then follows from the continuity of~$f$ and $t\mapsto P(M_j'\le t)$ along with the Bounded Convergence Theorem.
\end{proofsect}

\begin{proofsect}{Proof of Proposition~\ref{prop-3.3}}
We will write~$n$ instead of~$k$ as it is more consistent with the notations throughout this section. Let~$f\in\Cloc^+$ depend only on the coordinates in the subtrees rooted in~$x_0,\dots,x_\ell$, for some~$\ell\ge0$. Note that the law of~$h$ on the leaves in the subtree rooted at~$x_j$ with~$j\ge1$ has a bounded and continuous probability density. It follows that there exists a unique continuous function $\tilde f_{v,s}\colon\R^{\ell}\to\R$ such that
\begin{equation}
\label{E:6.78}
\begin{aligned}
\tilde f_{v,s}\bigl(\{h(x_0)-&\,h(x_j)\}_{j=0}^\ell\bigr)
\\
&= \frac{\wt E\Bigl(\texte^{-f(v,s,h(\zero)-h(\zero\cdot))}1_{\{h(x_0)=\max_{x\in\BbbL_n}h(x)\}}\,\Big|\,\sigma\bigl(h(x_0),\dots,h(x_n)\bigr)\Bigr)}
{\wt E\Bigl(1_{\{h(x_0)=\max_{x\in\BbbL_n}h(x)\}}\,\Big|\,\sigma\bigl(h(x_0),\dots,h(x_n)\bigr)\Bigr)}
\end{aligned}
\end{equation}
holds for almost every sample of~$h$.
Hence we get
\begin{equation}
%\label{}
\texte^{-g_{n,u}(v,s)}=\wt E\Bigl(\,\tilde f_{v,s}\bigl(\{h(x_0)-h(x_j)\}_{j=0}^\ell\bigr)\Big|\, h(\zero) = \max_{y\in\BbbL_n}h(y)=m_n+u\Bigr).
\end{equation}
Invoking the definition of~$Q_{n,s,u}$ we then rewrite this further as
\begin{equation}
\label{E:6.80}
\texte^{-g_{n,u}(v,s)}=\frac{E_{Q_{n,\tilde u,\tilde u}}(\, \tilde f_{v,s} (\{h(x_0) - h(x_j)\}_{j=0}^{\ell}) )}{E_{Q_{n,\tilde u,\tilde u}}(1)},
\end{equation}
where $\tilde u:=u+\wt m_n-m_n$.

Let now~$\tilde f$ be a generic function in~$C_\cc(\R^\ell)$. Lemma~\ref{lemma-6.10a} permits us to write the ratio in \eqref{E:6.80} as
\begin{equation}
%\label{}
\frac{E_{Q_{n,\tilde u,\tilde u}}(\, \tilde f (\{h(x_0) - h(x_j) \}_{j=1}^{\ell}) )}{E_{Q_{n,\tilde u,\tilde u}}(1)} = \frac{4\Xi_k(\tilde f)+\epsilon_{n,k}(\tilde u,\tilde f)}
{4\Xi_k(1)+\epsilon_{n,k}(\tilde u,1)}.
\end{equation}
Noting that the left-hand side is at most  $\Vert\tilde f\Vert$, we can consider a subsequence of $n\to\infty$ along which it converges. By \eqref{E:6.57}, taking this limit on the right-hand side wipes out the $\epsilon_{n,k}$-terms. With the help of \eqref{E:6.58}, this allows us to take~$k\to\infty$ along another subsequence to get
\begin{equation}
\label{E:6.82}
\lim_{n\to\infty}
\frac{E_{Q_{n,\tilde u,\tilde u}}(\, \tilde f (\{h(x_0) - h(x_j) \}_{j=1}^{\ell}) )}{E_{Q_{n,\tilde u,\tilde u}}(1)}
=\lim_{k\to\infty}\frac{\Xi_k(\tilde f)}{\Xi_k(1)}
\end{equation}
along any subsequence on the left and any subsequence of on the right. It follows that  both limits exist and are equal. The limit is uniform in~$\tilde u\in[ n^{\delta},n^{1-\delta}]$.

In order to express the result as an integral with respect to a measure, note that 
\begin{equation}
%\label{}
\tilde f\mapsto \lim_{k\to\infty}\frac{\Xi_k(\tilde f)}{\Xi_k(1)}
\end{equation}
is a positive linear functional on bounded continuous functions~$\R^{\ell}\to\R$ with norm one. Moreover, Lemmas~\ref{lemma-6.7}-\ref{lemma-6.8} imply that, for~$\tilde f$  supported outside $[-a,a]^{\ell+1}$, the value of the functional on~$\tilde f$ can be made as small as desired by taking~$a$ sufficiently large. The Riesz Representation Theorem then gives existence of a unique Borel probability measure~$\mu_\ell$ on~$\R^{\ell}$ such that
\begin{equation}
%\label{}
\lim_{k\to\infty}\frac{\Xi_k(\tilde f)}{\Xi_k(1)} = \int\tilde f\textd\mu_\ell
\end{equation}
for all bounded continuous $\tilde f\colon\R^{\ell}\to\R$. Since the measures $\{\mu_\ell\}_{\ell\ge0}$ are consistent, the Kolmogorov Extension Theorem implies that they are restrictions of a unique probability measure~$\mu$ on~$\R^{\N}$ concentrated on zero in the first coordinate.

To obtain~$\nu$ from~$\mu$, we need to undo the step that led from~$f(s,v,\cdot)$ to~$\tilde f$. Consider the canopy tree in Fig.~\ref{fig1} and, given a sample $\{S_k\}_{k\ge0}$ from~$\mu$, where~$S_0=0$, sample the Branching Random Walk with step distribution~$\NN(0,1/2)$ starting from $h(x_k):=S_k$, conditional on the leaf values to be non-negative everywhere. Identifying the leaves of the canopy tree with~$\N$ the distribution of~$h$ on the leaves is then~$\nu$.

The last item to address is the uniformity of the limit in parameters~$v$ and~$s$. Here we use that, for~$f$ continuous with compact support, for each~$\epsilon>0$ there exist $m\ge1$ and pairs $(v_1,s_1),\dots,(v_m,s_m)\in[0,1]\times\R$ such that for all~$(v,s)\in[0,1]\times\R$,
\begin{equation}
%\label{}
 \max_{i=1,\dots,m} \,\bigl\Vert f(v,s,\cdot)-f(v_i,s_i,\cdot)\bigr\Vert<\epsilon
\end{equation}
where~$\Vert\cdot\Vert$ is the supremum norm on~$\R^\N$. An elementary estimate then shows
\begin{equation}
%\label{}
 \max_{i=1,\dots,m}\Vert \tilde f_{v,s}-\tilde f_{v_i,s_i}\Vert<2\epsilon.
\end{equation}
 The convergence \eqref{E:6.82} is thus uniform in $\tilde f\in\{\tilde f_{v,s}\colon v\in[0,1],s\in\R\}$ as desired.
\end{proofsect}

%\newpage
\section*{Acknowledgments}
\nopagebreak\nopagebreak\noindent
This project has been supported in part by JSPS KAKENHI Grant Numbers 22K13927 and the NSF awards DMS-1954343 and DMS-2348113.

\bibliographystyle{abbrv}

\end{document}